\documentclass[final,3p,times]{elsarticle}
\usepackage{lineno,hyperref}
\usepackage{amssymb}
\usepackage{subfigure}
\usepackage{amsmath}
\usepackage{amsthm}
\usepackage{xcolor}
\usepackage{mathrsfs}
\usepackage[utf8]{inputenc}
\usepackage{textcomp}
\usepackage[font=small,labelfont=bf]{caption}
\usepackage{graphicx}
\usepackage[version=4]{mhchem}
\usepackage{longtable,tabularx}
\usepackage{comment}
\usepackage{multirow}
\usepackage[section]{placeins}
\usepackage{stmaryrd}
\usepackage{bm}
\usepackage{float}

\modulolinenumbers[200]
\journal{Journal}

\begin{document}

\begin{frontmatter}

\title{Physics-informed neural networks for inverse problems in supersonic flows}


\author{Ameya D. Jagtap$^{1,\dagger}$, Zhiping Mao$^{2,\dagger}$, Nikolaus Adams$^3$, and George Em Karniadakis$^{1}$}
\cortext[mycorrespondingauthor]{ Corresponding author Email: george$\_$karniadakis@brown.edu \newline $\dagger$ First two authors contributed equally.}

\address{$^1$~Division of Applied Mathematics, Brown University, 182 George Street, Providence, RI, 02912, USA}
\address{$^2$~School of Mathematical Sciences Xiamen University, Xiamen, Fujian 361005, China.}
\address{$^3$~Department of Mechanical Engineering, Technical University of Munich, 85748 Garching, Germany.}

\begin{abstract}
Accurate solutions to inverse supersonic compressible flow problems are often required for designing specialized aerospace vehicles. In particular, we consider the problem where we have data available for density gradients from Schlieren photography  as well as data at the inflow and part of wall boundaries. These inverse problems are notoriously difficult and traditional methods may not be adequate to solve such ill-posed inverse problems. To this end, we employ the physics-informed neural networks (PINNs) and its extended version, extended PINNs (XPINNs), where domain decomposition allows to deploy locally powerful neural networks in each subdomain, which can provide additional expressivity in subdomains, where a complex solution is expected. Apart from the governing compressible Euler equations, we also enforce the entropy conditions in order to obtain viscosity solutions. Moreover, we enforce positivity
conditions on density and pressure.
We consider inverse problems involving two-dimensional expansion waves, two-dimensional oblique and bow shock waves. We compare solutions obtained by PINNs and XPINNs and invoke some theoretical results that can be used to decide on the generalization errors of the two methods.
\end{abstract}


\end{frontmatter}

\linenumbers

\section{Introduction}
In recent years physics-informed neural networks (PINNs) \cite{Raissi2019JCP} emerged as an alternative simple method to solve many problems in computational science and engineering, see, for example \cite{raissi2018deep, mao2020physics, raissi2020hidden, jagtap2020extended, jagtap2020conservative, Shukla2022A, shukla2021parallel, chen2020physics, bin2021pinneik, raissi2020hidden, haghighat2021physics, yang2019adversarial, cai2021flow, kharazmi2021hp, jagtap2022deep2}. In particular, PINNs do not require meshes and can efficiently solve forward problems and even ill-posed inverse problems, which are otherwise difficult or sometimes even impossible to solve using traditional numerical methods. Moreover, PINN can easily handle noisy, sparse and multi-fidelity data sets. The main advantage of the PINN methodology is that it can seamlessly incorporate all the given information like governing equation, experimental data, initial/boundary conditions into the loss function, thereby recasting the original problem into an equivalent optimization problem. Recently, Shin et al. \cite{shin2020convergence} established the mathematical foundation of PINNs for linear partial differential equations, whereas in \cite{mishra2021estimates}, Mishra and Molinaro presented estimate on the generalization error of the PINN methodology.

In this work, we consider inverse problems of the shock wave problems in supersonic compressible flows. The governing equations of such flows are compressible Euler equations, which admit discontinuities or shocks, even though the initial states are smooth. Such ill-posed inverse problems are difficult or even sometimes impossible to solve using the traditional numerical solvers.
Moreover, for the shock wave problem, the traditional numerical methods usually require boundary conditions (BCs) for all field variables. However, for the simulations of high-speed flows, traditional numerical methods are not only expensive, but also it is usually difficult to determine the BCs, for example, as mentioned in \cite{anderson1995computational}, it is particularly vital to couple the surface BCs into the flow-field calculation for the expansion wave problem. Far-field BCs or symmetric BCs were employed for the compressible flows ~\cite{bayliss1982far, anderson1995computational}, but applying these kinds of BCs results in large computational domains, and therefore becomes computationally expensive. Other BCs like numerical boundary condition may cause numerical instability \cite{hall1981implementation}. We investigate in the present work the possibility of using the idea of PINN to resolve this issue by extending the algorithm from the one-dimensional case \cite{mao2020physics} to the two-dimensional shock wave problems. Another advantage of PINNs is that PINNs can be employed in any arbitrarily shaped domain of physical interest, while the traditional numerical methods need a specified large computational domain. 
In the PINNs literature, Mao et al. \cite{mao2020physics} first solved both forward and inverse one and two-dimensional Euler equations involving shocks in the Cartesian domain. Recently, Patel et al. \cite{patel2022thermodynamically} introduced thermodynamically consistent physics-informed neural networks, which enforce the entropy conditions for the scalar conservation laws and the one-dimensional Euler equations. In \cite{fuks2020limitations} Fuks and Tchelepi proposed a stable way for PINNs to handle scalar conservation laws involving discontinuities.
Moreover, Jagtap et al. proposed the spatial domain decomposition based PINN approach for conservative laws called as conservative PINN (cPINN) \cite{jagtap2020conservative}, and more general space-time domain decomposition based PINN approach for any type of differential equations, named eXtended PINNs (XPINNs) \cite{jagtap2020extended}. Apart from parallel implementation \cite{shukla2021parallel}, such domain decomposition approaches can allow local powerful networks in subdomains where complex solutions can be expected. In \cite{hu2021extended} Hu et al. provide a prior generalization bound via the complexity of the target functions in the partial differential equation (PDE) problem, and a posterior generalization bound via the posterior matrix norms of the networks after optimization. Based on these bounds, one can analyze the conditions under which XPINNs improve generalization. The main feature of XPINN, namely the spatio-temporal domain decomposition, introduces a  generalization trade-off. The XPINNs decompose the complex solution of differential equations into several simple parts, which decreases the complexity needed to learn by the networks in the respective subdomains and thereby boosts generalization. In contrast, such decomposition leads to less training data being available in each subdomain, and hence such a model is typically vulnerable to over-fitting, due to which it may become less generalizable. In the present study we used bounds in \cite{hu2021extended} for the decomposition of computational domain.

Apart from PINNs, the other frameworks which employed deep neural networks for tackling shock wave problems can be found in \cite{magiera2020constraint, bezgin2021data, mao2021deepm, monfort2017deep, he2020inverting}.

\begin{figure}[http]
    \centering
    \subfigure[Attached waves: oblique shock and expansion waves]{\label{intro:oblique:expansion}
    \includegraphics[width=0.4\textwidth]{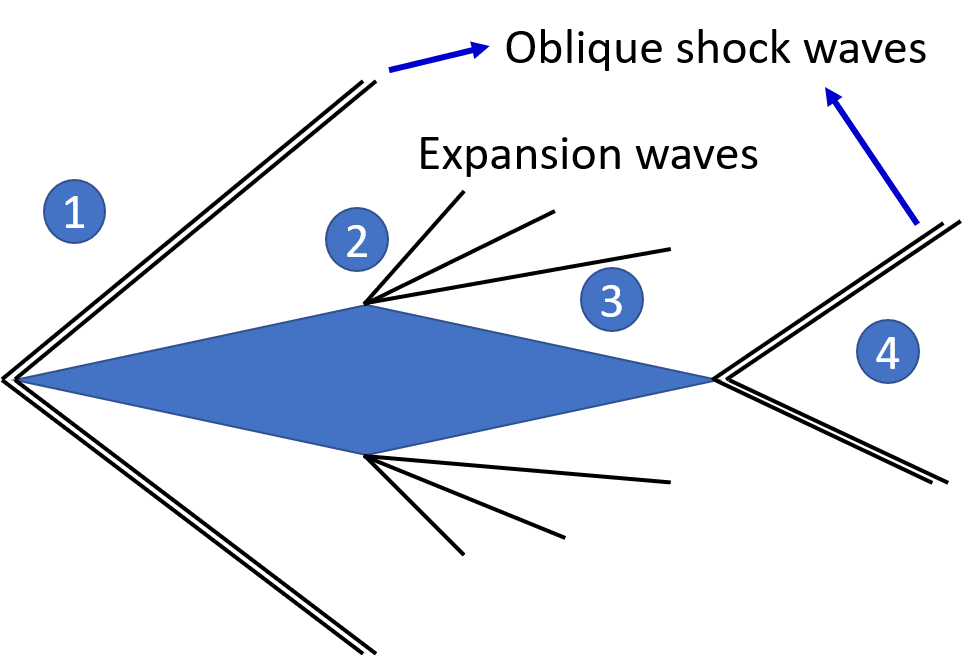}}
    ~~~~~~~~~~~~
    \subfigure[Detached waves: bow shock waves]{\label{intro:bow}
    \includegraphics[width=0.4\textwidth]{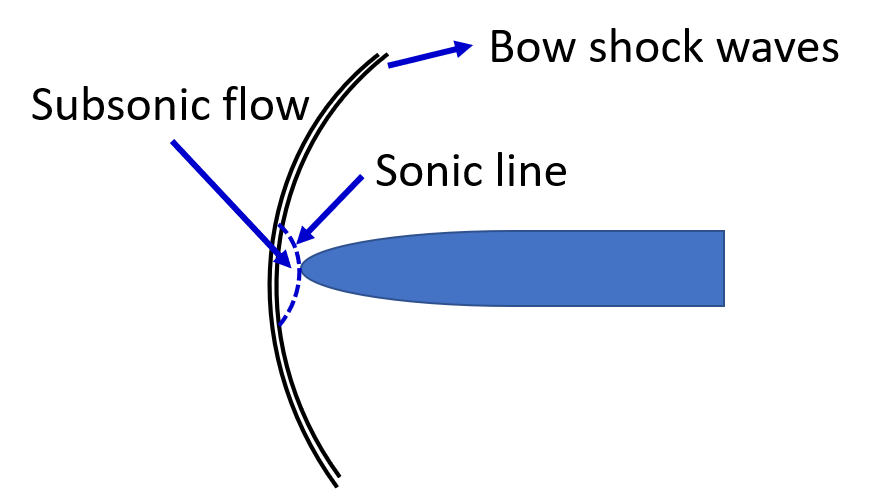}}
    \caption{Schematic of the three types of shock waves. Left: Oblique shock (regions \textcircled{1} and \textcircled{4}) and expansion waves (regions \textcircled{2} and \textcircled{3}), right: bow shock waves.}
    \label{fig:schematic:shock}
\end{figure}

Inverse problems in supersonic compressible flows are often encountered in designing specialized high-speed vehicles. Apart from the governing Euler equations, which  admit multiple weak solutions, we also enforce the entropy conditions that are satisfied by the viscosity solutions as well as positivity conditions on density and pressure fields.
In this work, we shall infer all the states (density, velocity and pressure) for the shock wave problems, including an oblique shock wave, expansion wave and bow shock wave (see Fig. \ref{fig:schematic:shock}), in the two-dimensional case. 
In particular, we use the idea of hidden fluid mechanics based on PINNs to solve steady state problems.  More precisely, we use the information of density gradients, a small data set on primitive variables, as well as global constraints to predict the density, velocity, and pressure for these inverse problems. In particular, we use both PINNs and XPINNs to infer unknown states involving shocks and expansion waves. 

This paper is arranged as follows: in Section 2 we present the governing compressible Euler equations and the entropy conditions. In Section 3 we describe the methodology, including the mathematical setup for the feed-forward neural networks, with a brief introduction to PINNs and XPINNs methodology, and the optimization methods. In section 4, various computational examples such as expansion wave, oblique shock wave and the bow shock wave are solved and various comparisons are made. Finally, we summarize our findings in Section 5.

\section{Governing Equations}
The governing compressible Euler equations in the conservative form~\cite{Courant1999, Liepmann2001, Zucker2002} can be written as:
\begin{equation}\label{NCL}
   \partial_t U + \nabla\cdot G(U) = 0, ~~~ \bm{x} \coloneqq \{x,y\} \in \Omega\subset \mathbb{R}^2, ~~ t\in (0, T], 
\end{equation}
where
$$ U =  \left[ {\begin{array}{c}
\rho\\
\rho u\\
\rho v\\
\rho E\\
\end{array} } \right],
\; G = \{G_1, G_2\}, \;
\text{with }
G_1(U) =\left[ {\begin{array}{c}
\rho u\\
p+ \rho u^2 \\
 \rho u v\\
 p u+ \rho u E\\
\end{array} } \right], 
\;
G_2(U) =\left[ {\begin{array}{c}
\rho v\\
\rho uv \\
p+ \rho v^2\\
 p v+ \rho v E\\
\end{array} } \right].
$$
Here, $\rho$ is the density, $p$ is the pressure, $(u,\,v) $ are the velocity components,  and $E$ is the total energy.
We use the additional equation of state, which describes the relation of the pressure and energy, to close the above Euler equations. In particular, we consider the equation of state for a polytropic gas given by
\begin{equation}\label{eq:EOS}
    p = (\gamma-1)\left(\rho E - \frac{1}{2}\rho \|{\bm u}\|^2\right),
\end{equation}
where $\gamma$ is the ratio of specific heats, whose value for air is 1.4, and $\bm{u} = (u, v)$.
Here, we solve the steady state Euler equations given as  
\begin{equation}\label{NCL:steady}
    \nabla\cdot G(U) = 0, \; \bm{x} \in \Omega.
\end{equation}

\textbf{Remark}: It is important to note that for neural network training we employ the non-dimensional form of Euler equations, which look exactly the same as the dimensional ones.

\subsection{Entropy conditions}
In general, equation \eqref{NCL} admits several weak solutions \cite{Dafermos2016} and thus additional conditions have to be imposed to select the `physically' relevant solution, among the others, the so-called \textit{Entropy condition}. The viscosity solutions \cite{bianchini2005vanishing} of conservation laws are associated with the zero-viscosity limit given by $U =  \lim_{\epsilon \rightarrow 0^+} U_{\epsilon} $, where $ \partial_t U_{{\epsilon}} + \nabla\cdot G(U_{\epsilon}) = \epsilon \nabla^2 U_{\epsilon}$.
These viscosity solutions also satisfy entropy conditions. 

A pair $(\eta, \phi)$ with a convex function $\eta(U)$ and an associated entropy flux (a vector valued function) $\boldsymbol{\phi}(U) = [\phi_1(U), \phi_2(U)]^T$ are called entropy-entropy flux pair \cite{godlewski1996numerical}, if it satisfies the relations
$$ \eta'(U)G_1'(U) = \phi'_1(U),  \ \ \ \ \ \ \eta'(U)G_2'(U) = \phi'_2(U). $$
Thus, an additional conservation law can be added to the governing compressible Euler equations
\begin{equation}
\eta_t + \phi_{1_x}+ \phi_{2_y} = 0
\end{equation}
for smooth solutions. At discontinuities, the entropy is not always conserved thus the admissible solution is chosen based on the following inequality.
\begin{equation}
\eta_t + \phi_{1_x}+ \phi_{2_y} \leq 0.
\end{equation}
All the admissible, bounded and non-smooth solutions satisfy the Entropy condition in a distributional sense.
It is important to note that the given hyperbolic conservation laws may not have a unique entropy pair \cite{godlewski1996numerical}. 
For Euler equations, the most common choice of an entropy-entropy flux pair is 
\begin{equation}
\eta \triangleq \frac{-\rho s}{\gamma-1} \ \ \text{and} \ \ \boldsymbol{\phi} \triangleq\frac{-\rho \bm{u} s}{\gamma-1}.
\end{equation}
where the specific entropy $s = \log(p/\rho^{\gamma})$.Note that during the network training, we also enforced the positivity conditions on density and pressure. This will avoid the problem of handling negative values for pressure and density. This can be done by choosing the following output from the neural networks \textit{max}$(\alpha, \rho)$ and \textit{max}$(\alpha, p)$ for density and pressure, respectively. Here $0<\alpha <<1$ is a small positive constant.

\section{Methodology}
\subsection{Mathematical setup for feed-forward neural network (FFNN)}

Let $\mathcal{N}^L: \mathbb{R}^{N_0} \rightarrow \mathbb{R}^{N_L}$ be a FFNN of $L$ layers and $N_k$ neurons in the $k^{th}$ layer, where the input layer has $N_0$ and the output layer $N_L$ neurons. The weights and biases in the $k^{th}$ layer  ($1 \leq k \leq L$) are given as $\bm{W}^k \in \mathbb{R}^{N_k \times N_{k-1}}$ and $\bm{b}^k \in \mathbb{R}^{N_k}$, respectively. If $\bm{z}\in \mathbb{R}^{N_0}$ is the input vector at the $k^{th}$ layer, then the output vector at $k^{th}$ layer is denoted by $\mathcal{N}^k(\bm{z})$. Thus, the input vector can be expressed as $\mathcal{N}^0(\bm{z}) = \bm{z}$. We employ the layer-wise adaptive activation function, denoted by $\sigma(n a^k)$, where $a^k$ are the trainable parameters, $n$ is the predefined scaling factor and $\sigma$ is the activation function,  see Jagtap et al. \cite{jagtap2020adaptive, jagtap2020locally}. The introduction of the additional activation slope parameters $a^k$ dynamically changes the slope of activation function during training, thereby increases the training speed of the neural networks, see  \cite{jagtap2020adaptive, jagtap2020locally}. Recently, Jagtap et al. \cite{jagtap2022deep} proposed the \textit{Rowdy activation function}  that is designed to eliminate any saturation region by injecting sinusoidal
noise like effects, which further increases the training speed. In this paper, we employ layer-wise locally adaptive activation functions with scaling factor $n = 10$ for all hidden-layers and initialize $n a^k = 1, ~\forall k$. The $(L-1)$-hidden layer FFNN is defined as
\begin{equation}
\begin{aligned}
& \mathcal{N}^1(\bm{z}) = \bm{W}^1 \bm{z} + \bm{b}^1, \quad \text{ for }\quad k=1,\\
& \mathcal{N}^k(\bm{z}) = \bm{W}^k \sigma(a^{k-1} \mathcal{N}^{k-1}(\bm{z})) + \bm{b}^k \in \mathbb{R}^{N_k},\quad \text{ for }\quad 2\leq k \leq L\ ,
\end{aligned}
\end{equation} where in the activation function is identity in the last layer.
Let $\tilde{\boldsymbol{\Theta}} = \{ \bm{W}^k, \bm{b}^k, a^k \}_{k=1}^L \in \mathcal{V}$ be the collection of all trainable parameters including weights, biases, and activation slopes, and $\mathcal{V}$ is the parameter space, then, the neural network output can be written as
$$\xi_{\tilde{\boldsymbol{\Theta}}}(\bm{z}) = \mathcal{N}^L(\bm{z}; \tilde{\boldsymbol{\Theta}})\ ,$$
where $\mathcal{N}^L(\bm{z}; \tilde{\boldsymbol{\Theta}})$ emphasizes the dependence of the neural network output $\mathcal{N}^L(\bm{z})$ on the trainable parameters $ \tilde{\boldsymbol{\Theta}}$. In general, weights and biases are initialized from known probability distributions. In particular, we used the Xavier initialization \cite{glorot2010understanding} for all the trainable parameters. 

\subsection{Overview of Physics-Informed Neural Networks (PINNs)}
\begin{figure}[!h]
\begin{center}
\includegraphics[scale=0.75,angle=0]{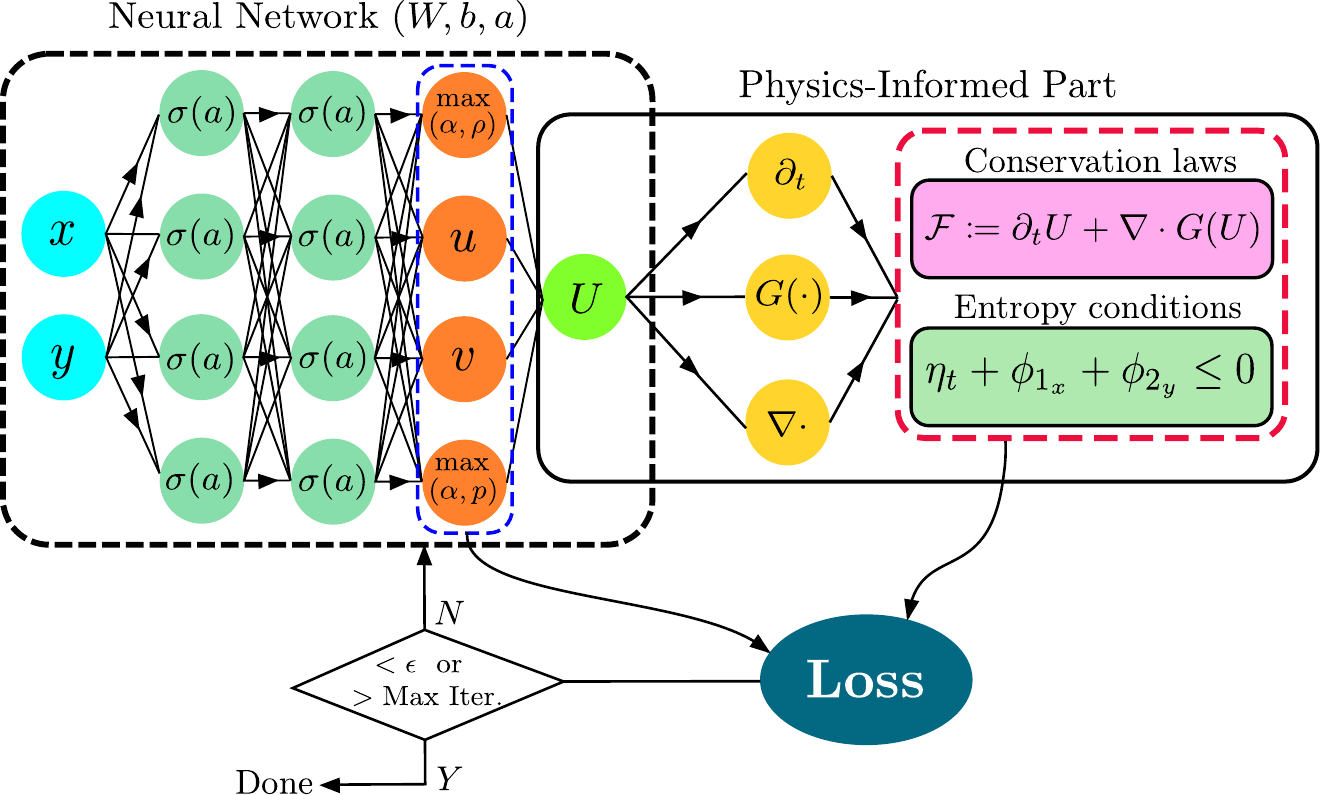}
\end{center}
\caption{
Schematic representation of the PINNs for the Euler equations. The left part is the 2 hidden-layer deep feed-forward neural networks with input $x$ and $y$, and with the primitive variables as the output. Here, $\sigma(\cdot)$ is the nonlinear activation function, whereas $a$ is the tunable activation slope. The right part is the physics-informed part, where the network output is forced to satisfy conservation laws as well as entropy conditions. Both parts contribute to loss function. Note that during the network training the positivity condition on the density and the pressure is enforced by choosing the maximum value between them and the small positive constant $\alpha$. }
\label{fig:Schematic}
\end{figure}

For a steady compressible inviscid flow governed by the Euler equations \eqref{NCL:steady}, we assume that apart from the inflow and wall boundary data on partial or full primitive variables, there is a set of experimental data (Schlieren data) available on density gradient. By using such data and combining them with the Euler equations \eqref{NCL:steady}, we would like to design PINNs to infer all the states of interest, i.e., the density $\rho(x,y)$, the velocity $\bm{u}(x,y)$ and the pressure $p(x,y)$.
Fig. \ref{fig:Schematic} shows the schematic representation of the PINNs for the Euler equations. The left part is the feed-forward neural networks with input $x$ and $y$ and the primitive variables as the output. The right part is the physics-informed part where  network output is forced to satisfy conservation laws. Both parts contribute to loss function.
In the present work, we shall learn the density, velocity, and pressure fields by using a feed-forward neural network. In particular, we shall approximate the density, velocity and pressure by setting the neural network prior on these primitive variables denoted jointly as $\xi_{\tilde{\boldsymbol{\Theta}}}(x,t)$.

Let $\{\bm{x}^{(i)}_{\xi}\}_{i=1}^{N_{\xi}}$ and $\{\bm{x}^{(i)}_F\}_{i=1}^{N_F}$ be the set of randomly selected training data and residual points, respectively. These points are drawn from a distribution, which is usually not known a priori and often needs to be chosen from the given input training data. The PINN algorithm aims to learn a surrogate $\xi \approx \xi_{\tilde{\boldsymbol{\Theta}}}$ for predicting the solution $\xi$ of the given PDE. The PINN loss function is given as
\begin{equation}\label{Loss}
 \mathcal{J}(\tilde{\boldsymbol{\Theta}})  =   W_{\xi} ~ [ \text{MSE}_{\xi}(\tilde{\boldsymbol{\Theta}}; \{\bm{x}^{(i)}_{\xi}\}_{i=1}^{N_{\xi}} ) + \text{MSE}_{\nabla \rho} (\tilde{\boldsymbol{\Theta}}; \{\bm{x}^{(i)}_{\nabla \rho}\}_{i=1}^{N_{\nabla \rho}} )] + W_{\mathcal{F}} ~ \text{MSE}_{\mathcal{F}}(\tilde{\boldsymbol{\Theta}}; \{\bm{x}^{(i)}_F\}_{i=1}^{N_F} ) \ ,
\end{equation}
 where $W_{\xi}$ and $W_{\mathcal{F}}$ are the weights for the data mismatch and residual terms, respectively. The mean squared errors (MSEs) are given by
\begin{align*}
\text{MSE}_{\xi}(\tilde{\boldsymbol{\Theta}}; \{\bm{x}^{(i)}_{\xi}\}_{i=1}^{N_{\xi}} ) & = \frac{1}{N_{\xi}} \sum_{i=1}^{N_{\xi}}\left|\xi^{(i)} - \xi_{{\tilde{\boldsymbol{\Theta}}}}(\bm{x}^{(i)}_{\xi})\right|^2\ , \ \\  \text{MSE}_{\mathcal{F}}(\tilde{\boldsymbol{\Theta}}; \{\bm{x}^{(i)}_F\}_{i=1}^{N_F} )  &= \frac{1}{N_{F}} \sum_{i=1}^{N_{F}}\left|\mathcal{F}_{{\tilde{\boldsymbol{\Theta}}}}(\bm{x}^{(i)}_{F})\right|^2,
\\ \text{MSE}_{\nabla \rho} (\tilde{\boldsymbol{\Theta}}; \{\bm{x}^{(i)}_{\nabla \rho}\}_{i=1}^{N_{\nabla \rho}} )&= \frac{1}{N_{\nabla}} \sum_{i=1}^{N_{\nabla}}| \nabla \rho^i - \nabla \rho_{{\tilde{\boldsymbol{\Theta}}}}(\bm{x}^i_{\nabla})|^2.
\end{align*}
 The term $\text{MSE}_{\xi}$ is the MSE for data mismatch term for primitive variables, which enforces the given experimental and/or synthetic data. $\text{MSE}_{\nabla \rho}$ represent the data mismatch term for density gradient. The term $\text{MSE}_{\mathcal{F}}$ is the MSE for the governing PDE residue  where $\mathcal{F}_{{\tilde{\boldsymbol{\Theta}}}}(\bm{x}_{F}) \coloneqq \mathcal{F}(\xi_{{\tilde{\boldsymbol{\Theta}}}}(\bm{x}_{F}))$ represents the residual of the governing system of conservation laws given by equation \eqref{NCL:steady} and the entropy conditions. The notations $N_{\xi}$, $N_F$ and $N_{\nabla}$ are the number of training data, residual points and number of data points for density gradient, respectively. The trainable parameters of the neural networks $\xi_{\tilde{\boldsymbol{\Theta}}}$ can be estimated by minimizing the loss function given by Eq. \eqref{Loss}.

The PDE residual $\mathcal{F}_{{\tilde{\boldsymbol{\Theta}}}}$ in the loss function is constructed using automatic differentiation \cite{baydin2018automatic}, which is a graph-based differentiation method for calculating derivatives of field variables accurately in a computational graph. Moreover, compared to traditional numerical differentiation, automatic differentiation does not suffer from truncation and round-off errors. Highly popular machine learning libraries such as TensorFlow \cite{abadi2016tensorflow} and Pytorch \cite{paszke2017automatic} provide automatic differentiation package.

\subsection{Extended Physics-Informed Neural Networks (XPINNs)}
The extended physics-informed neural networks (XPINNs) methodology \cite{jagtap2020extended} is a recently developed generalized space-time domain decomposition approach for deep learning of PDEs. It overcomes many limitations of the vanilla PINN method such as parallel implementation capacity, deployment of locally powerful networks, large representation capacity, etc. Here, we shall briefly discuss the XPINN method.

The computational domain $\Omega$ is divided into $N_{sd}$ number of arbitrary, non-overlapping subdomains.
Let $\{\bm{x}^{(i)}_{\xi_q}\}_{i=1}^{N_{\xi_q}} $, $\{\bm{x}^{(i)}_{F_q}\}_{i=1}^{N_{F_q}} $ and $\{\bm{x}^{(i)}_{I_q}\}_{i=1}^{N_{I_q}}$ be the set of randomly selected training, residual, and the common interface points, respectively in the $q^{th}$ subdomain.  The $N_{\xi_q}, N_{F_q}$, and $ N_{I_q}$ represent the number of training data points, the number of residual points, and the number of points on the common interface in the $q^{th}$ subdomain, respectively.
Similar to PINN, the XPINN algorithm aims to learn a surrogate $\xi_q \approx \xi_{\tilde{\boldsymbol{\Theta}}_q}$, $q = 1,2,\cdots, N_{sd}$ for predicting the solution $\xi \approx \xi_{\tilde{\boldsymbol{\Theta}}} = \sum_{q=1}^{N_{sd}} \xi_{\tilde{\boldsymbol{\Theta}}_q} \cdot 1_q$ of the given PDE over the entire computational domain. Here $1_q$ denotes the indicator function; see \cite{jagtap2020extended} for more details. The loss function of XPINN is defined subdomain-wise, which has a similar structure as the PINN loss function in each subdomain but is endowed with the interface conditions for stitching the subdomains together. For the forward problem, the loss function in the $q^{th}$ subdomain is defined as
\begin{align}\label{XPINNloss}
 \mathcal{J}(\tilde{\boldsymbol{\Theta}}_q)  =  ~& W_{\mathcal{F}_q} ~(\text{MSE}_G^{Mass}+\text{MSE}_G^{MomX}+\text{MSE}_G^{MomY}+\text{MSE}_G^{Ene}+\text{MSE}_G^{Entropy}) \nonumber 
  \\  +& W_{\xi_q} ~ (\text{MSE}_{\xi}^D  + MSE_{\nabla \rho}    ) 
 +  W_{I_q} ~ \text{MSE}_{\xi_{avg}} + W_{I_{\mathcal{F}_q}} ~  \text{MSE}_{\mathcal{R}}, ~~ q = 1,2,.\cdots, N_{sd}, 
\end{align}
where the MSE for mass, momentum ($x$ and $y$ directions), energy conservation laws, and entropy condition are given as
\begin{align*}
\text{MSE}_G^{Mass} & = \frac{1}{N_{F_q}} \sum_{i=1}^{N_{F_q}}| (G_1^{Mass}(\bm{x}^i_{F_q}))_{x} +(G_{2}^{Mass}(\bm{x}^i_{F_q}))_{y}|^2,
\\ \text{MSE}_G^{MomX} & = \frac{1}{N_{F_q}} \sum_{i=1}^{N_{F_q}}|(G_{1}^{MomX}(\bm{x}^i_{F_q}))_{x} +(G_{2}^{MomX}(\bm{x}^i_{F_q}))_{y}|^2,
\\ \text{MSE}_G^{MomY} & = \frac{1}{N_{F_q}} \sum_{i=1}^{N_{F_q}}|(G_{1}^{MomY}(\bm{x}^i_{F_q}))_{x} +(G_{2}^{MomY}(\bm{x}^i_{F_q}))_{y}|^2,
\\ \text{MSE}_G^{Ene} & = \frac{1}{N_{F_q}} \sum_{i=1}^{N_{F_q}}|(G_{1}^{Ene}(\bm{x}^i_{F_q}))_{x} +(G_{2}^{Ene}(\bm{x}^i_{F_q}))_{y}|^2,
\\ \text{MSE}_G^{Entropy} & = \frac{1}{N_{F_q}} \sum_{i=1}^{N_{F_q}}|(-\phi_{1}^{Ene}(\bm{x}^i_{F_q}))_{x} +(-\phi_{2}^{Ene}(\bm{x}^i_{F_q}))_{y} + \epsilon |^2,
\end{align*}
where $\epsilon$ is a small positive constant. The MSE for non-zero Dirichlet boundary condition is given as
\begin{align*}
  \text{MSE}_{\xi}^D &= \frac{1}{N_{\xi_q}} \sum_{i=1}^{N_{\xi_q}}| \xi^i - \xi_{{\tilde{\boldsymbol{\Theta}}}_q}(\bm{x}^i_{\xi_q})|^2.
\end{align*}
Along the common interface, the MSE average term for all primitive variables represented by a common symbol $\xi$ and residual continuity condition ($\text{MSE}_{\mathcal{R}}$) are given as
\begin{align*}
 \text{MSE}_{\xi_{avg}}(\tilde{\boldsymbol{\Theta}}_q;\{\bm{x}^{(i)}_{I_q}\}_{i=1}^{N_{I_q}} )  &= \sum_{\forall q^+} \left( \frac{1}{N_{I_q}} \sum_{i=1}^{N_{I_q}}\left| \xi_{{\tilde{\boldsymbol{\Theta}}}_q}(\bm{x}^{(i)}_{I_q})-  \left\{\!\left\{\xi_{{\tilde{\boldsymbol{\Theta}}}_q}(\bm{x}^{(i)}_{I_q})\right\}\!\right\}\right|^2 \right),
\\ \text{MSE}_{\mathcal{R}}(\tilde{\boldsymbol{\Theta}}_q;\{\bm{x}^{(i)}_{I_q}\}_{i=1}^{N_{I_q}} ) & = \sum_{\forall q^+} \left( \frac{1}{N_{I_q}} \sum_{i=1}^{N_{I_q}}\left|  \mathcal{F}_{{\tilde{\boldsymbol{\Theta}}}_q} (\bm{x}^{(i)}_{I_q}) - \mathcal{F}_{{\tilde{\boldsymbol{\Theta}}}_{q^+}}(\bm{x}^{(i)}_{I_q})\right|^2 \right),
\end{align*}
respectively. Notation $\{\!\{ \cdot \}\!\}$ represent the average value of primitive variables along the common interface.
The term $\text{MSE}_{\nabla \rho} = \frac{1}{N_{\nabla_q}} \sum_{i=1}^{N_{\nabla_q}}| \nabla \rho^i - \nabla \rho_{{\tilde{\boldsymbol{\Theta}}}_q}(\bm{x}^i_{\nabla_q})|^2$ represents the data mismatch term for density gradient obtained from the Schlieren photography experimental technique traditionally used in the study of high-speed aerodynamics.

As discussed in \cite{jagtap2020extended}, the XPINN can reduce the generalization error by carefully selecting residual points in each subdomain. Moreover, with the XPINN methodology, the computational domain can be decomposed in any arbitrary way according to the required number and the location of residual points. This can be achieved with some prior knowledge about the solution behavior such that the localized powerful neural network can be employed in that subdomain. The XPINN can also reduce the approximation error by selecting the network size, and other hyperparameters in each subdomain depending on the complexity of the solution. Furthermore, XPINN can also reduce optimization error by using the different network size for each subdomain.

\subsection{Optimization Methods}
By minimizing the loss function, we seek the optimal parameters of the network $\tilde{\boldsymbol{\Theta}}^*$ from the parameter space.  For this purpose there are several optimization algorithms available that can be used to minimize loss function. The gradient-based methods are one of the most popular class of optimization methods that can be employed for the minimization of loss functions. In the basic form, given an initial value of tunable parameters $\tilde{\boldsymbol{\Theta}}$, these parameters are updated as
$ \tilde{\boldsymbol{\Theta}}^{m+1} = \tilde{\boldsymbol{\Theta}}^m - \eta_l \left. \frac{\partial \mathcal{J}(\tilde{\boldsymbol{\Theta}}) }{\partial \tilde{\boldsymbol{\Theta}}} \right|_{\tilde{\boldsymbol{\Theta}} = \tilde{\boldsymbol{\Theta}}^m}$,
where the learning rate is denoted by $\eta_l$. In particular, we shall use the Adam optimizer \cite{kingma2014Adam} followed by limited-memory Broyden–Fletcher–Goldfarb–Shanno (L-BFGS) optimizer \cite{byrd1995limited}, which is a quasi-Newton optimization method. For smooth and regular solutions the L-BFGS optimizer can find a better solution with a small number of iterations compared to the Adam optimizer, due to second-order accuracy as opposed to Adam, which is first-order accurate but in general more robust. Hence, we often start with the Adam optimizer and after we reach a small value of the loss function we switch to L-BFGS.

\section{Computational results}
In this section we shall solve several examples including expansion and shock waves problems. In appendix A we present a pedagogical example involving the inverse problem for the two-dimensional Euler equations with smooth solutions. We provide all the details for this simple example so that the interested reader can repeat our results.

In this work, the comparison study has been performed between the PINN and XPINN methods based on  the following points.

\begin{itemize}
    \item In both PINNs and XPINNs, we used unity weights for the different MSE terms in the loss function for all the problems, except for the case where dynamics weights are employed. 
    \item It is important to note that in this study we do not aim to compare the computational cost associated with PINNs and XPINNs, because such study has already been performed in \cite{shukla2021parallel}. The main objective of this work is to solve the ill-posed inverse problems in supersonic flow, which has many real-world applications in designing aerospace vehicles. 
    \item We are solving an inverse problem with given data inside the domain as well as on the boundary of the domain. Such data availability allows us to decompose the computational domain into several smaller subdomains as per the compleixity of the data, especially for the XPINN method. 
    \item For all the test cases, the primitive variable data, which is available in the dimensional form is converted into the non-dimensional one for the network training.
 \end{itemize}

\subsection{Expansion wave problem}
\begin{figure}[http]
    \centering
    \subfigure{
    \includegraphics[width=0.4\textwidth]{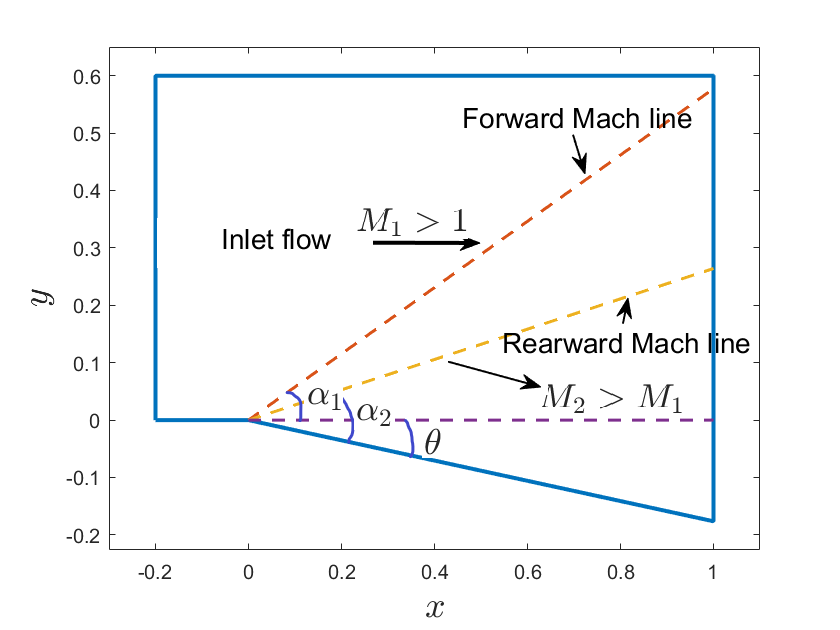}
    }
    \subfigure{
    \includegraphics[width=0.42\textwidth]{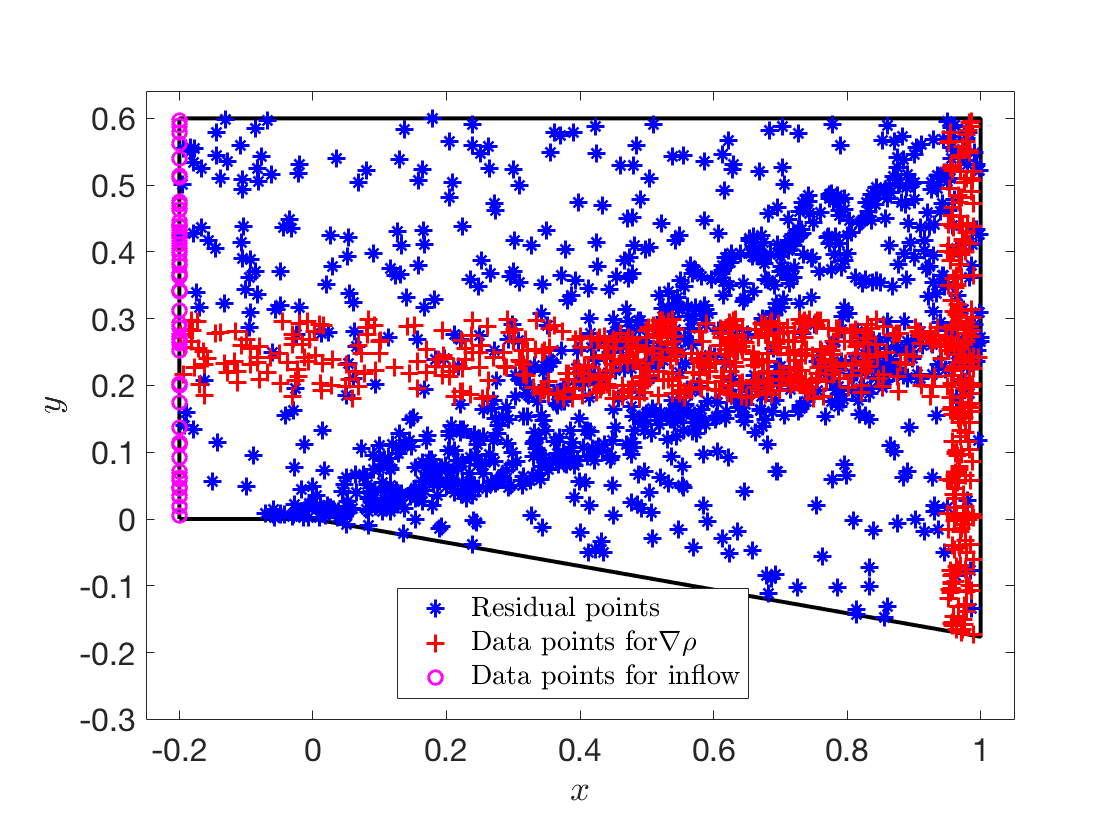}
    }
    \caption{Expansion wave problem: (Left) schematic of the expansion waves, and  (right) distributions of the 1200 residual points (blue $\ast$ points, where we compute the residuals), the 340 data points for the density gradient (red $+$ points) and the 100 data points for the inflow conditions (magenta $\circ$ points).}
    \label{fig:distribution:expansion}
\end{figure}

A supersonic expansion fan is a centered expansion process that occurs when a supersonic flow turns around a convex corner. The expansion fan consists of an infinite number of Mach waves, diverging from a sharp corner; see the left plot of Figure \ref{fig:distribution:expansion} for the schematic.
To overcome the aforementioned boundary conditions issue, we are going to use PINNs and XPINNs to solve the expansion wave problem with the inlet flow conditions as well as the density gradient information from Schlieren photography; see Figure \ref{fig:distribution:expansion} (right).

We solve the inverse problem of the expansion wave problem by using PINN without using the outflow boundary conditions. Instead, we use the information of density gradient $\nabla \rho$ and a small number of pressure data points on the surface of the wedge, following the one-dimensional problem~\cite{mao2020physics}. In addition, we use the inlet conditions here. However, unlike the one-dimensional case where we use the density gradient in the entire computational domain, here we only use the density gradients in a small region of the computational domain.
By combining the mathematical model, we have the weighted loss function of PINN given by 
\begin{equation*}
   \mathcal{J} = \omega_{1}~\text{MSE}_{\mathcal{F}} +
    \omega_{2} ~\text{MSE}_{\nabla \rho|_{D}} + \omega_{3}~\text{MSE}_{Inflow} +  \omega_{4}~\text{MSE}_{p^{*}},
\end{equation*}
where $\text{MSE}_{\mathcal{F}}$ corresponds to the Euler equation \eqref{NCL:steady} including the entropy conditions, $\text{MSE}_{Inflow}$ corresponds to the inflow boundary conditions on primitive variables, $\text{MSE}_{\nabla \rho|_{D}}$ corresponds to the density gradient and $D\subset \Omega$ is a small region of the computational domain, $\text{MSE}_{p^*}$ corresponds to the pressure data on the wall surface. We also employ dynamic weights ~\cite{wang2021understanding} denoted as $\omega_i, i =1,2,\cdots$ in front of all the MSE terms in the loss function. Hereafter, for brevity, the loss function $\mathcal{J}$ is denoted without showing its dependence on the network parameters $\tilde{\boldsymbol{\Theta}}$.

Consider the following example:
\textit{Let $\gamma = 1.4$. Assume the inlet flow for a tuning corner with $\theta = 10^\circ$ is given by 
$$M_{\infty} = 2,\; \rho_\infty = 1.23kg/m^3,\; u_\infty = 678.1m/s,\; v_\infty = 0,\; p_\infty = 1.01\times 10^{5} Pa, \; T_\infty = 286.1 K.$$
}
The exact solutions can be obtained by using the isentropic relations~\cite{anderson1995computational}.
The distributions of the residual points and the data points of the density gradient are shown in Fig. \ref{fig:distribution:expansion} (right). For the pressure, we use 50 randomly distributed data points on the surface of the wedge, i.e., $y = -tanh(\theta) x,\, x\in(0, 1).$

\begin{figure}[http]
    \centering
    \subfigure[Exact solution]{
    \includegraphics[trim= 0cm 0cm 0cm 0cm, scale=0.16]{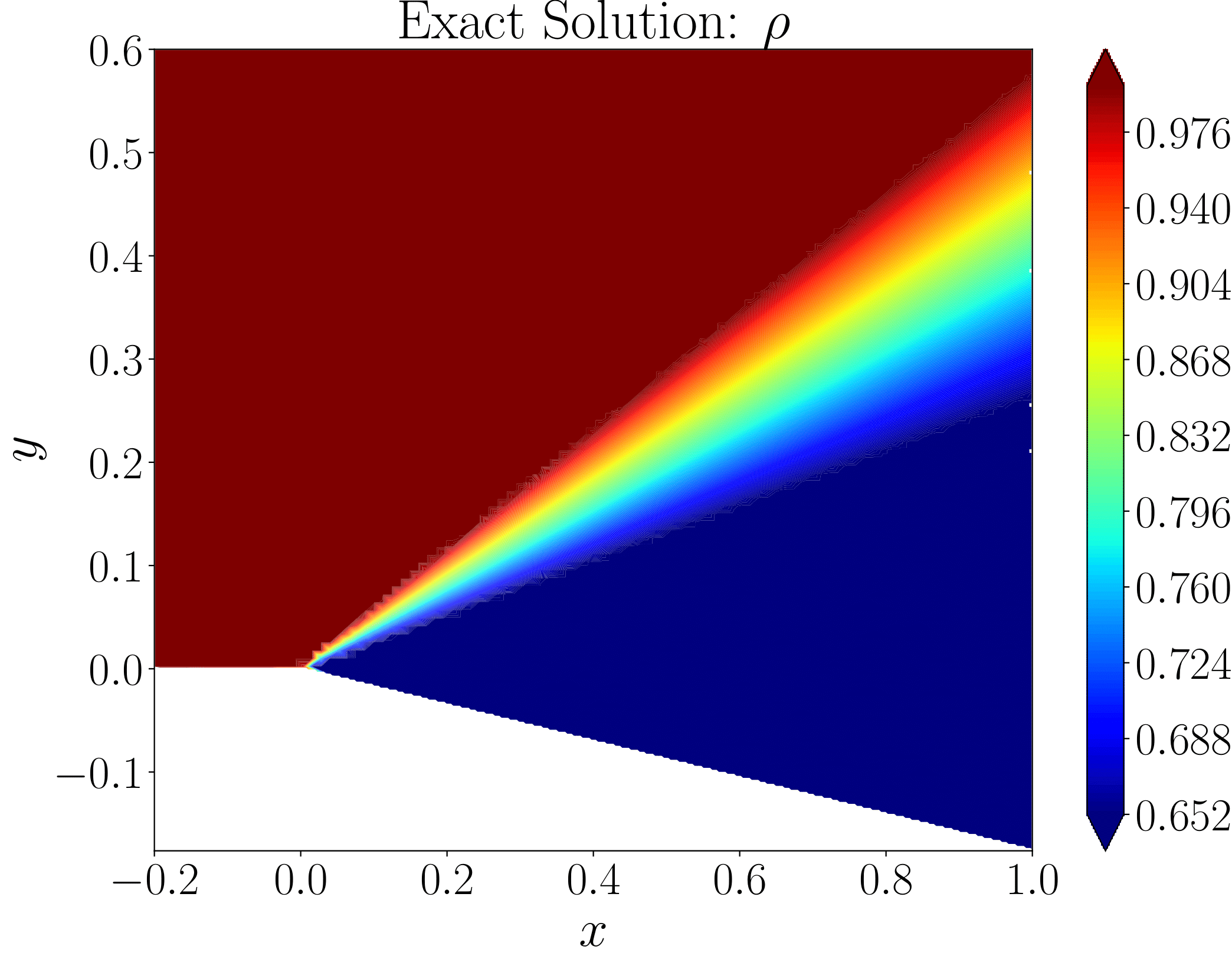}
    }
    \subfigure[PINN solution with DW]{
    \includegraphics[trim= 0cm 0cm 0cm 0cm, scale=0.16]{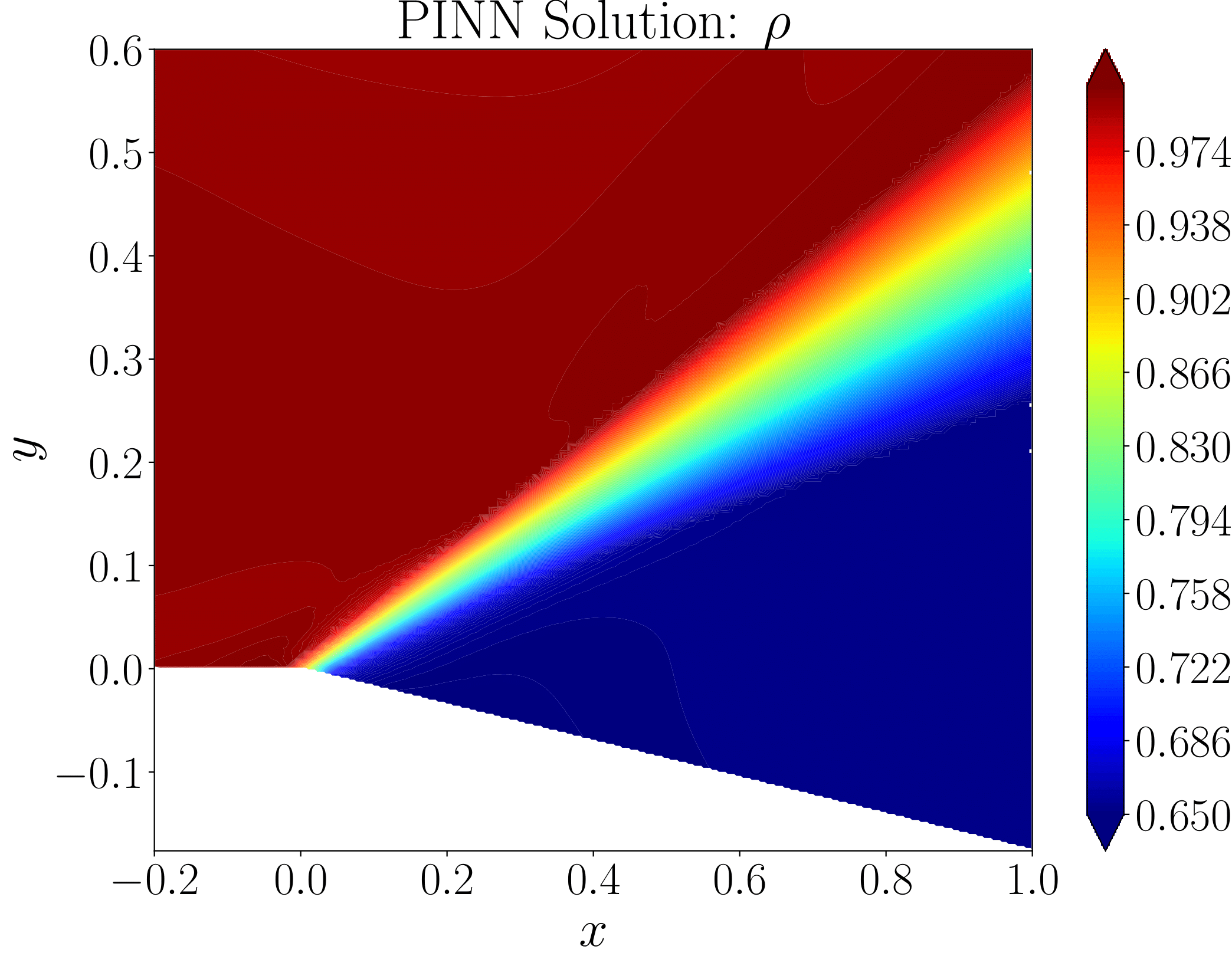}
    }
    \subfigure[Relative error with DW]{
   \includegraphics[trim= 0cm 0cm 0cm 0cm, scale=0.16]{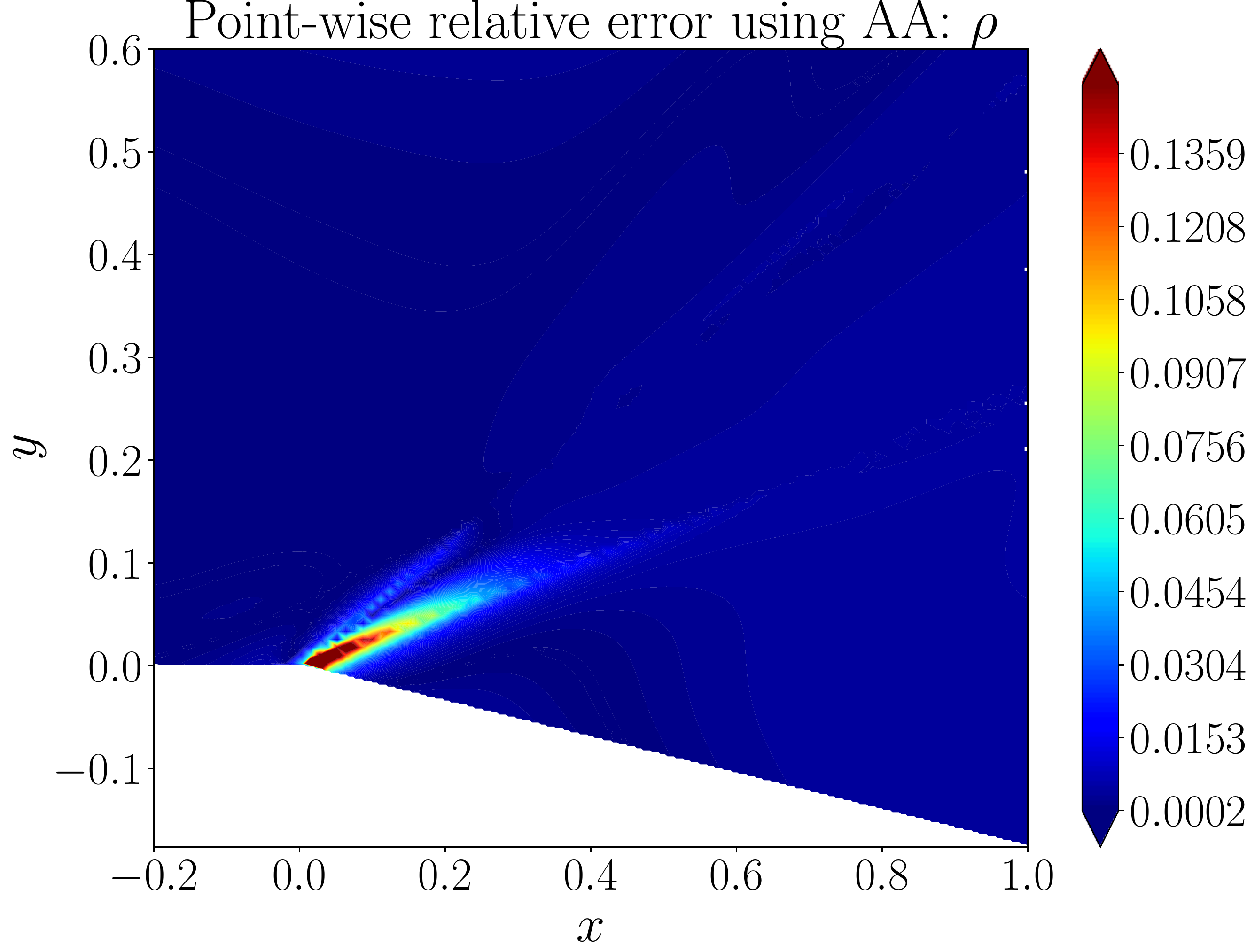}
    }
    \subfigure[Relative error with AA]{
    \includegraphics[trim= 0cm 0cm 0cm 0cm,scale=0.16]{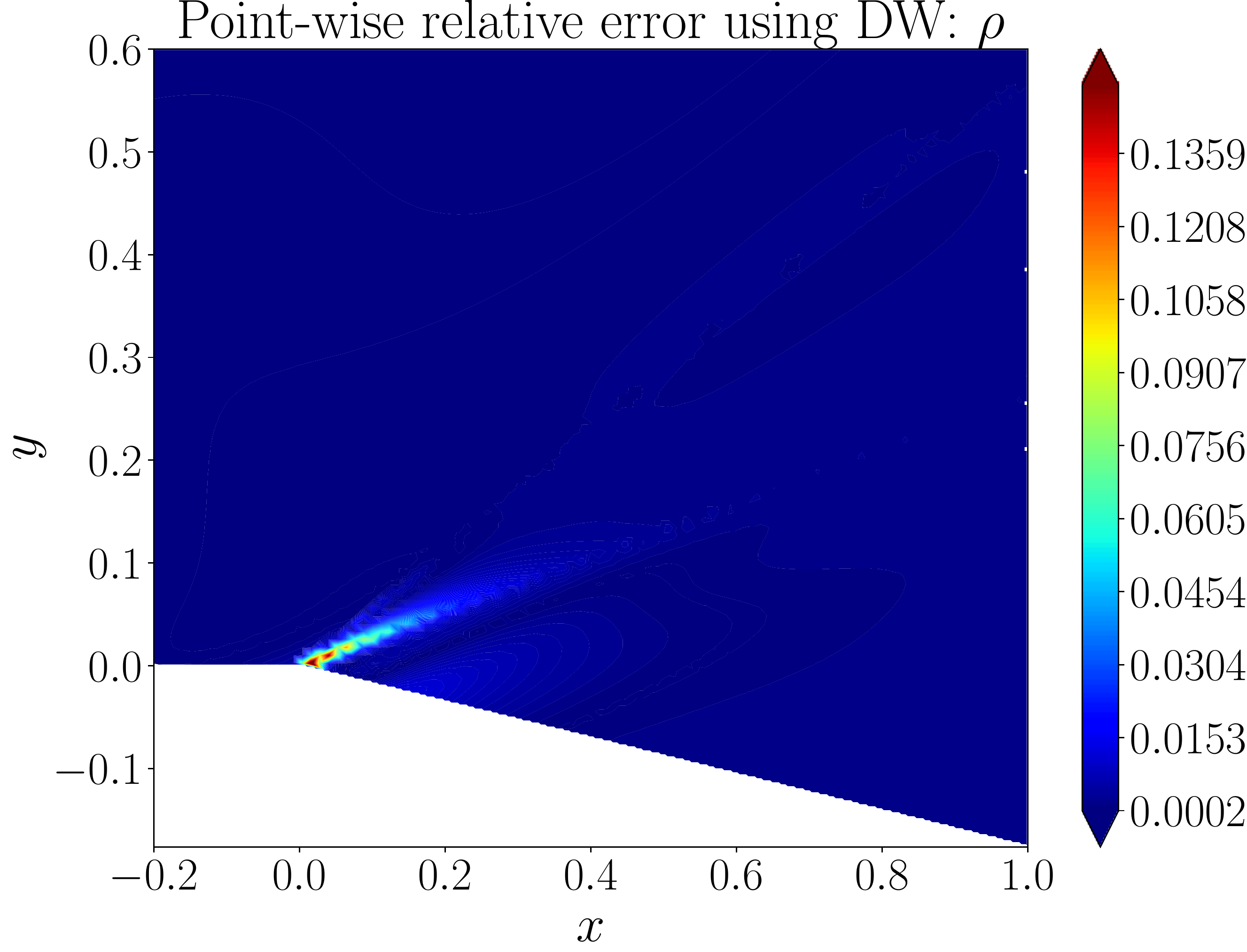}
    }
    \subfigure[Relative error without DW/AA]{
   \includegraphics[trim= 0cm 0cm 0cm 0cm,scale=0.16]{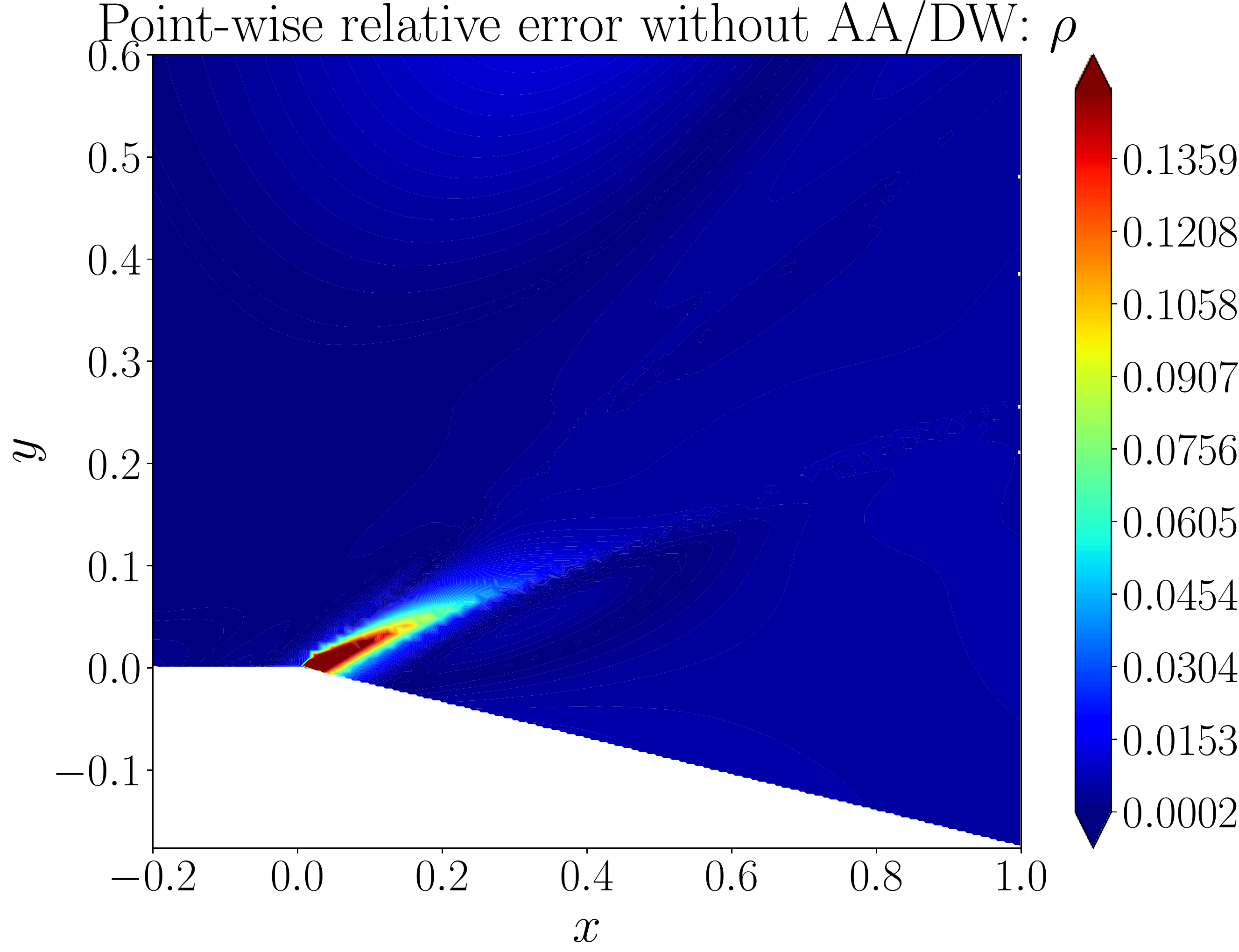}
    }
    \subfigure[ Relative error with XPINN.]{
    \includegraphics[trim= 0cm 0cm 0cm 0cm,scale=0.16]{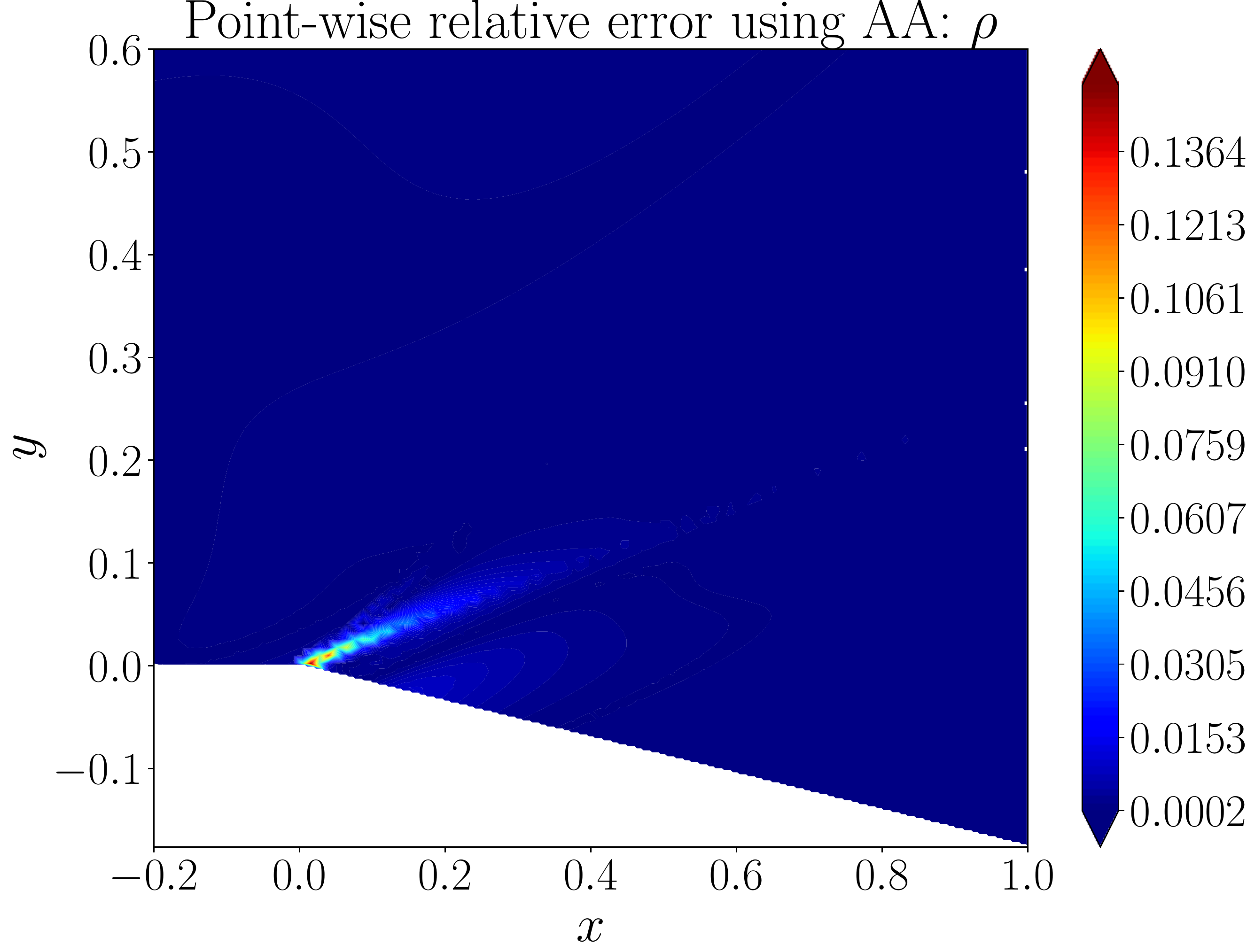}
    }
    \caption{Expansion wave problem: Comparison of the {\emph{density}} between PINN solutions with/without adaptive activation function (AA) or dynamic weights (DW). Here we use `tanh' as the basic activation function. (a) Exact solution. (b) PINN solution using dynamic weights. (c) Relative point-wise error with dynamic weights. (d) Relative point-wise error with adaptive activation function. (e) Relative point-wise error without dynamic weights or adaptive activation function. (f)  Relative point-wise error with XPINN using adaptive activations.}
    \label{fig:expansion:rho:compa}
\end{figure}

\begin{figure}[http]
    \centering
    \subfigure[Exact solution]{
    \includegraphics[scale=0.16]{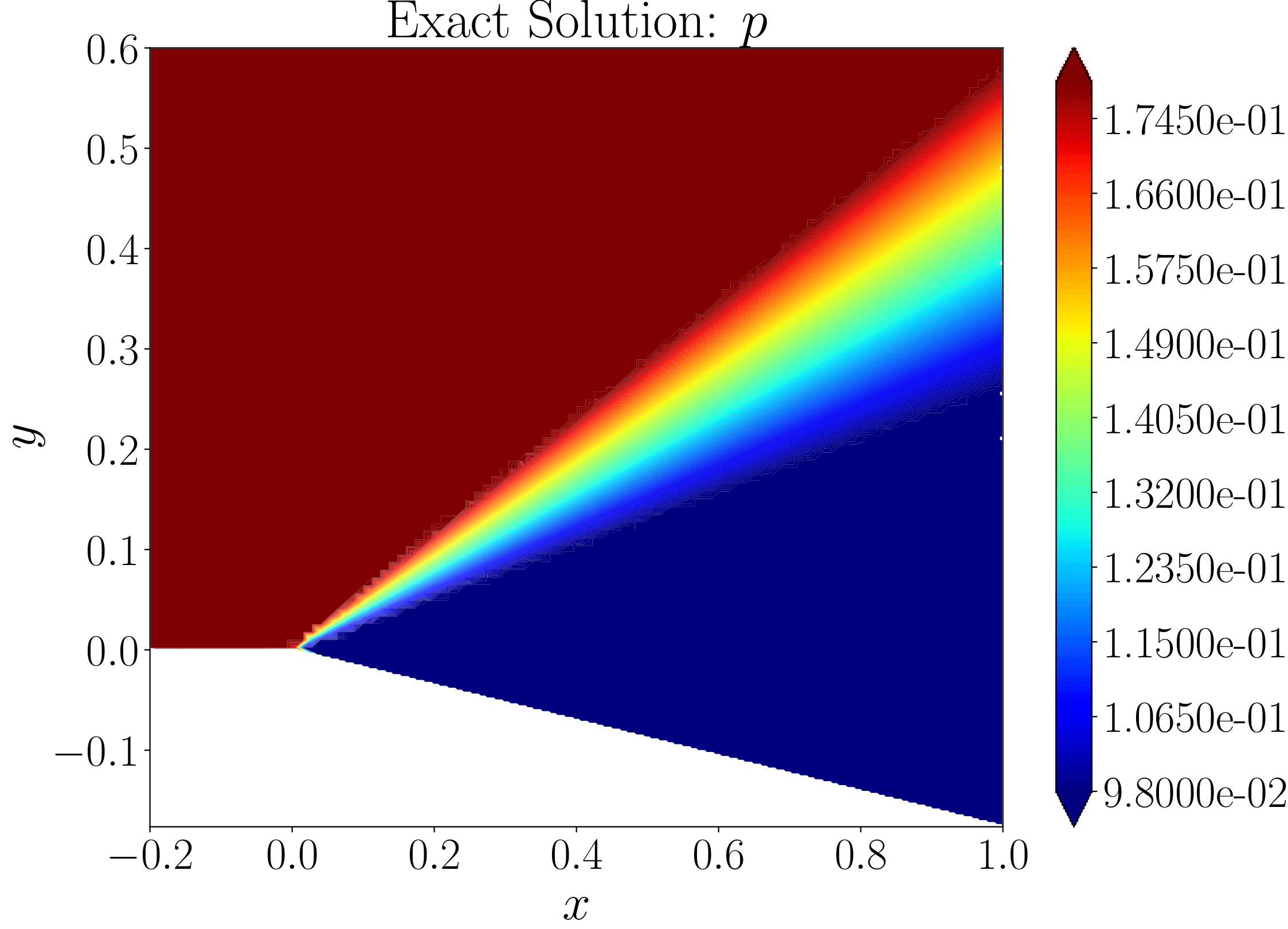}
    }
    \subfigure[PINN solution with DW]{
    \includegraphics[scale=0.16]{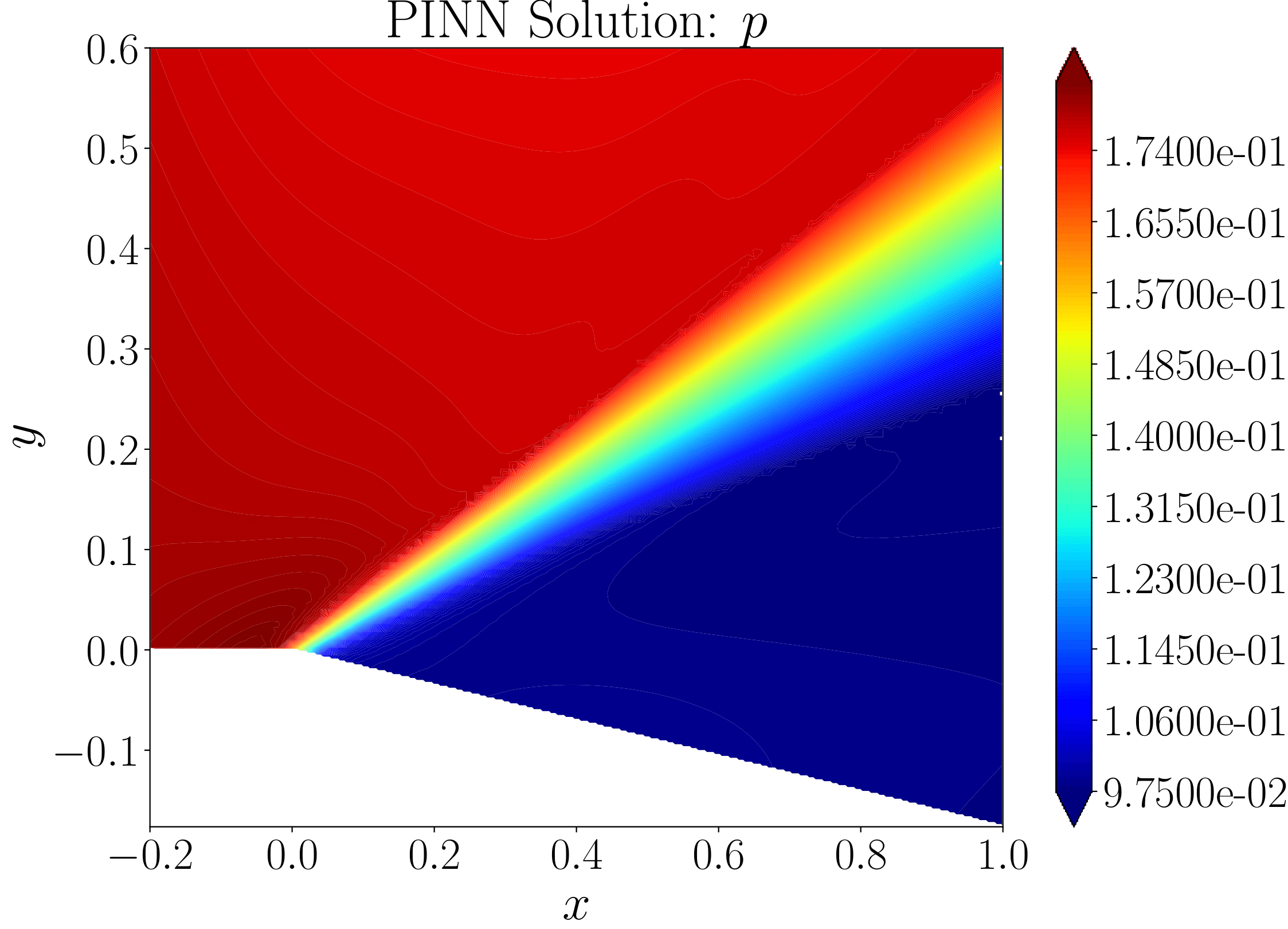}
    }
    \subfigure[Relative error with DW]{
    \includegraphics[scale=0.16]{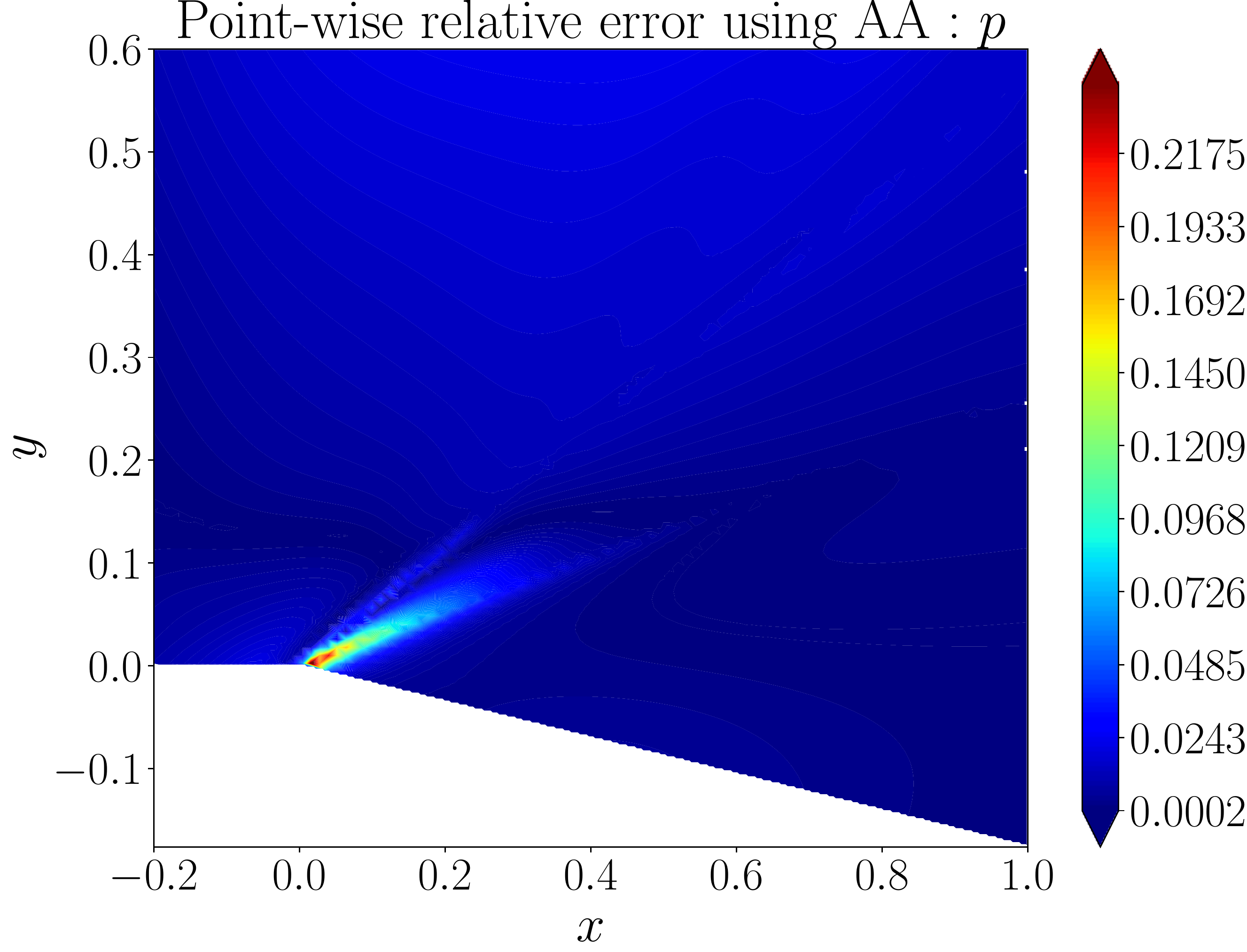}
    }
    \subfigure[Relative error with AA]{
    \includegraphics[scale=0.16]{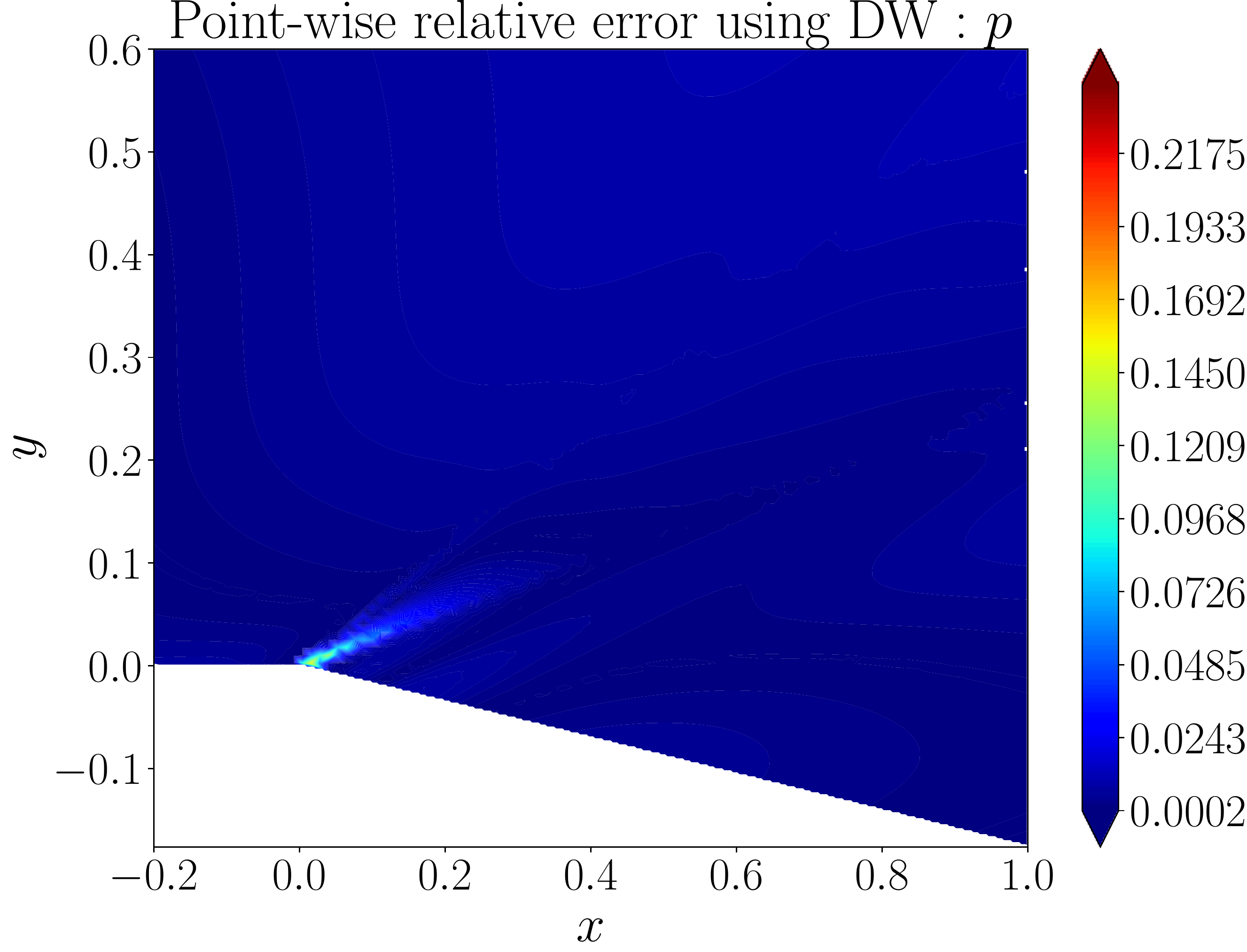}
    }
    \subfigure[Relative error without DW/AA]{
    \includegraphics[scale=0.16]{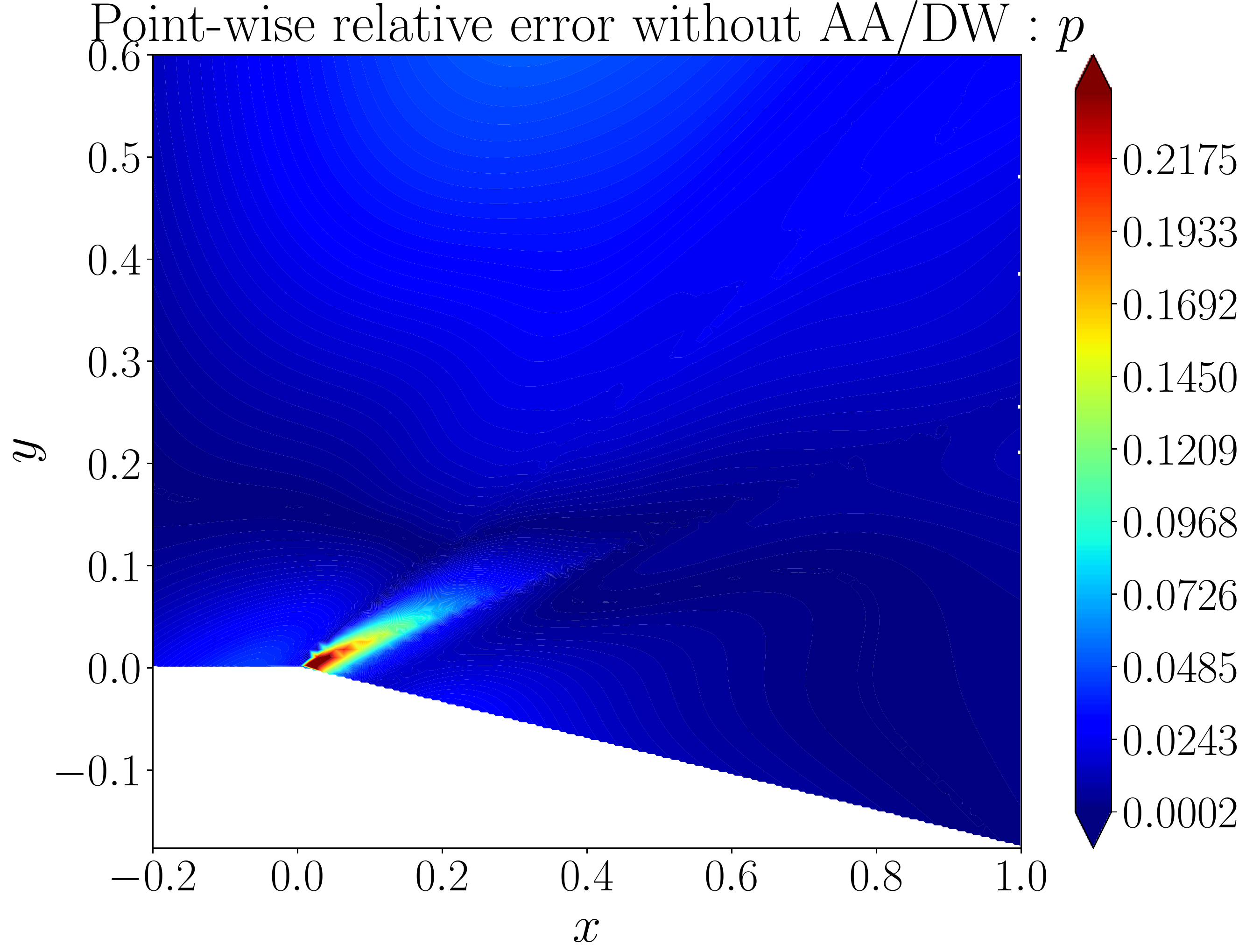}
    }
    \subfigure[ Relative error with XPINN.]{
    \includegraphics[scale=0.16]{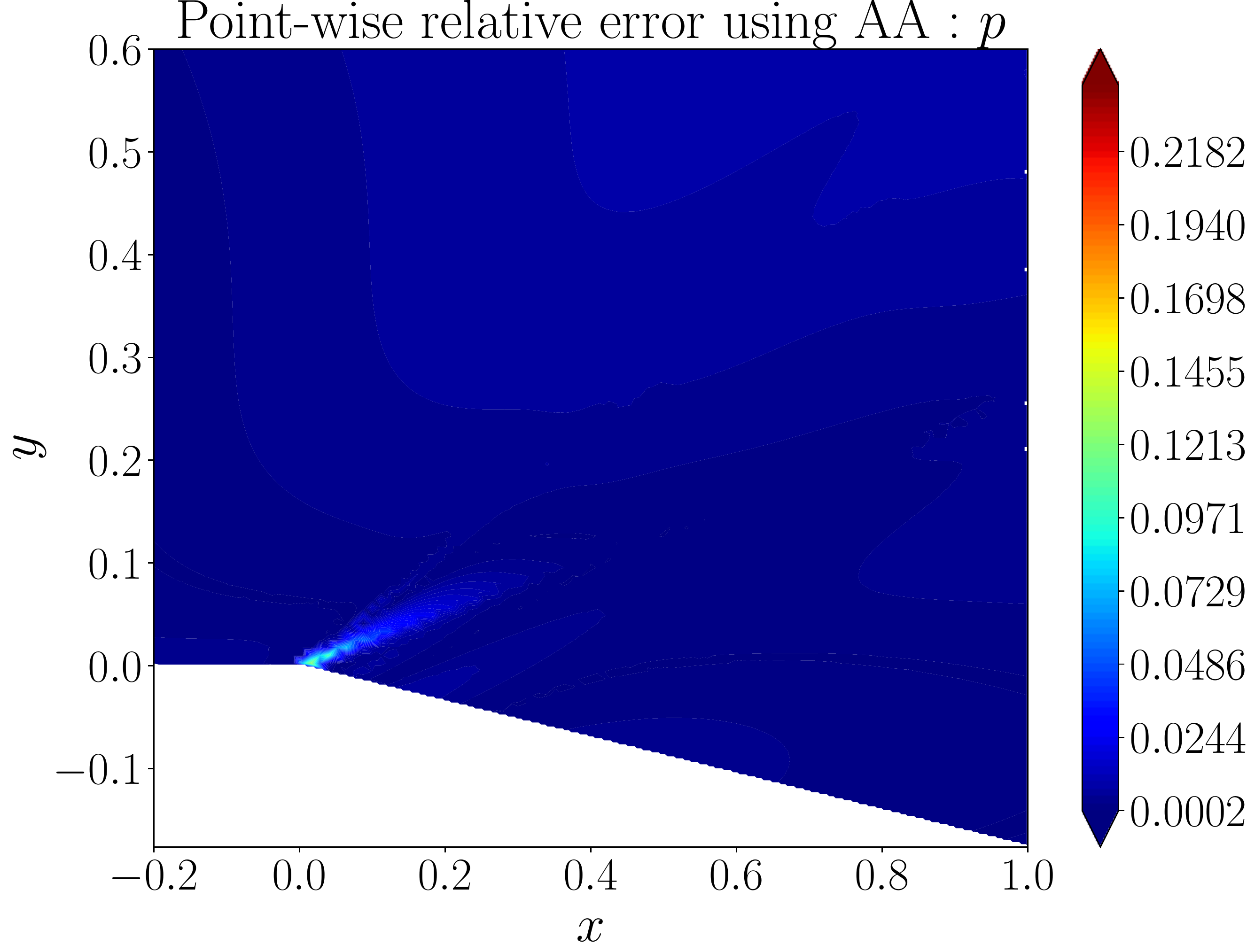}
    }
    \caption{Expansion wave problem: Comparison of the {\emph{pressure}} between PINN solutions with/without adaptive activation function (AA) or dynamic weights (DW). Here we use `tanh' as the basic activation function. (a) Exact solution. (b) PINN solution using dynamic weights. (c) Relative point-wise error with dynamic weights. (d) Relative point-wise error with adaptive activation function. (e) Relative point-wise error without dynamic weights or adaptive activation function. (f)  Relative point-wise error with XPINN using adaptive activations.}
    \label{fig:expansion:p:compa}
\end{figure}

\begin{figure}[http]
    \centering
    \subfigure[Exact solution]{
    \includegraphics[scale=0.16]{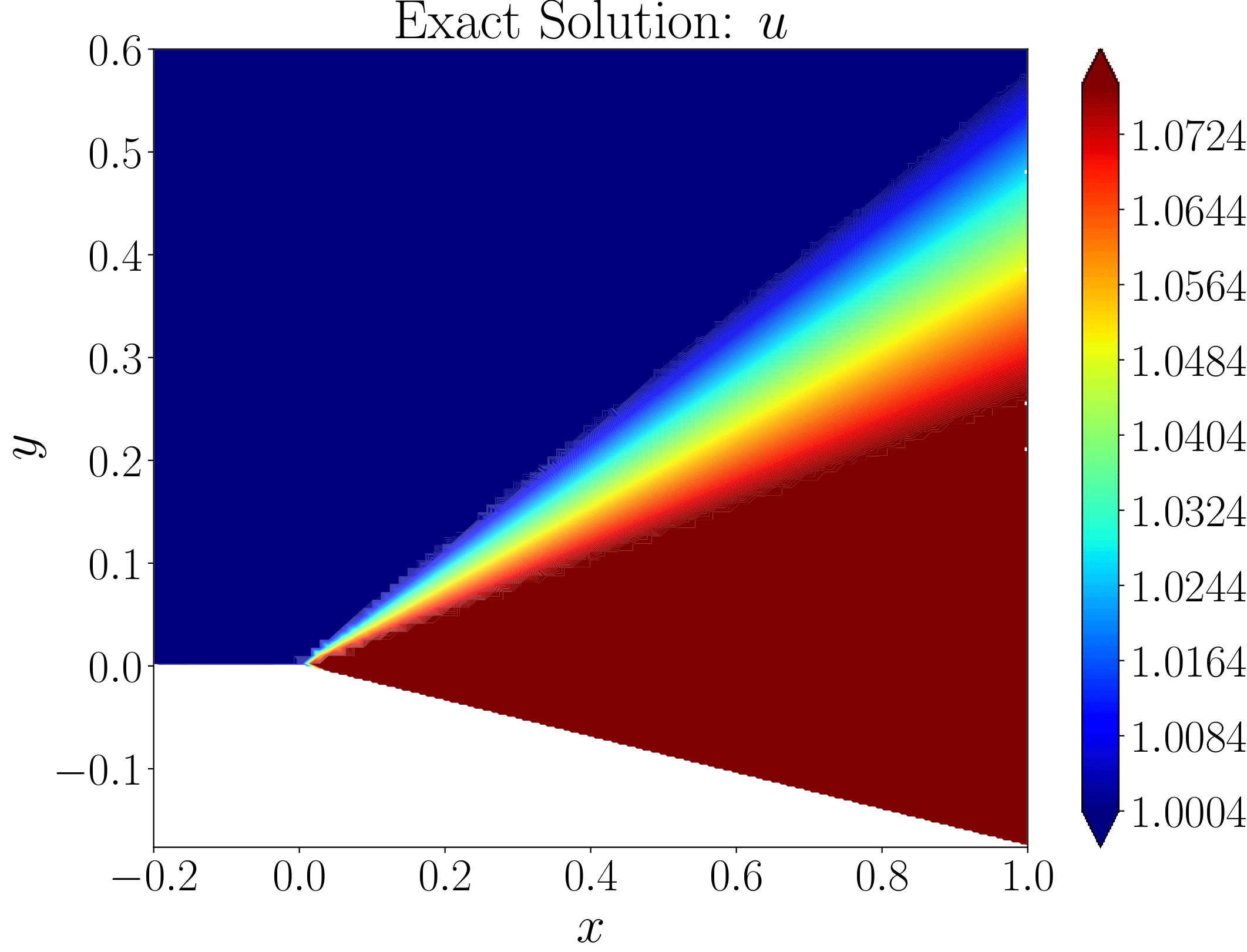}
    }
    \subfigure[PINN solution with DW]{
    \includegraphics[scale=0.16]{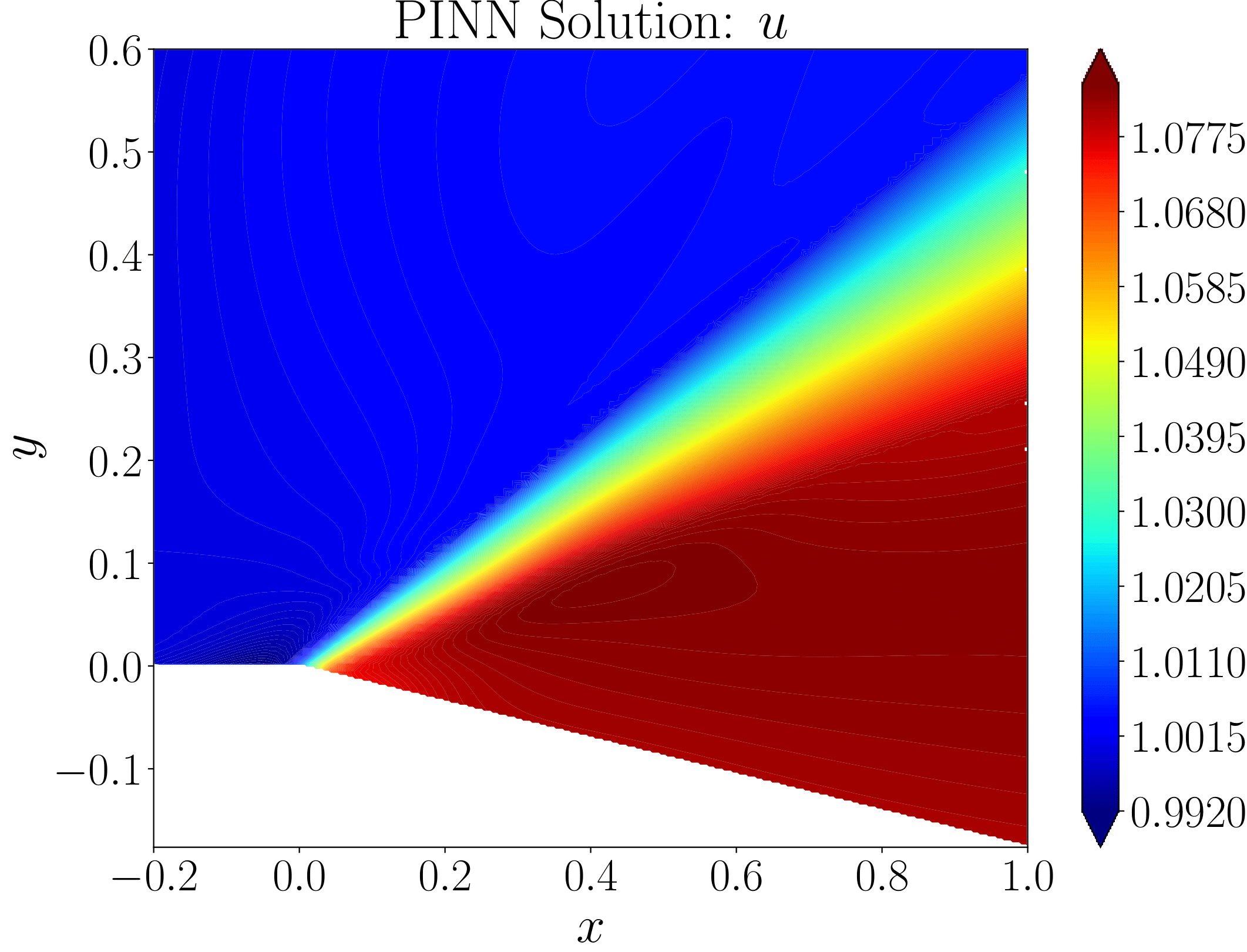}
    }
    \subfigure[Relative error with DW]{
   \includegraphics[scale=0.16]{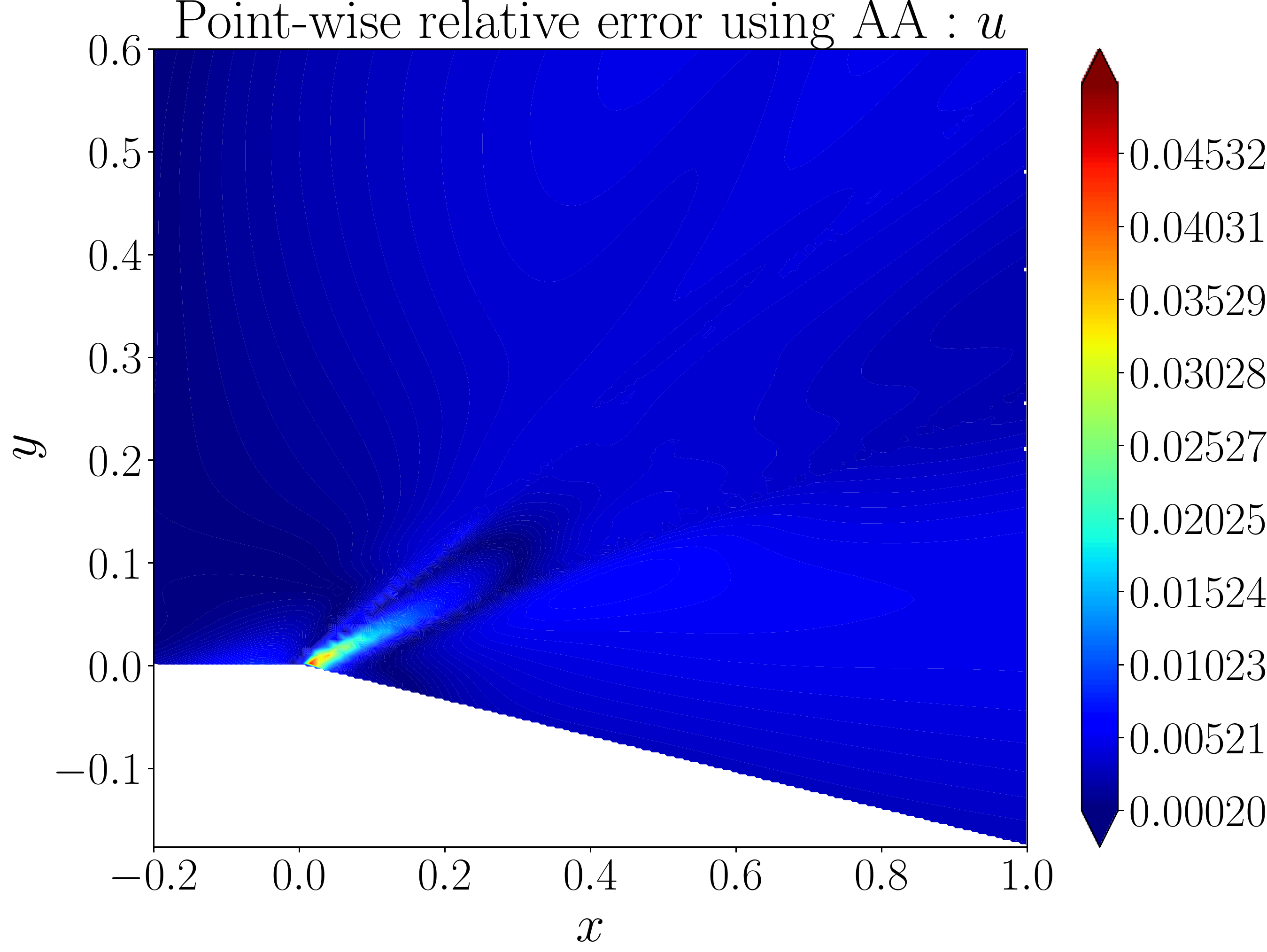}
    }
    \subfigure[Relative error with AA]{
    \includegraphics[scale=0.16]{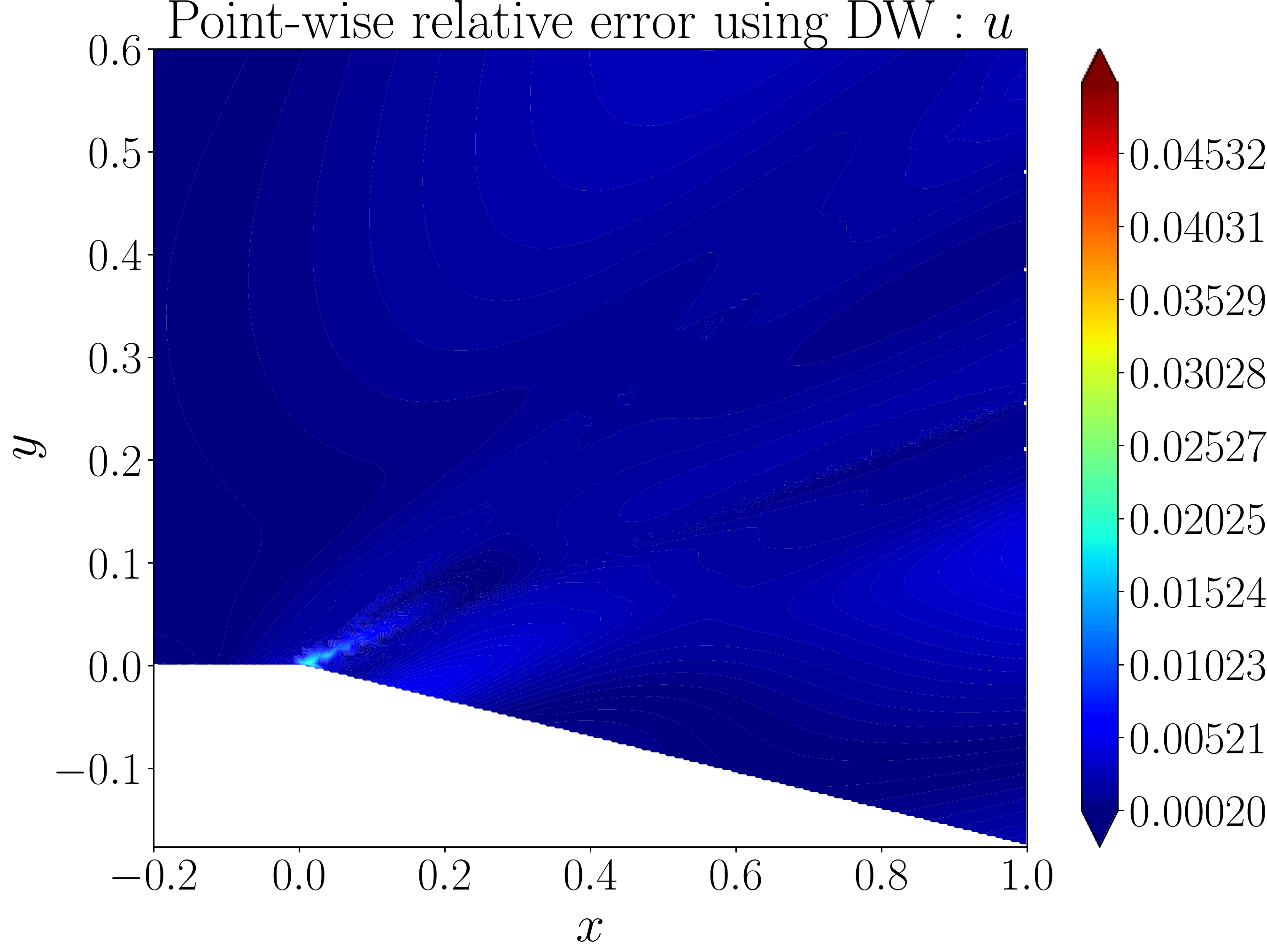}
    }
    \subfigure[Relative error without DW/AA]{
    \includegraphics[scale=0.16]{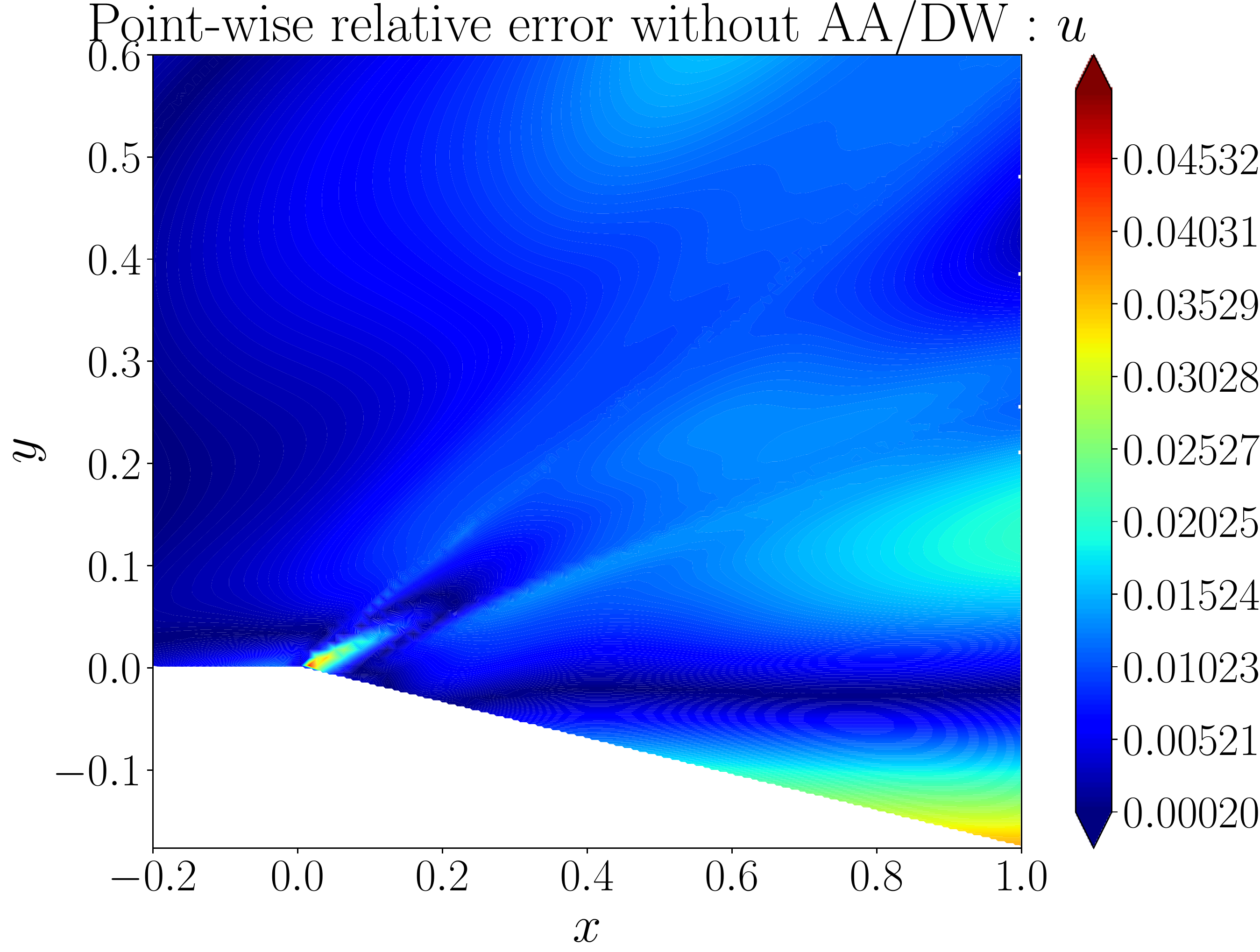}
    }
    \subfigure[ Relative error with XPINN.]{
    \includegraphics[scale=0.16]{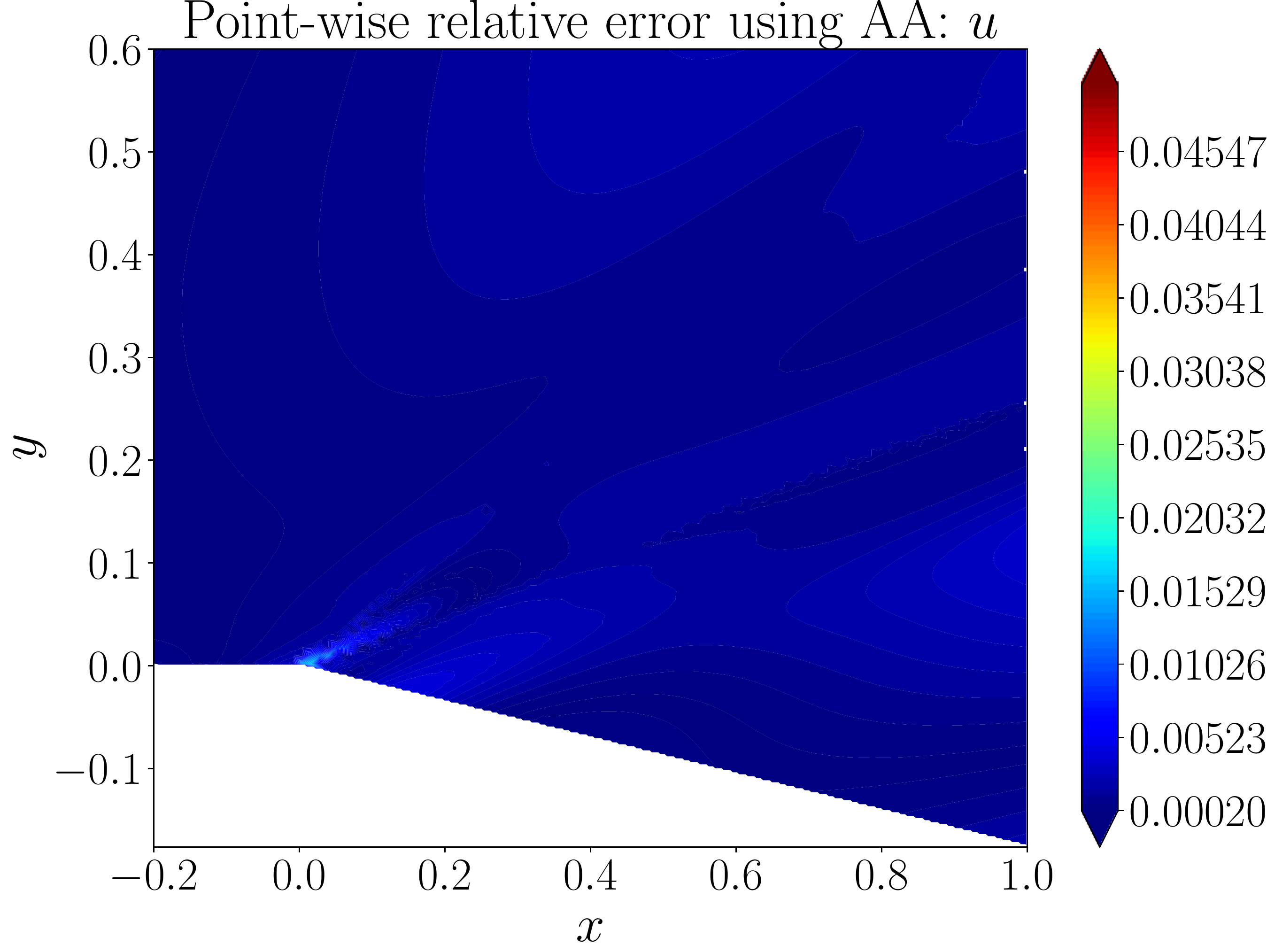}
    }
    \caption{Expansion wave problem: Comparison of the {\emph{velocity $u$}} between PINN solutions with/without adaptive activation function (AA) or dynamic weights (DW). Here we use `tanh' as the basic activation function. (a) Exact solution. (b) PINN solutions using dynamic weights. (c) Relative point-wise error with dynamic weights. (d) Relative point-wise error with adaptive activation function. (e) Relative point-wise error without dynamic weights or adaptive activation function. (f)  Relative point-wise error with XPINN using adaptive activations.}
    \label{fig:expansion:u:compa}
\end{figure}

\begin{figure}[http]
    \centering
    \subfigure[Exact solution]{
    \includegraphics[scale=0.17]{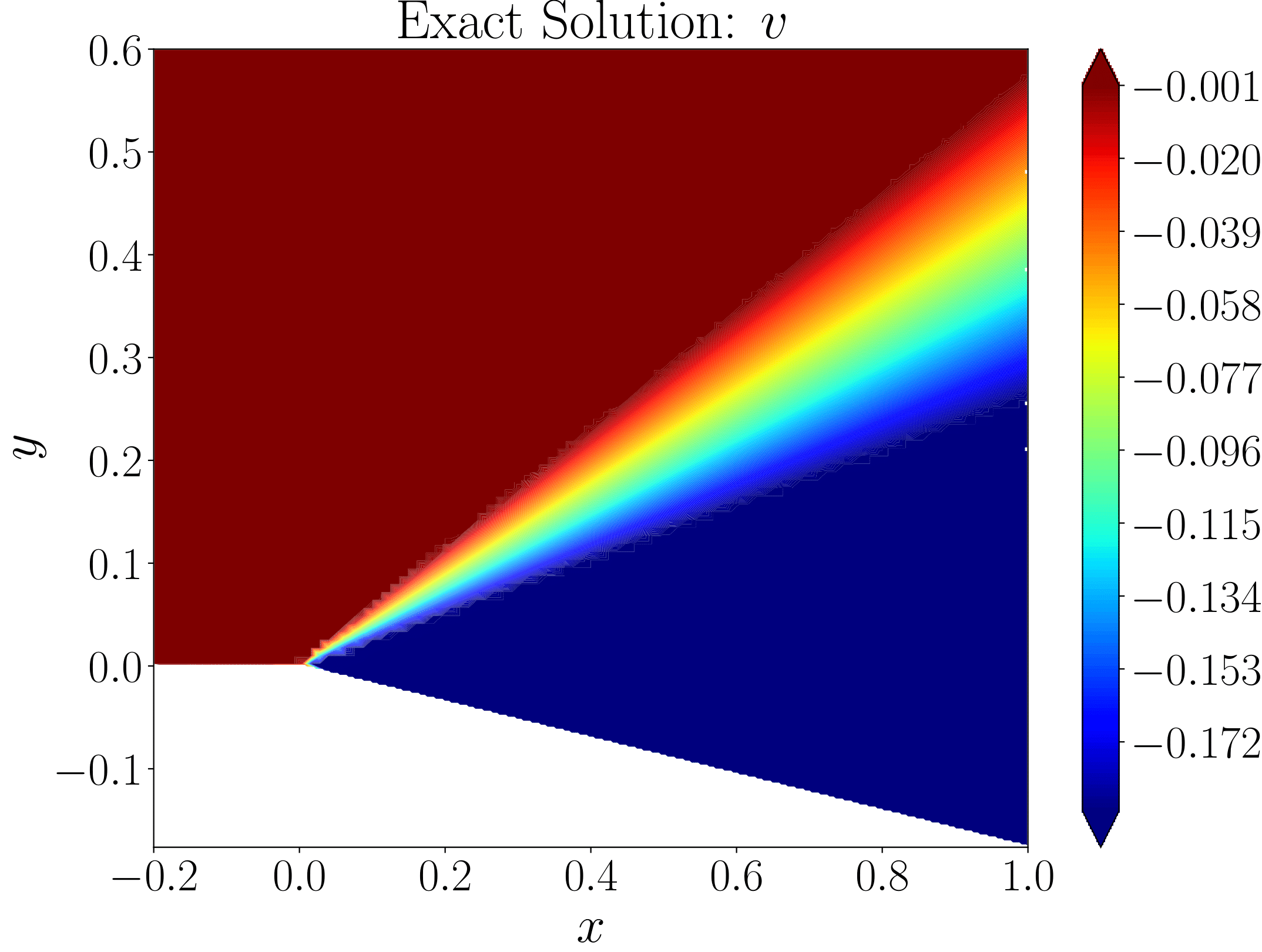}
    }
    \subfigure[PINN solution with DW]{
    \includegraphics[scale=0.17]{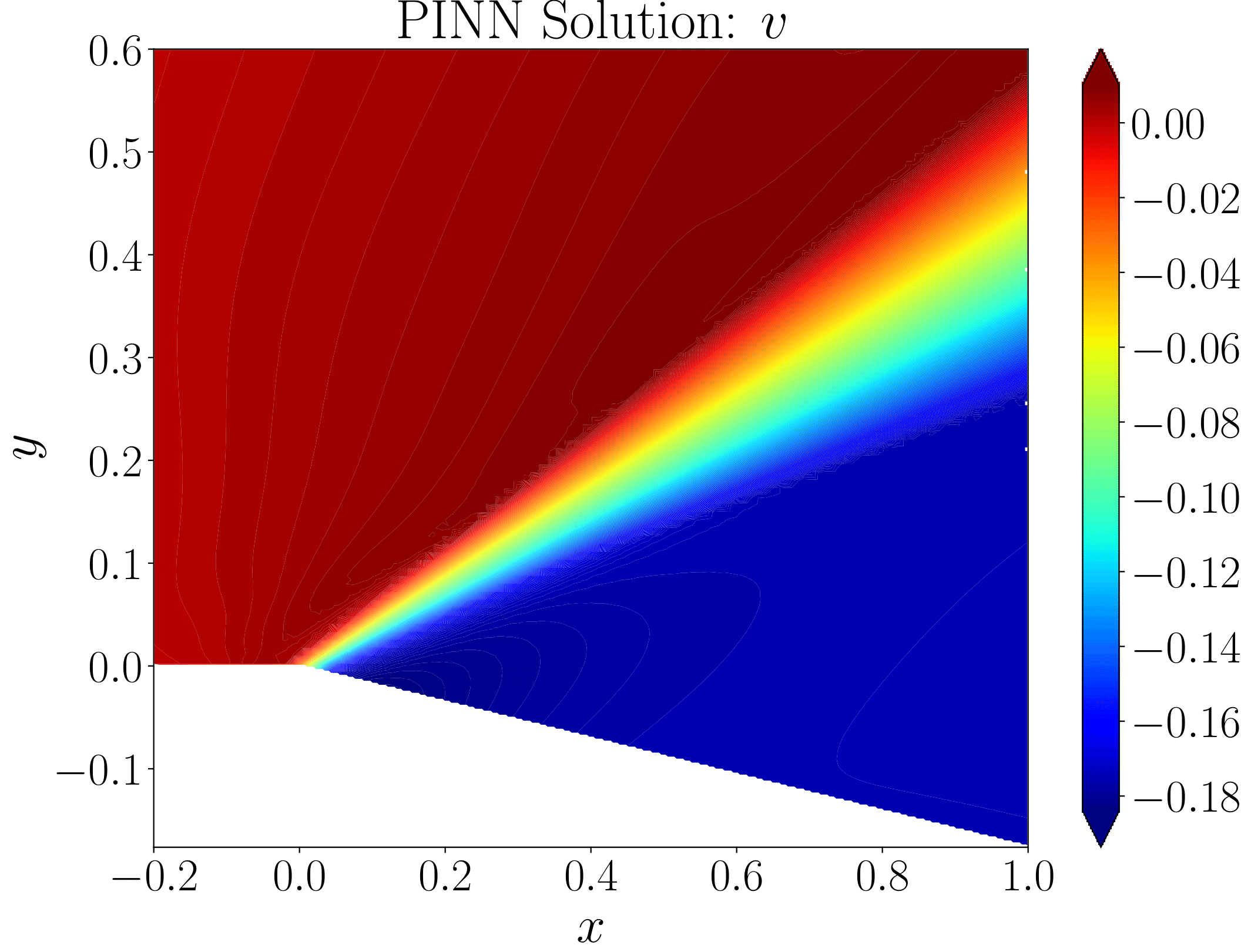}
    }
    \subfigure[Relative error with DW]{
    \includegraphics[scale=0.17]{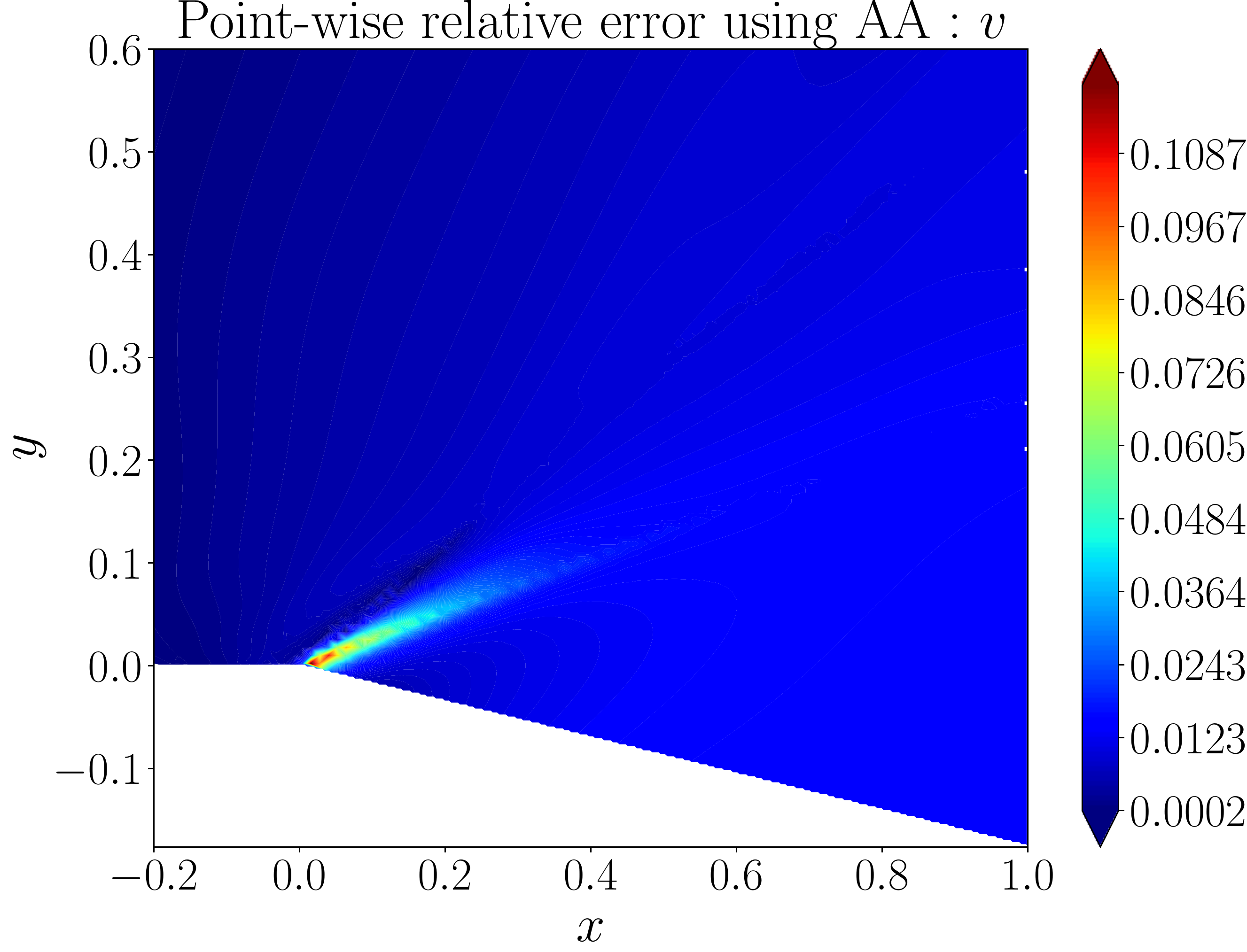}
    }
    \subfigure[Relative error with AA]{
    \includegraphics[scale=0.17]{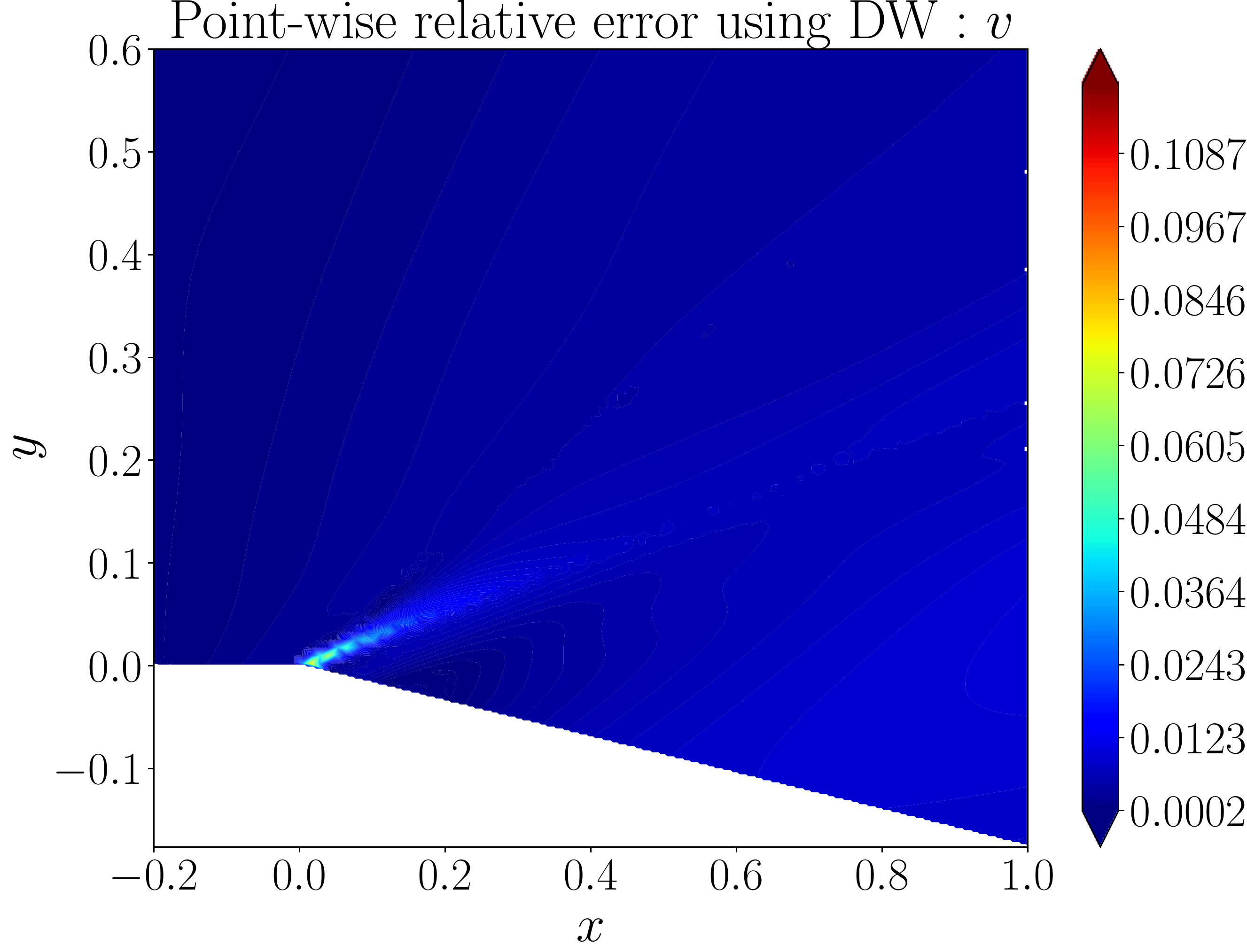}
    }
    \subfigure[Relative error without DW/AA]{
    \includegraphics[scale=0.17]{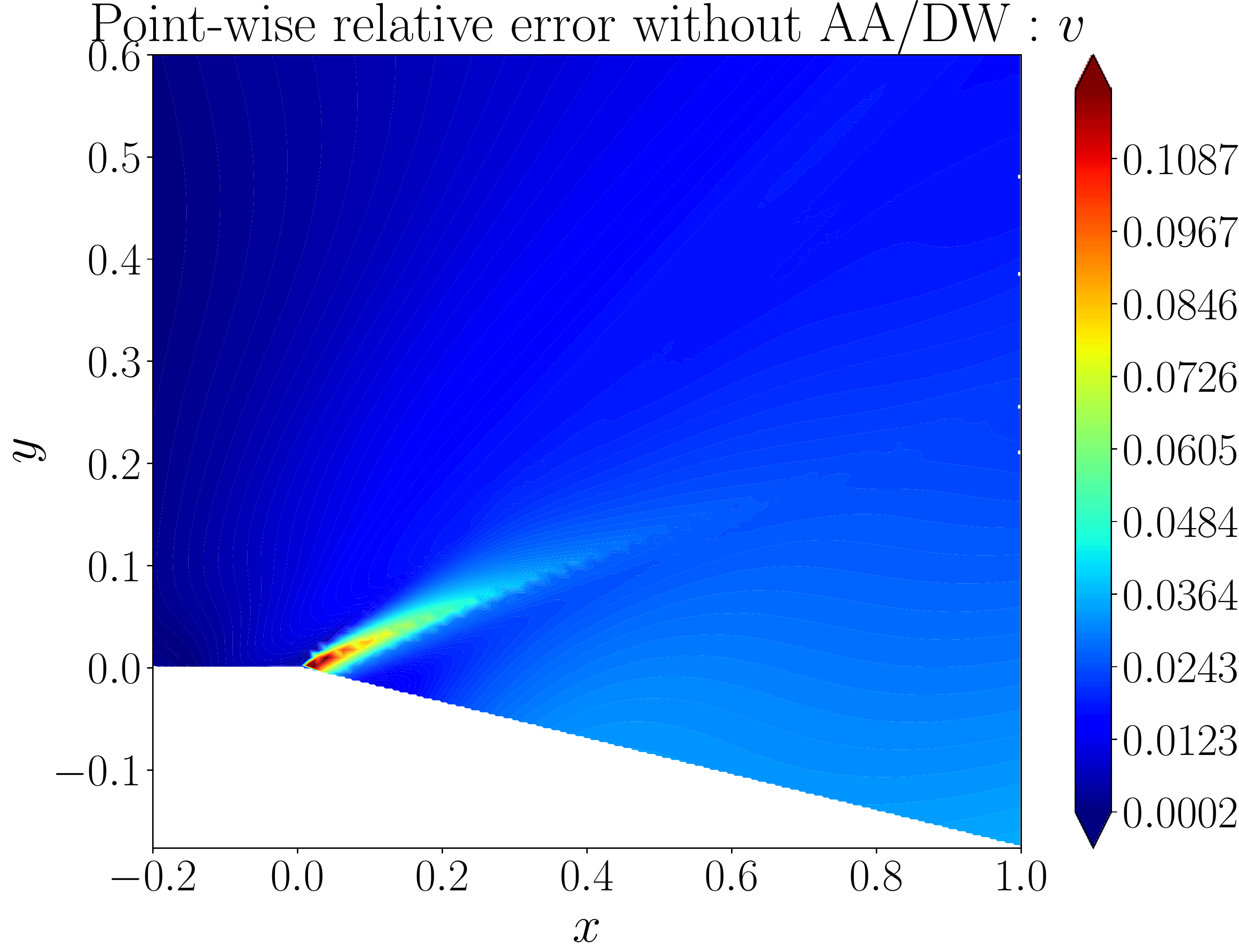}
    }
    \subfigure[Relative error with XPINN.]{
    \includegraphics[scale=0.17]{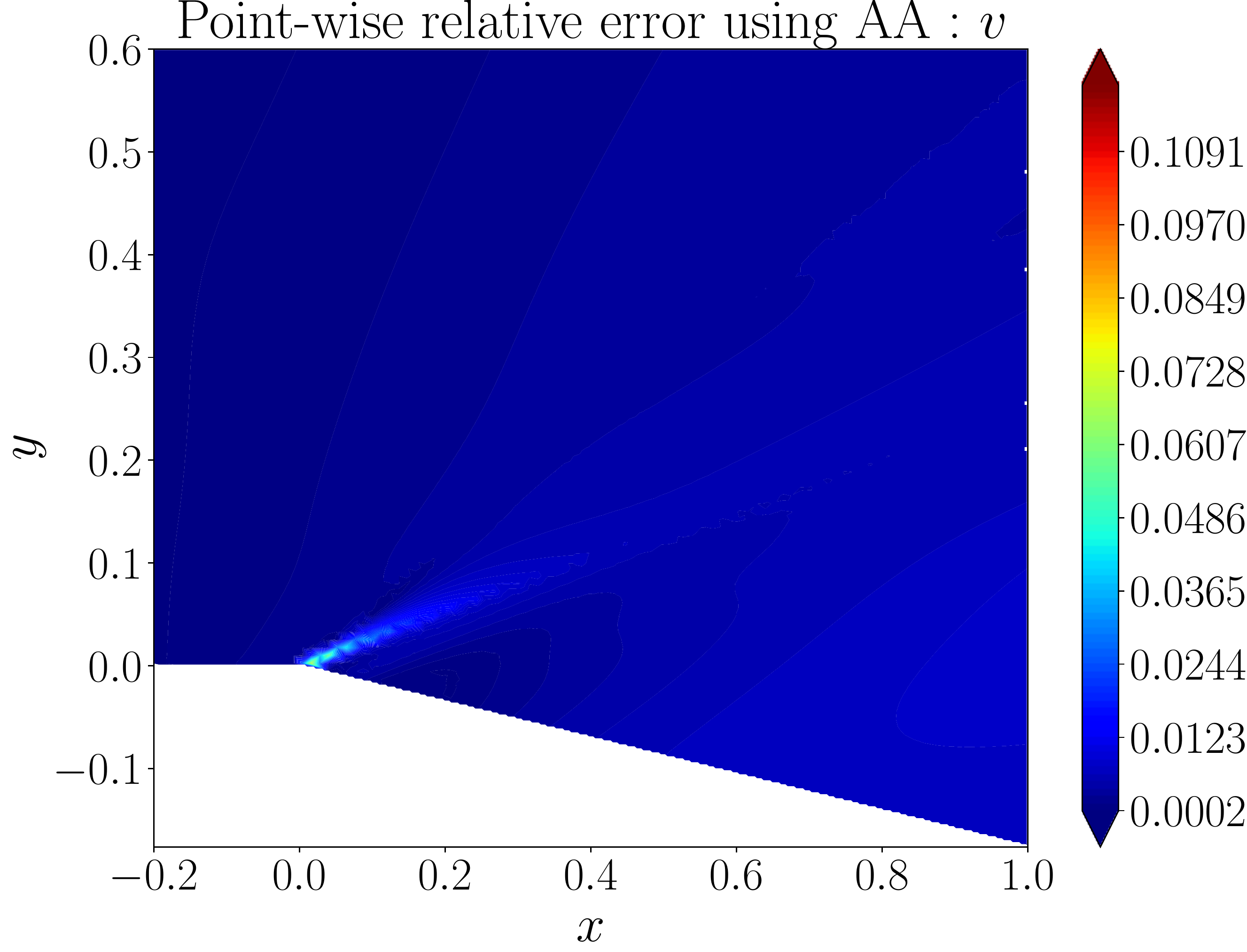}
    }
    \caption{Expansion wave problem: Comparison of the {\emph{velocity $v$}} between PINN solutions with/without adaptive activation function (AA) or dynamic weights (DW). Here we use `tanh' as the basic activation function. (a) Exact solution. (b) PINN solutions using dynamic weights. (c) Relative point-wise error with dynamic weights. (d) Relative point-wise error with adaptive activation function. (e) Relative point-wise error without dynamic weights or adaptive activation function. (f) Relative point-wise error with XPINN using adaptive activations.}
    \label{fig:expansion:v:compa}
\end{figure}

Besides using the typical hyperbolic tangent activation function and fixed weights, we additionally investigate the effectiveness of using the adaptive `tanh' activation function~\cite{jagtap2019adaptive}. For PINNs, we use 6 hidden-layers with 40 neurons in each layer for all cases.
The results are shown in Figs. \ref{fig:expansion:rho:compa}-\ref{fig:expansion:v:compa}. We observe that all the predictions of the density obtained with or without adaptive activation function/dynamic weights are in good agreement with the exact solution. However, for the predictions of the pressure and the velocity, the PINN solutions obtained without activation functions or dynamic weights are not as good as the ones obtained with activation functions or dynamic weights, especially for the velocity (see Figs. \ref{fig:expansion:u:compa}-\ref{fig:expansion:v:compa}). This means that we predict all the states by using  PINNs with adaptive activation function or dynamic weights along with the above setup, and the use of adaptive activate functions or dynamical weights helps to improve the accuracy of PINNs. This is very important when we consider the oblique shock or bow shock wave problems later. Additionally, we observe that the results obtained by using dynamical weights have supreme accuracy among all the results of PINNs. This is also verified by presenting the relative $L_2$ error (see Table \ref{tab:ExpRelErr}). Specifically, by using the above setup, we obtain good predictions with PINN without using the outflow boundary condition and the wall boundary condition.
Moreover, we verify the condition 
\begin{equation}\label{bc:nU}
    \tilde{n}\cdot {\bm u} = 0
\end{equation}
on the wall, where $\tilde{n}$ is the unit outer normal vector. This will be used as additional constraint for the bow shock problem considered later.

\begin{figure} 
\centering
\includegraphics[trim=0cm 0cm 0cm 0cm, clip=true, scale=0.5, angle = 0]{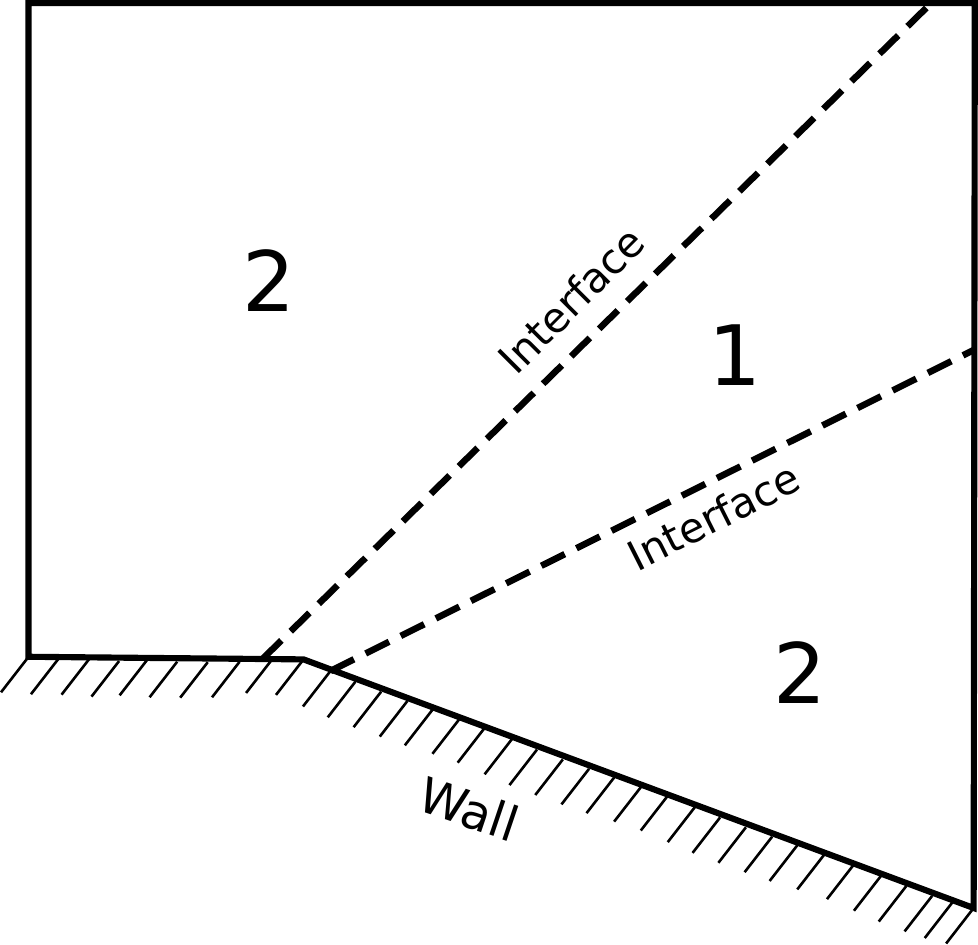}
\captionof{figure}{Expansion wave problem: Schematic presentation of subdomains and interfaces for the XPINN method.}
\label{fig:VB}
\end{figure}

\begin{figure}[http]
    \centering
    \subfigure{
    \includegraphics[trim=0cm 0cm 0cm 0cm, clip=true, scale=0.46, angle = 0]{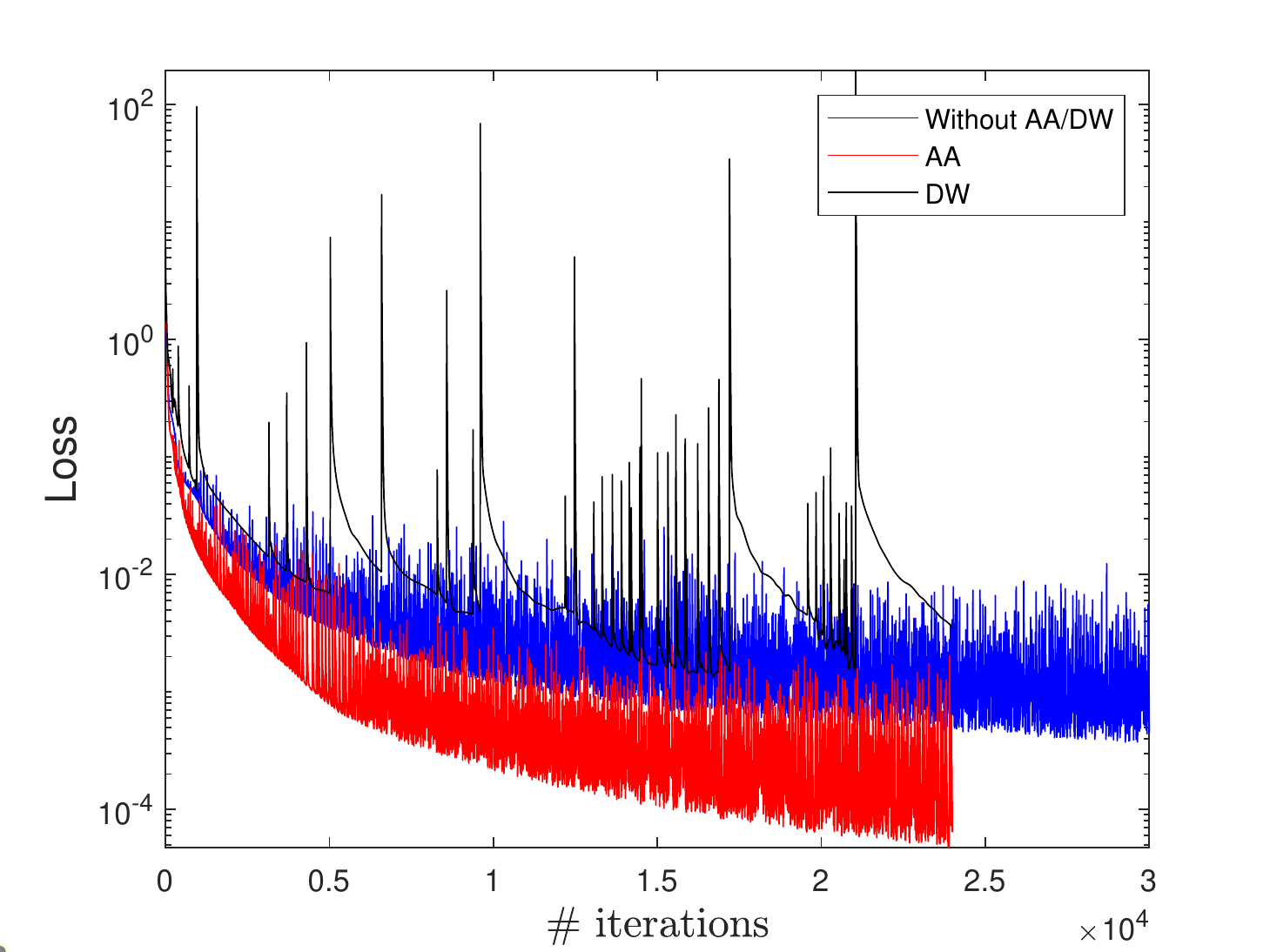}
    }
    \subfigure{
    \includegraphics[trim=0cm 0cm 0cm 0cm, clip=true, scale=0.46, angle = 0]{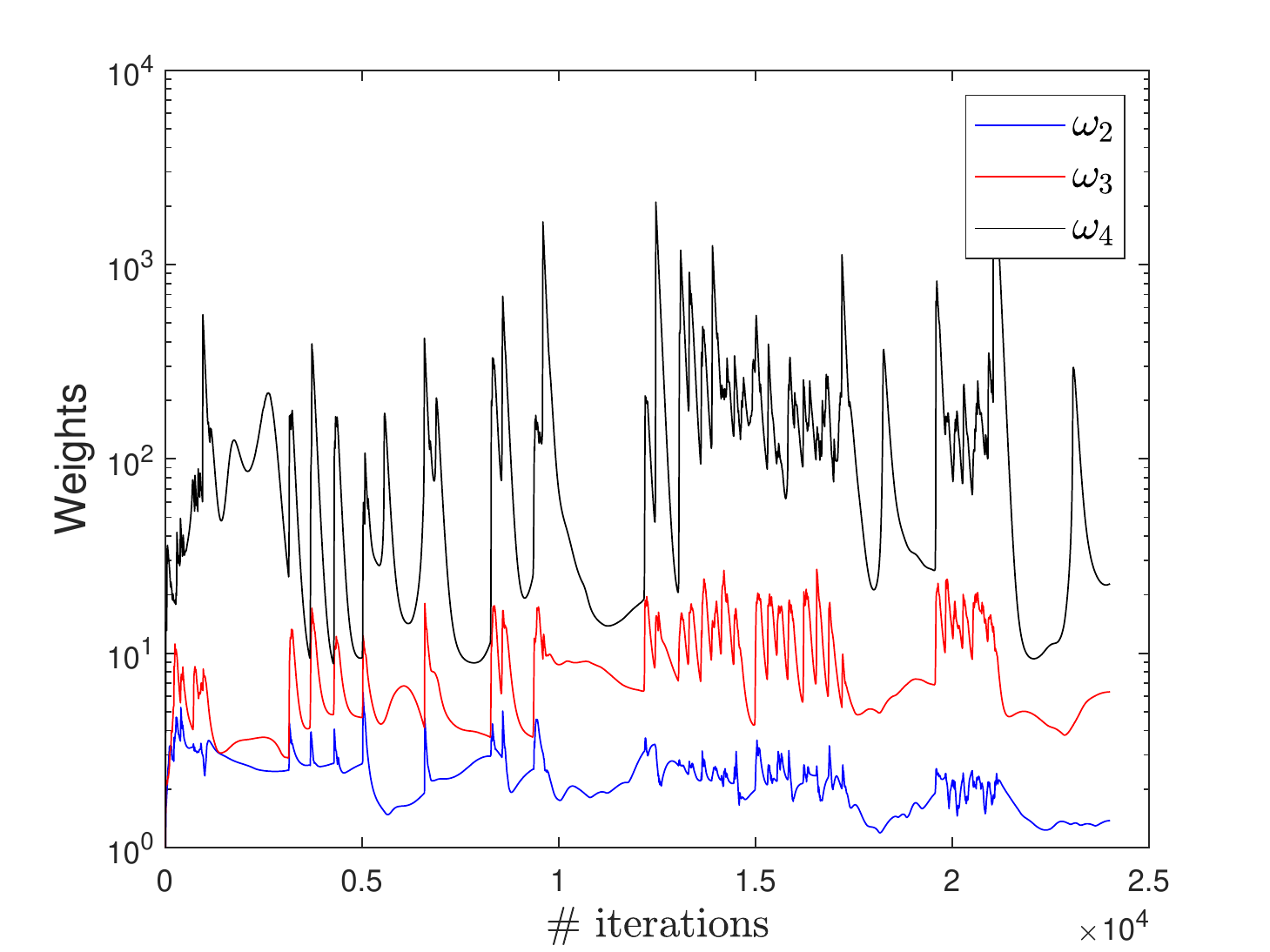}
    }
    \caption{Expansion wave problem: Loss function variation (left) and the variation of dynamic weights (right) using PINNs.}
    \label{fig:loss:expansion}
\end{figure}

\begin{figure}[http]
    \centering
    \subfigure{
    \includegraphics[trim=0cm 0cm 0cm 0cm, clip=true, scale=0.5, angle = 0]{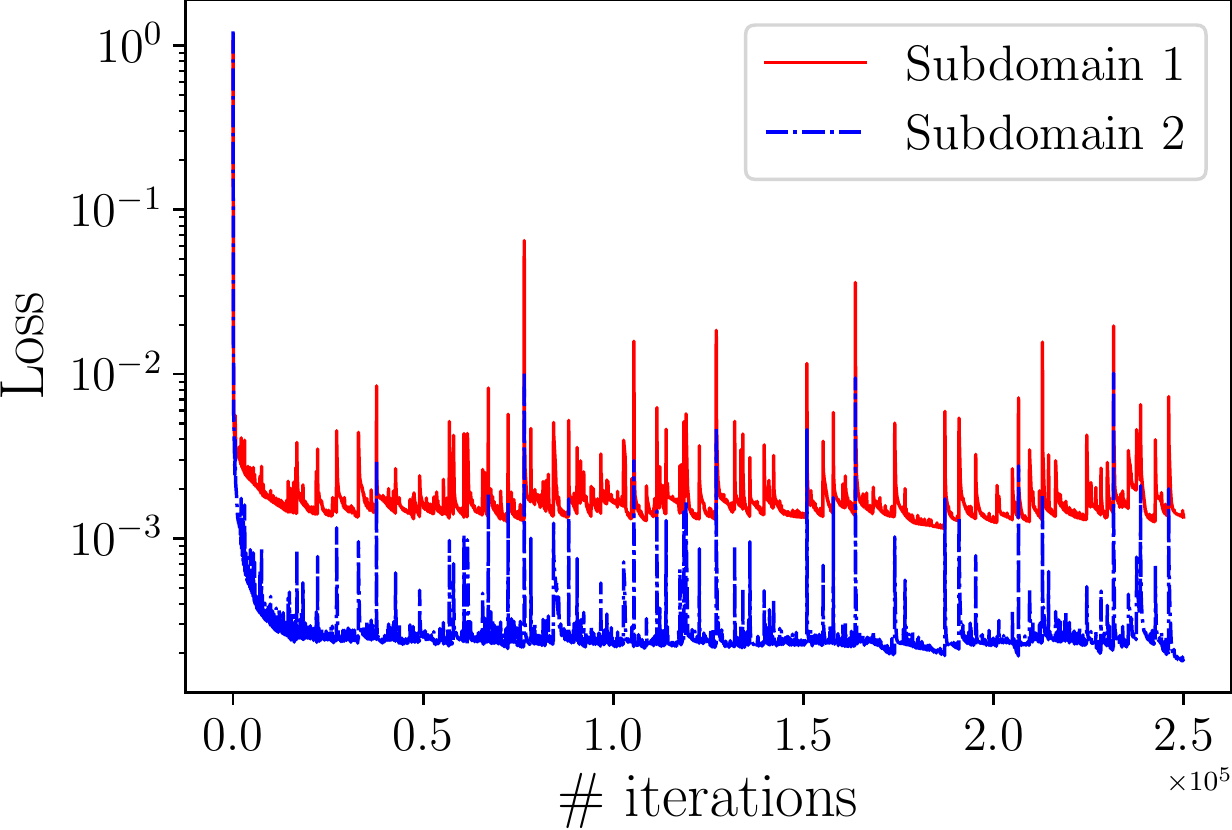}
    }
    \subfigure{
    \includegraphics[trim=0cm 0cm 0cm 0cm, clip=true, scale=0.47, angle = 0]{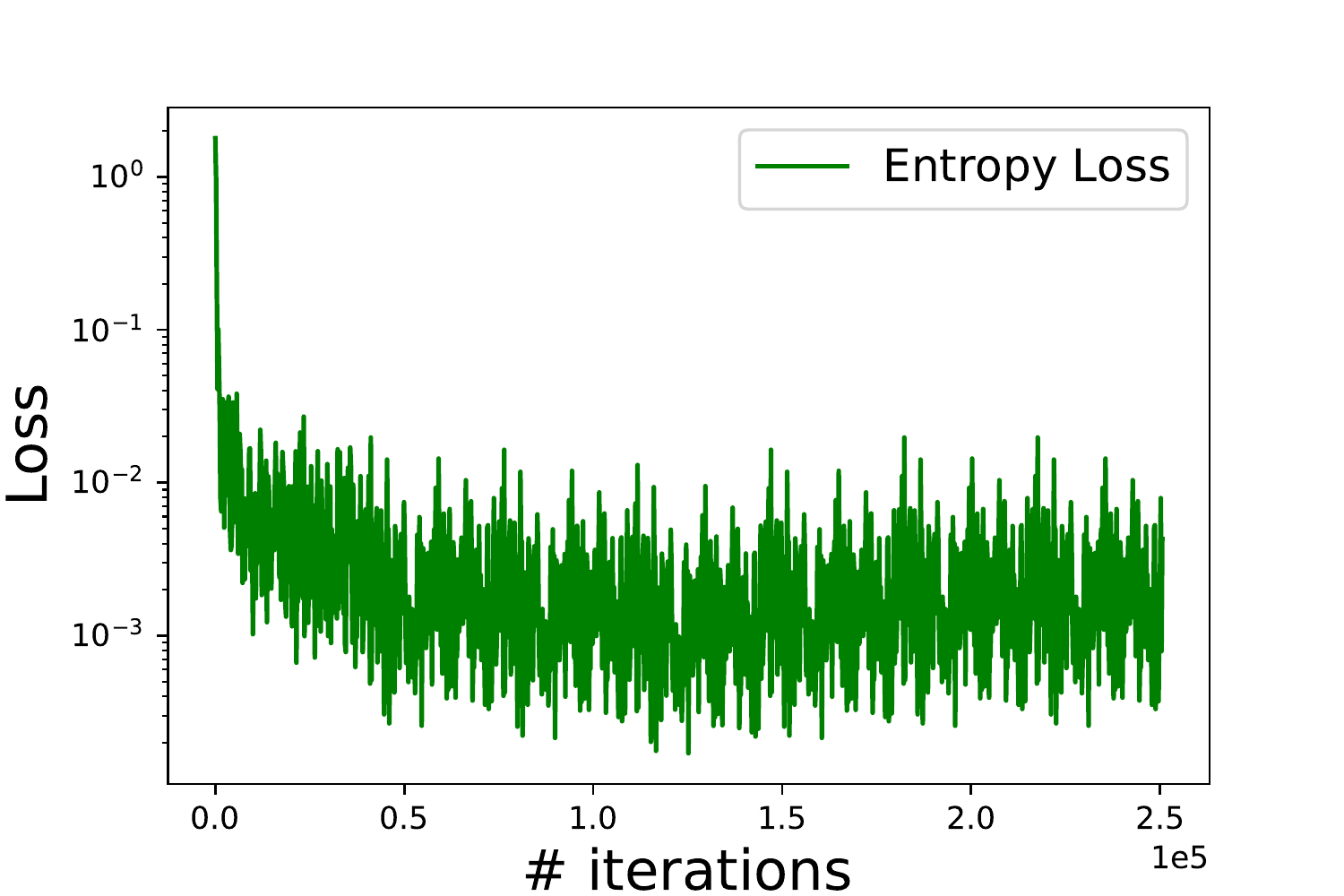}
    }
    \caption{Expansion wave problem: Loss function variation in two subdomains (left) and the combined entropy loss function for both subdomains (right) using XPINNs.}
    \label{fig:loss:expansion1}
\end{figure}

We also solved this problem using the XPINN methodology, and compared it with PINN results. For that, we divide the domain into two subdomains as shown in Fig. \ref{fig:VB}, keeping the network architecture same as PINNs in these subdomains. The number of evenly distributed interface points on both interfaces is 300.
The adaptive activation functions are used for the XPINN method. It can be seen that the XPINN results are slightly better than PINN. The  relative $L_2$ errors in all the primitive variables for all cases are given in Table \ref{tab:ExpRelErr}. In the case of XPINNs, we present the total relative error from all the subdomains.  We also plot the loss function versus iterations for PINNs; see Fig. \ref{fig:loss:expansion} (left). Fig. \ref{fig:loss:expansion}  (right plot) shows the evolution of the dynamic weights ($\omega_2, \omega_3, \omega_4$) with the number of iterations for PINNs. Fig.  \ref{fig:loss:expansion1} shows the total loss function plot for XPINNs (left) and the combined entropy loss (right) from all the subdomains.

\begin{table}[]
\centering
\begin{tabular}{|c|c|c|c|c|}
\hline
 & PINNs without AA/DW & PINNs with AA & PINNs with DW  & XPINNs with AA \\ \hline
$\rho$ & 1.3494e-02 & 6.5783e-03& 4.0532e-03& 3.3864e-03\\ 
$u$  &  7.4464e-02 & 2.7493e-02& 1.5532e-02& 1.1842e-02 \\ 
$v$  &  6.9054e-02 & 2.8942e-02& 1.6578e-02& 1.3864e-02 \\
$p$  &  2.2646e-02 & 5.5783e-03& 4.3532e-03& 3.0864e-03 \\
\hline
\end{tabular}
\caption{Expansion wave problem: Relative $L_2$ errors in all the primitive variables for all cases. For XPINNs, we present the total relative error from all the subdomains.}
\label{tab:ExpRelErr}
\end{table}

Next, we show the theoretical generalization bound in XPINN and PINN models given in \cite{hu2021extended}.
\begin{table}[]
\centering
\begin{tabular}{|c|c|c|c|c|}
\hline
Method   & Rel. $L_2$ error in $\rho$& Norms & Theoretical Bounds\\ \hline
PINN using AA  &  6.5783e-03 &100.00\% & 100.00\% \\ \hline
XPINN-1  & \multirow{2}{*}{3.3864e-03}  & 125.31\% & \multirow{2}{*}{83.23\%} \\ \cline{1-1} \cline{3-3}
XPINN-2  & & 32.59\% & \\ \hline
\end{tabular}
\caption{The Rademacher complexity and the theoretical bounds for the Expansion wave problem. These are computed using the relative $L_2$ error in density.}
\label{tab:Exp}
\end{table}
These results are calculated using the relative $L_2$ error in density. We present the results including the relative $L_2$ error in density, as well as the calculated theoretical generalization bound in XPINN and PINN models in Table \ref{tab:Exp}. The results show that XPINN generalizes better than PINN (see the theoretical bound column).
In addition, the Rademacher complexity is presented in the column `Norms', where that of PINN is set 100\% for comparison. The complexity of PINN is the largest in subdomain 1, since PINN has to fit the entire target function into one network, which is very complex, whereas in subdomain 2, the XPINN complexity is slightly larger due to the reduced training data set.

\subsection{Oblique shock problem}
\begin{figure}[http]
    \centering
    \includegraphics[trim=0cm 0cm 0cm 0cm, clip=true, scale=0.145, angle = 0]{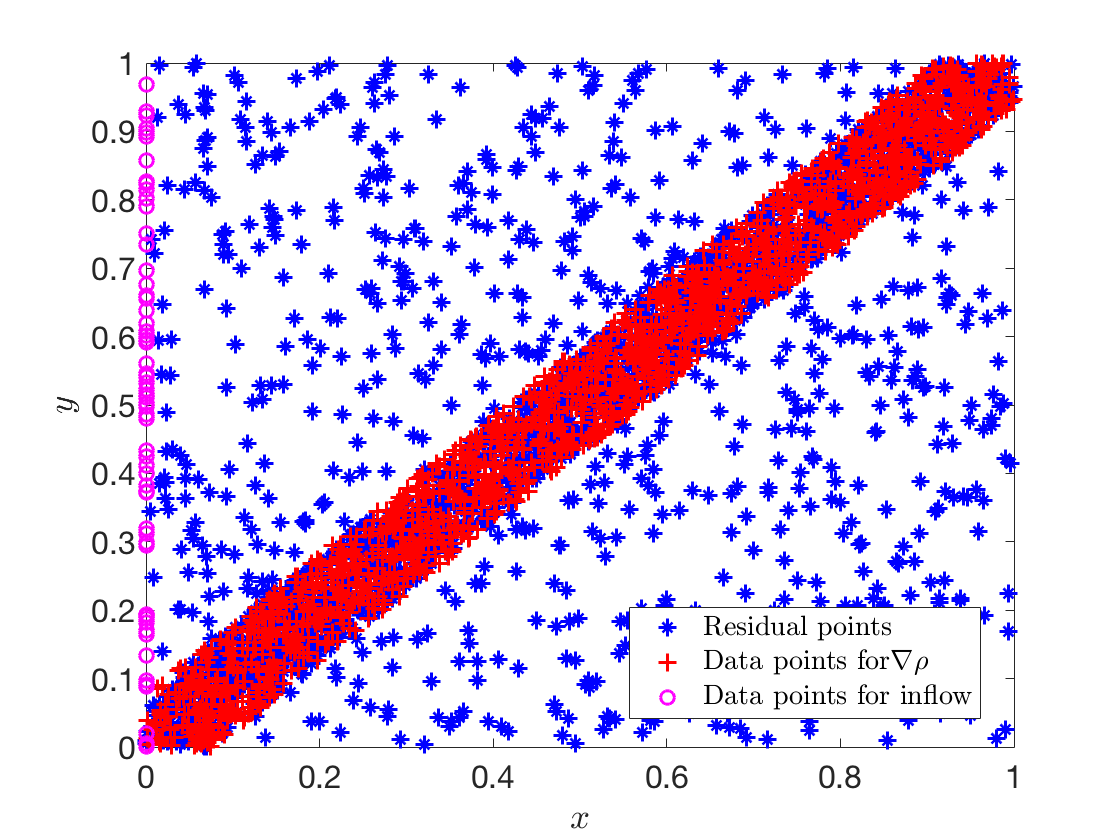}\,\,\,\,\,\,\,\,\,\ \ \ \ \ \ \ \ \ \ \ \ \ \ \ \ 
\includegraphics[trim=0cm 0cm 0cm 0cm, clip=true, scale=0.58, angle = 0]{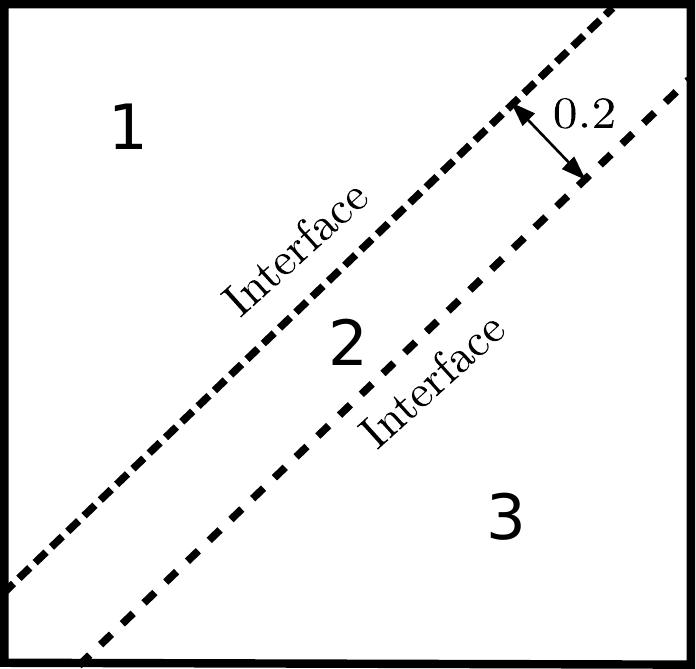}
    \caption{Oblique shock problem: (Left) distributions of the 3900 residual points (blue $\ast$ points), the 1200 data points for the density gradient (red $+$ points), and the 120 data points for the inflow conditions (magenta $\circ$ points). (Right) schematic presentation of subdomains and interfaces for the XPINN method.}
    \label{fig:data:Oblique}
\end{figure}

We now consider another type of attached shock wave problem, i.e., the oblique shock problem (see {regions \textcircled{1} and \textcircled{4}} of the left plot of Fig. \ref{fig:schematic:shock}).
Let $\gamma = 1.4$, consider the oblique shock problem with the exact solution given by 
\begin{equation}
(\rho, u,v, p) =  \left\{ {\begin{array}{cc}
(0.06688~kg/m^3,   738.2~m/s, 0,~ 9485~Pa),  ~~~~~~~ \text{ pre-shock,}\\
(0.09515~kg/m^3, 635.9\sin\theta~ m/s, 635.9\cos\theta~ m/s, 1.5\times 10^4~ Pa),  \text{ post-shock,}\\
\end{array} } \right.
\end{equation}
where $\theta = 10^\circ$.
The loss function of PINN is given as follows:
\begin{equation}\label{loss:oblique}
\begin{aligned}
   \mathcal{J} = & ~\omega_{1}~\text{MSE}_{\mathcal{F}} + \omega_{2}~\text{MSE}_{\nabla \rho|_{D}} + \omega_{3}~\text{MSE}_{Inflow} + \omega_{4}~\text{MSE}_{p^*} \\
    &+ \omega_{5}~(\text{MSE}_{Mass} + \text{MSE}_{Momentum} + \text{MSE}_{Energy}),
\end{aligned}
\end{equation}
where $\text{MSE}_{\mathcal{F}}$ corresponds to the Euler equation and the entropy condition, $\text{MSE}_{\nabla \rho|_{D}}$ corresponds to the density gradient and $D\subset \Omega: = [0,1]^2$ is a region shown in the Fig. \ref{fig:data:Oblique} (left) (domain with red $+$ points). Note that due to discontinuous solution we employed finite different formula to enforced the density gradient; see \cite{mao2020physics} for more details. $\text{MSE}_{Inflow}$ corresponds to the inflow boundary condition along the inflow boundary points shown by magenta circles in Fig. \ref{fig:data:Oblique} (left), $\text{MSE}_{p^*}$ corresponds to the pressure values at the point $(0.4,0)$, and $\text{MSE}_{Mass}$, $\text{MSE}_{Momentum}$ and $\text{MSE}_{Energy}$ correspond to the total mass, the total momentum and the total energy, respectively. For PINNs, we use 7 hidden-layers with 30 neurons in each layer for all cases.
In the case of XPINN, the domain is divided into 3 subdomains, and a separate neural network with same network architecture as PINNs is employed in each subdomain. The Fig. \ref{fig:data:Oblique} (right) shows the schematic of three XPINN subdomains along with common interfaces. For stitching these subdomains, the conservative flux continuity condition \cite{jagtap2020conservative} is also used along with residual continuity condition. Also, as discussed in \cite{jagtap2020conservative, jagtap2020extended}, an average solution is enforced along the common interface for faster convergence. The number of evenly distributed interface points on both interfaces is 200.

\begin{figure}[htpb] 
\centering
\includegraphics[trim=1.5cm 1cm 2cm 0cm, clip=true, scale=0.35, angle = 0]{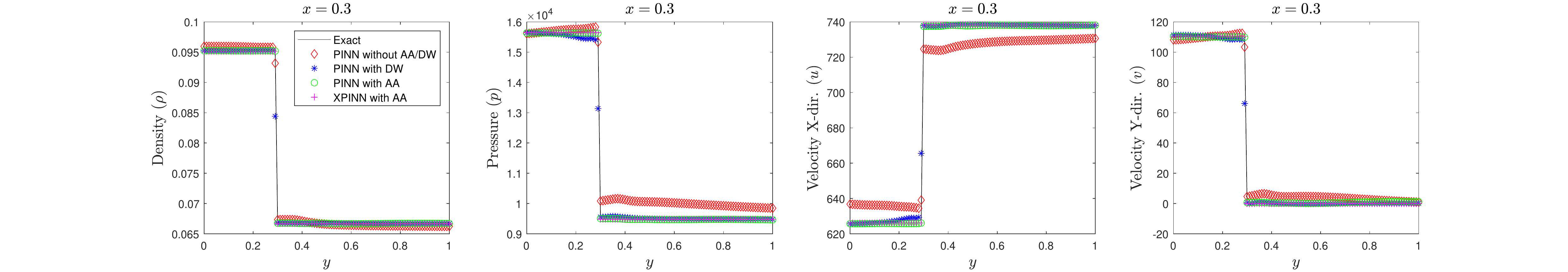}\\
\includegraphics[trim=1.5cm 1cm 2cm 0cm, clip=true, scale=0.35, angle = 0]{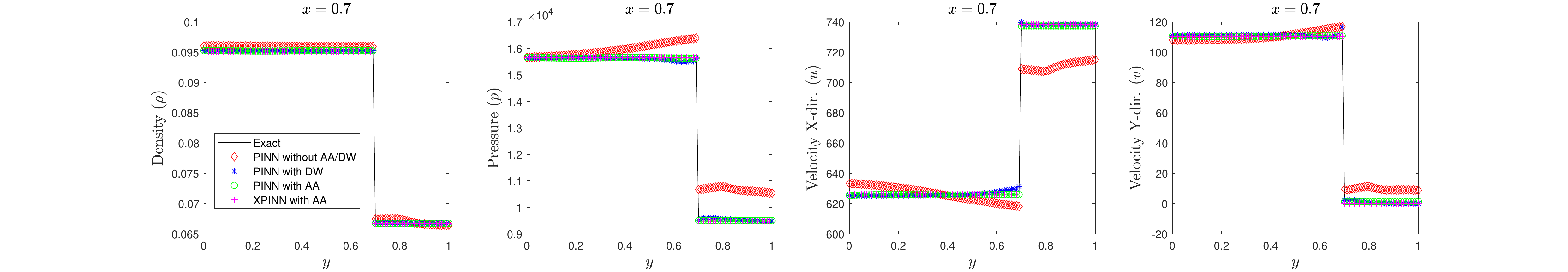}
\caption{Oblique shock problem: Primitive variables predicted by PINNs and XPINNs at $x = 0.3$ (top row) and $x = 0.7$ (bottom row).}
\label{fig:XPINNsOSW1}
\end{figure}

\begin{figure}[htpb] 
\centering
\includegraphics[trim=1.5cm 1cm 2cm 0cm, clip=true, scale=0.35, angle = 0]{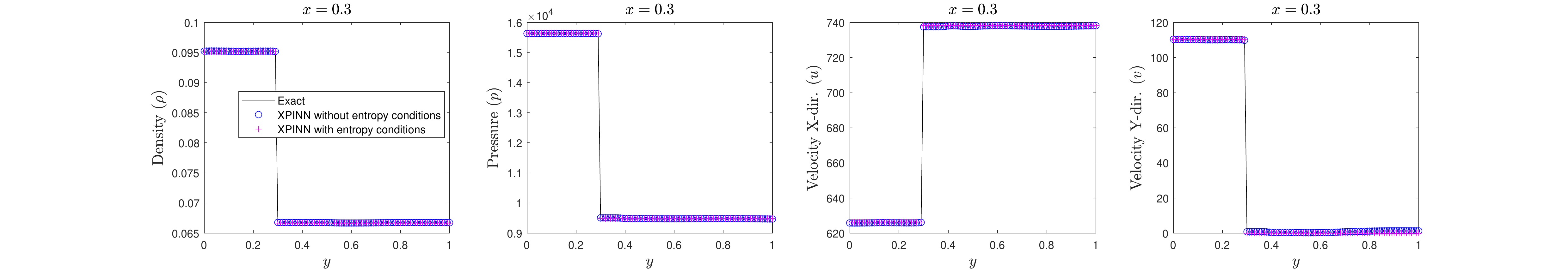}\\
\includegraphics[trim=1.5cm 1cm 2cm 0cm, clip=true, scale=0.35, angle = 0]{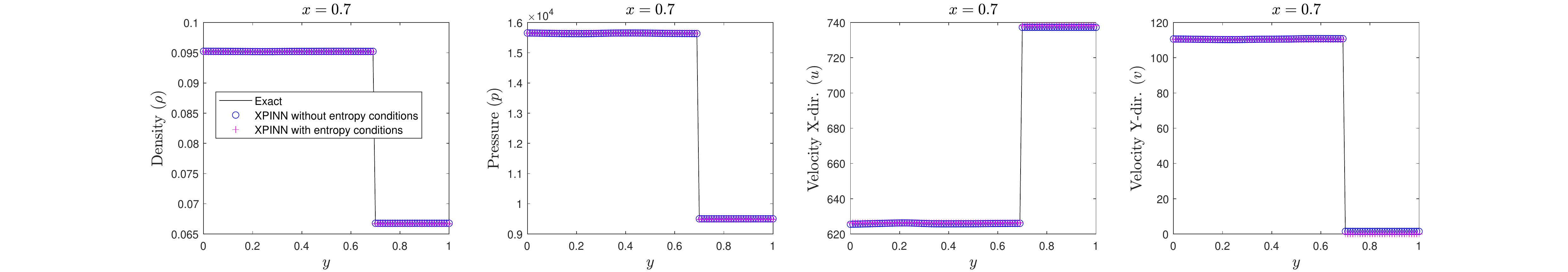}
\caption{Oblique shock problem: Primitive variables predicted by XPINNs with and without entropy condition at $x = 0.3$ (top row) and $x = 0.7$ (bottom row).}
\label{fig:XPINNsOSW2}
\end{figure}

The results are shown in Fig. \ref{fig:XPINNsOSW1}  at locations $x = 0.3$ and 0.7. The PINN is used with or without adaptive activations and dynamic weights, whereas XPINN is employed with adaptive activations. 
We note that the vanilla PINN solutions gives large errors. This implies that the vanilla PINNs cannot deal well with the problems with discontinuous solutions. On the other hand,
it can be observed that the predicted density, velocity and pressure match well with the exact solution using adaptive activations for both PINNs and XPINNs as well as the case of using dynamical weights for PINNs. The relative $L_2$ errors in all the primitive variabls for all the cases are given in Table \ref{tab:OSWRelErr}. In the case of XPINNs, the total relative error in all the subdomains is given together. This means that the use of adaptive activations or dynamical weights significantly improves the accuracy of the predictions. We also compare the results of XPINN with and without entropy condition as shown in Fig. \ref{fig:XPINNsOSW2}. In both cases, the oblique shock wave is captured exactly.
\begin{table}[]
\centering
\begin{tabular}{|c|c|c|c|c|}
\hline
 & PINNs without AA/DW & PINNs with AA & PINNs with DW  & XPINNs with AA \\ \hline
$\rho$ & 3.5373e-02 & 6.5533e-03& 7.1642e-03& 5.3864e-03\\ 
$u$  & 9.9734e-02 & 7.7493e-03& 1.8442e-02& 8.1842e-03 \\ 
$v$  &  8.1637e-02 & 7.8424e-03& 9.4577e-03& 7.3864e-03 \\
$p$  &  6.8367e-02 & 4.5363e-03& 6.7382e-03& 4.0864e-03 \\
\hline
\end{tabular}
\caption{Oblique shock problem: Relative $L_2$ errors in all the primitive variables for all cases. In the case of XPINNs, the total relative error from all the subdomains is given together.}
\label{tab:OSWRelErr}
\end{table}
\begin{figure}[http]
    \centering
    \subfigure{
    \includegraphics[trim=0cm 0cm 0cm 0cm, clip=true, scale=0.34, angle = 0]{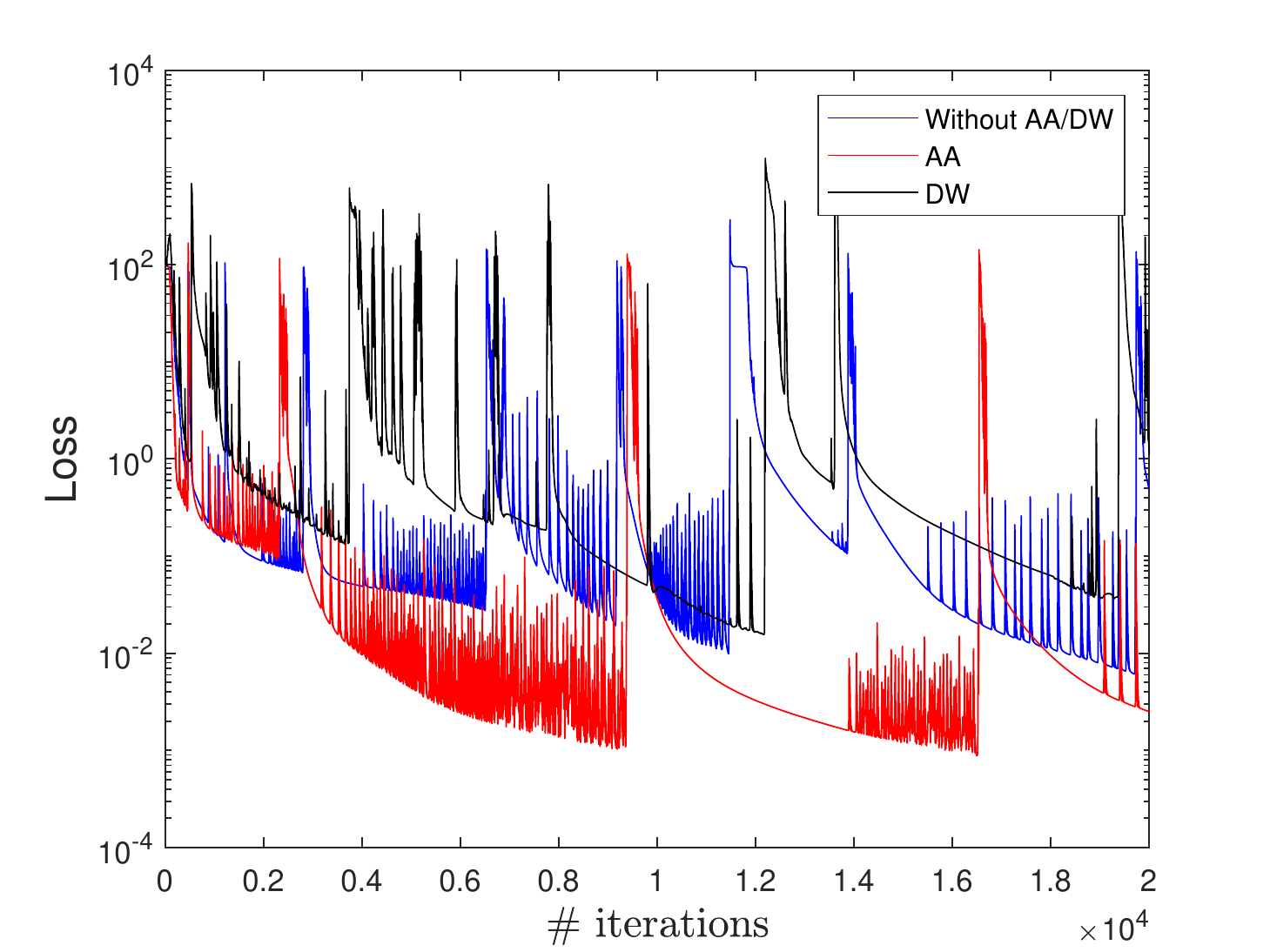}
    }
    \subfigure{
    \includegraphics[trim=0cm 0cm 0cm 0cm, clip=true, scale=0.34, angle = 0]{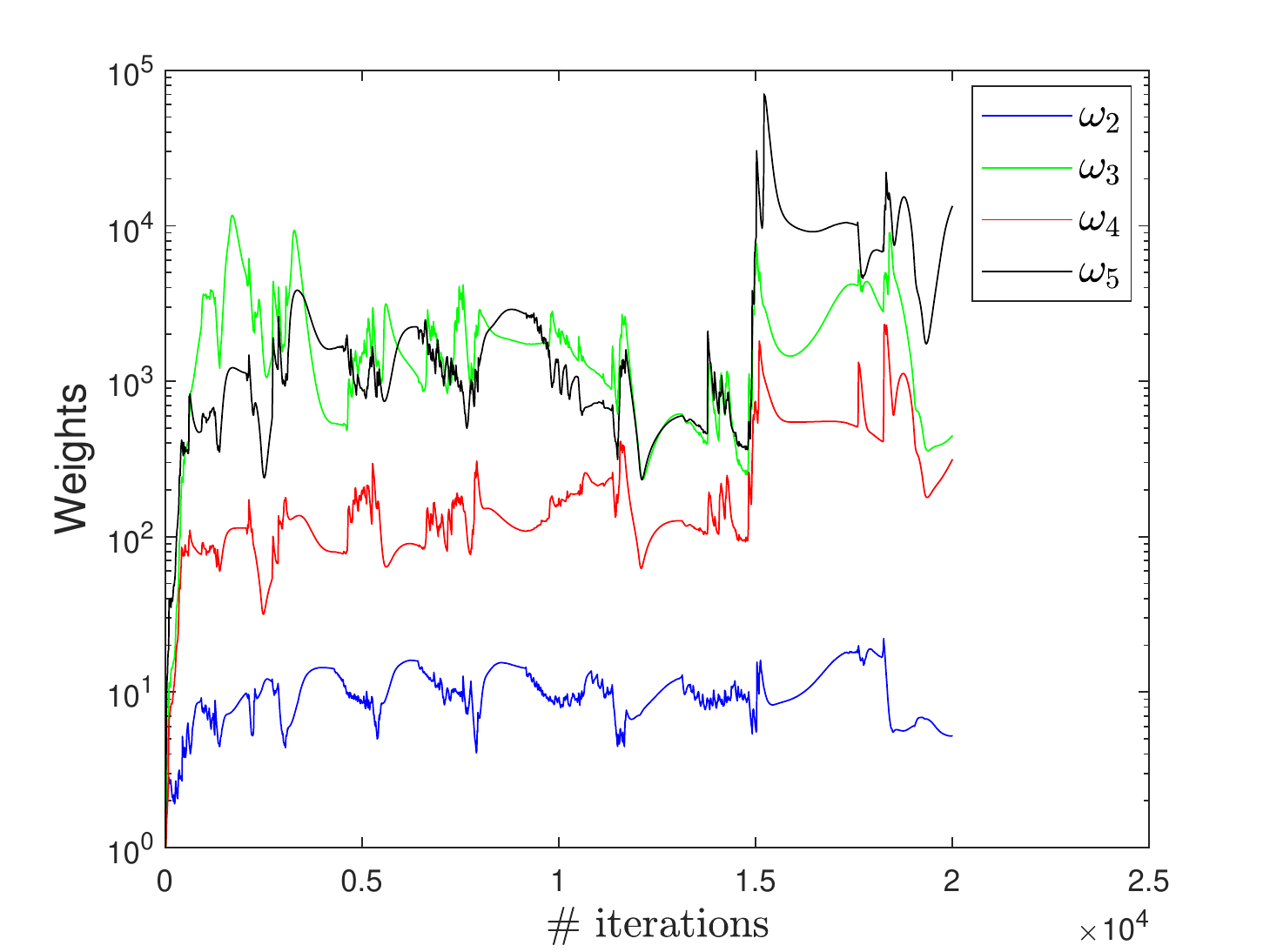}
    }
    \subfigure{
    \includegraphics[trim=0cm 0cm 0cm 0cm, clip=true, scale=0.43, angle = 0]{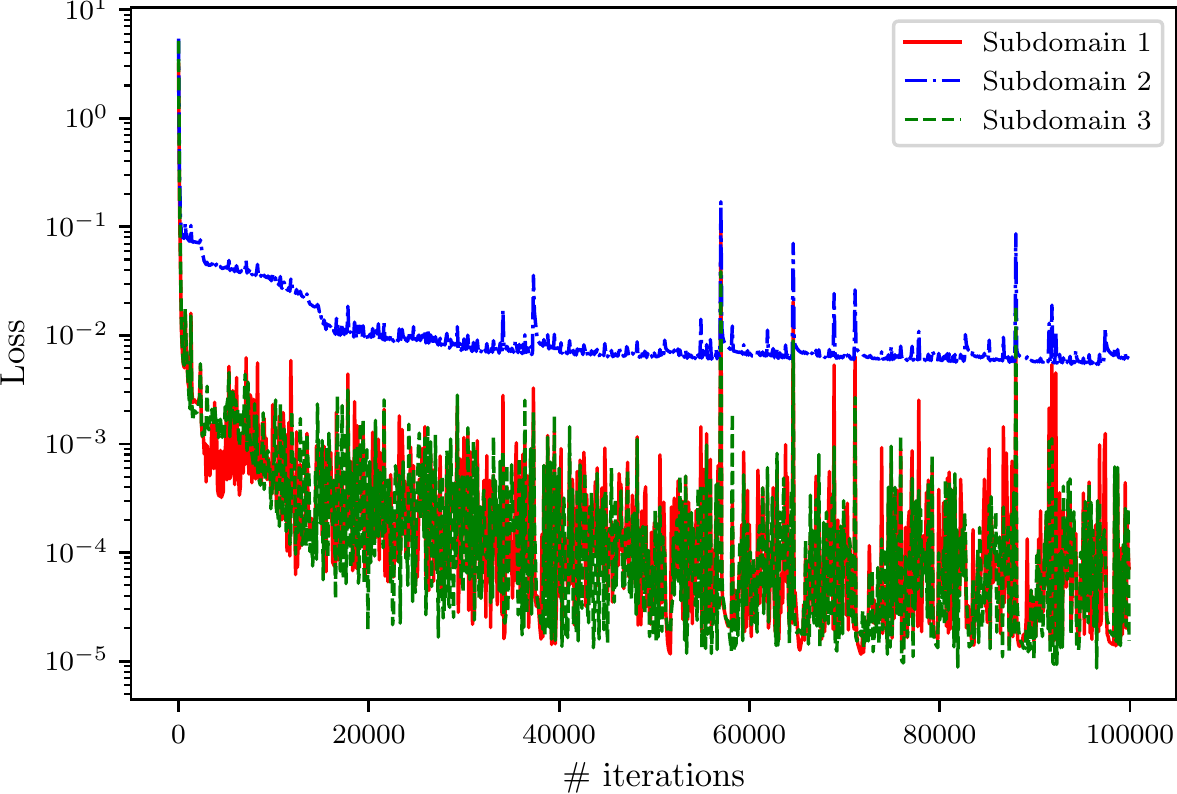}
    }
    \caption{Oblique shock problem: Loss versus iterations for PINNs and XPINNs given by the left and right figures, respectively. Middle figure shows the variation of PINN dynamic weights ($\omega_2$- $\omega_5$) with the number of iterations.}
    \label{fig:loss:oblique}
\end{figure}
Fig. \ref{fig:loss:oblique} shows the loss function variations with resect to  iterations for PINNs and XPINNs given by the left and right figures, respectively. The middle figure shows the variation of the PINN dynamic weights ($\omega_2$- $\omega_5$) with the number of iterations.

Now, we show the theoretical generalization bound in XPINN and PINN models for the oblique shock wave problem.
\begin{table}[]
\centering
\begin{tabular}{|c|c|c|c|}
\hline
Method   & Rel. $L_2$ error in $\rho$ & Norms & Theoretical Bounds\\ \hline
PINN using AA  & 6.5533e-03& 100.00\% & 100.00\% \\ \hline
XPINN-1 & \multirow{3}{*}{5.3864e-03}  & 41.34\% & \multirow{3}{*}{89.43\%} \\ \cline{1-1} \cline{3-3}
XPINN-2 &  & 121.69\% & \\ \cline{1-1} \cline{3-3}
XPINN-3 &  & 49.72\% & \\ \hline
\end{tabular}
\caption{The Rademacher complexity and the theoretical bounds for the oblique shock wave problems. These are computed using the relative $L_2$ error in density.}
\label{tab:Os}
\end{table}
Table \ref{tab:Os} shows the results including relative $L_2$ error in density and the calculated theoretical generalization bound in XPINN (for subdomains 1,2, and 3) and PINN models. Again, XPINN generalizes better than PINN (see the theoretical bound).  After observing the Rademacher complexity given in the column `Norms', the complexity of PINN is higher compared to XPINNs in subdomains 1 and 3, since PINN has to fit the entire discontinuous target function into one network. In subdomain 2, the complexity of XPINN is slightly higher than PINN due to the presence of oblique shock as well as the reduced training data set.

\subsection{Detached bow shock problem}
\begin{figure}[http]
    \centering
    \includegraphics[scale = 0.42]{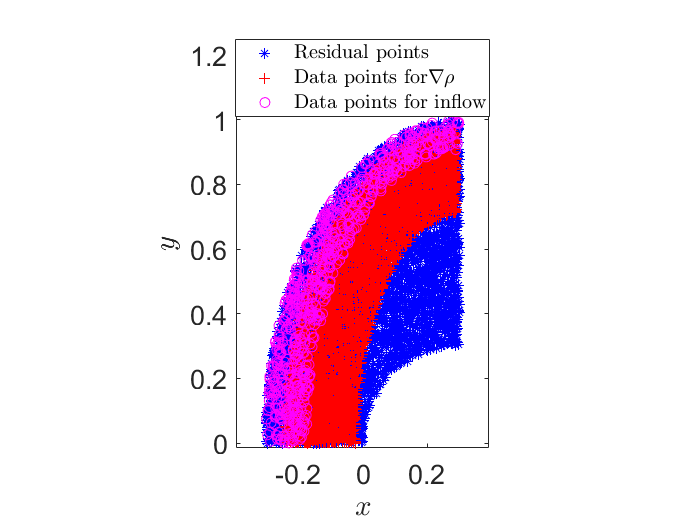}
    \includegraphics[scale = 0.58]{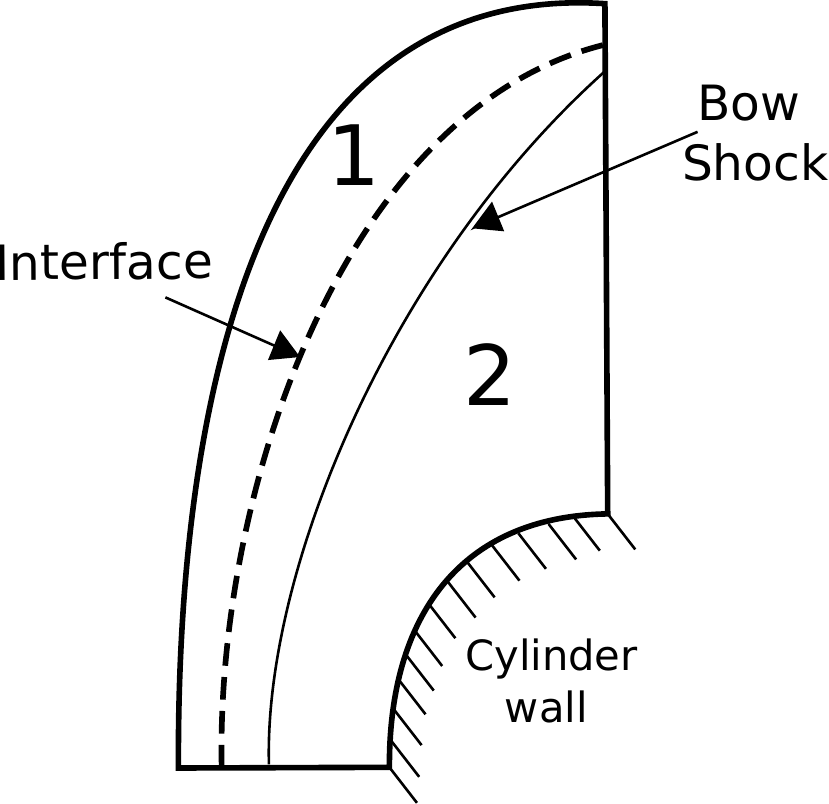}
    \caption{Bow shock problem: (Left) distributions of the 5600 residual points (blue $\ast$ points), the 700 data points for the density gradient (red $+$ points) and the 300 data points for the inflow conditions (magenta $\circ$ points). (Right) schematic presentation of subdomains and curved interface for the XPINN method.}
    \label{fig:distribution:bow:h}
\end{figure}

We now consider a detached shock problem, i.e., the bow shock problem; see Fig. \ref{fig:schematic:shock} (b).
Again, we would like to solve the inverse problem of the bow shock problem by using PINNs and XPINNs.
We now consider the bow shock problems with the following inlet flow conditions 
    $$M_\infty = 4,\; p_\infty = 101253.6Pa,\; \rho_\infty = 1.225kg/m^3,\; u_\infty = 1360.6963m/s,\; v_\infty = 0,\; T_\infty = 288K.$$
The distributions of the residual points, the data points for the density gradient and the inflow conditions are shown in Fig. \ref{fig:distribution:bow:h} (left). The data points for the pressure are located on the surface of the body.
The CFD data were obtained by using the `Openbowshock' code~\cite{Karpuk2021}.

In addition to the loss function \eqref{loss:oblique}, we also enforce the condition \eqref{bc:nU} along the surface of the cylinder, namely, we have the following loss function

\begin{equation}\label{loss:drho:p:inlet:U}
\begin{aligned}
   \mathcal{J} = &~ \omega_{1}~\text{MSE}_{\mathcal{F}} + \omega_{2}~\text{MSE}_{\nabla \rho|_{D}} + \omega_{3}~\text{MSE}_{Inflow} + \omega_{4}~\text{MSE}_{p^*} \\
    &+ \omega_{5}~(\text{MSE}_{Mass} + \text{MSE}_{Momentum} + \text{MSE}_{Energy}) + \omega_{6}~\text{MSE}_{\tilde{n}\cdot {\bm u}},
\end{aligned}
\end{equation}
where $\text{MSE}_{\tilde{n}\cdot {\bm u}} = \frac{1}{N_{\bm u}} \sum_{k = 1}^{N_{\bm u}} (\tilde{n}\cdot {\bm u})^2(x_k, y_k)$ and $(x_k, y_k),\; k= 1,\ldots, N_{\bm u}$ are points located on the surface of the body.
Moreover, according to previous investigation, PINNs with dynamic weights and adaptive activation functions usually perform better than the vanilla PINNs. Therefore, in this case, we  use PINNs with dynamic weights. Further, we used 5 hidden-layers with 160 neurons in each layer. We also employ XPINN and compare the results with the PINN results. To this end, we divide the domain into two subdomains as shown in figure \ref{fig:distribution:bow:h} (right) keeping the network architecture same in these two subdomains. The number of points on the interface is 200, which are uniformly distributed along the interface line.

\begin{figure}[http]
    \centering
        \subfigure{
    \includegraphics[trim=0cm 0cm 0cm 0cm, clip=true, scale=0.16, angle = 0]{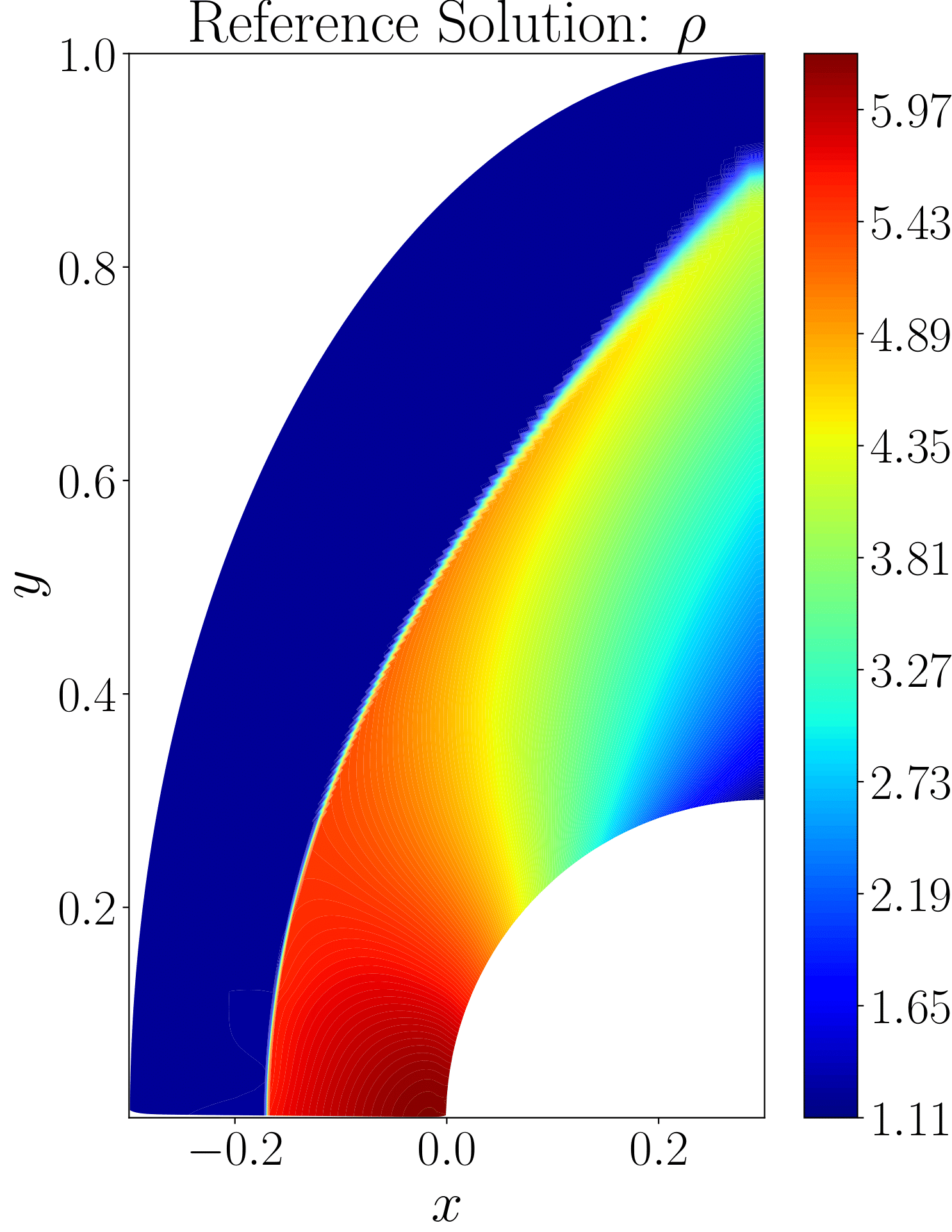}
    }
    \subfigure{
    \includegraphics[trim=0cm 0cm 0cm 0cm, clip=true, scale=0.16, angle = 0]{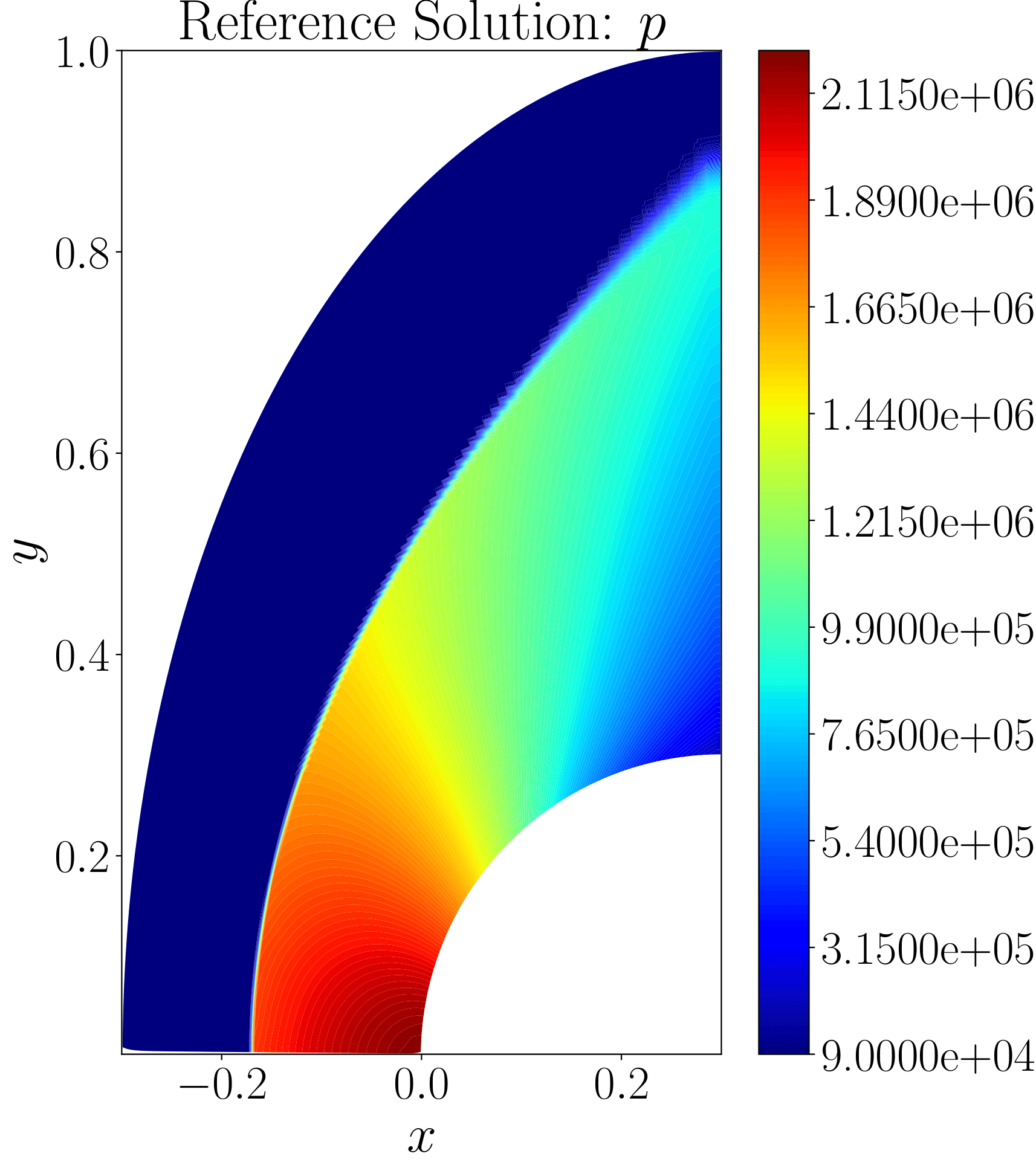}
    }
    \subfigure{
    \includegraphics[trim=0cm 0cm 0cm 0cm, clip=true, scale=0.16, angle = 0]{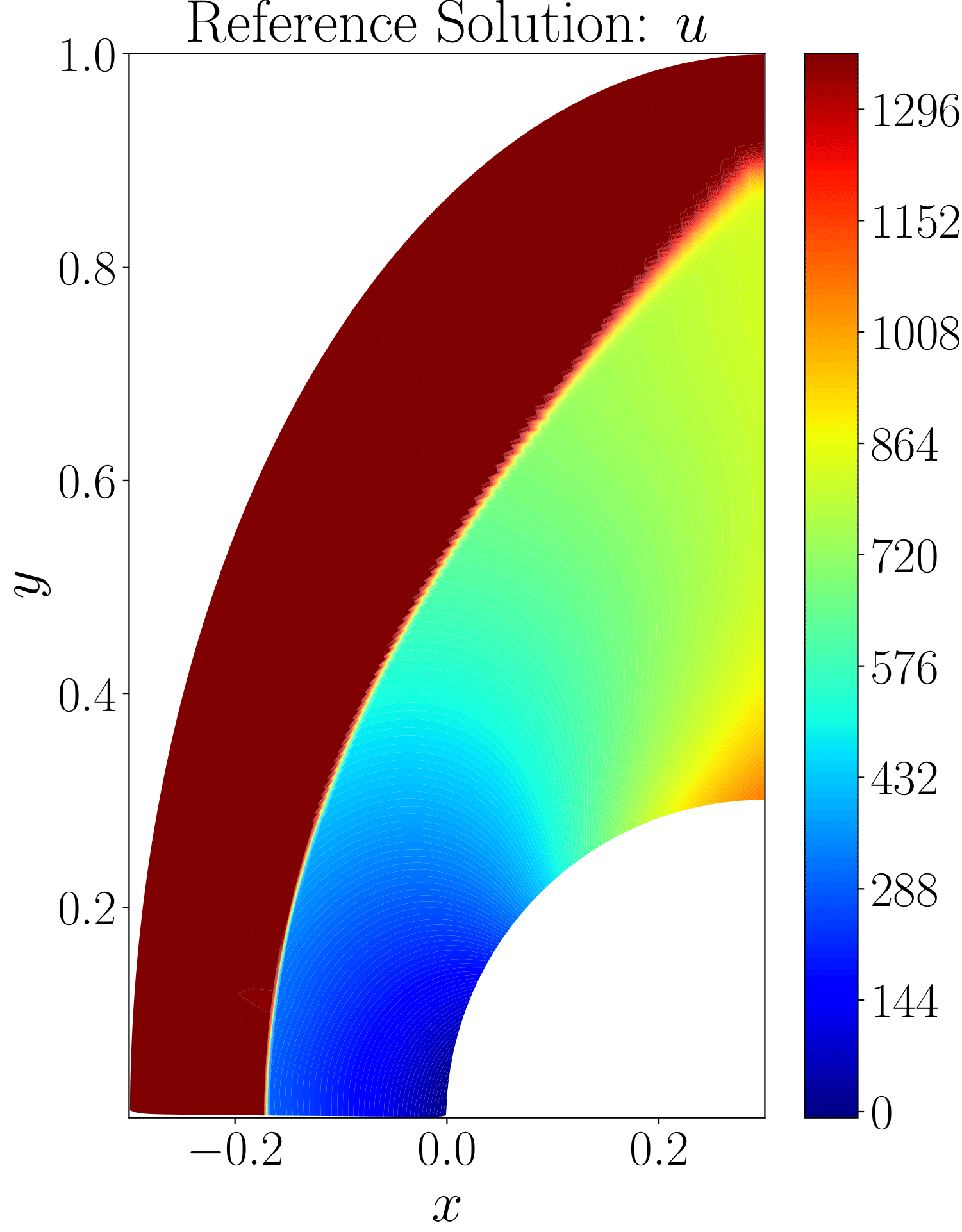}
    }
    \subfigure{
    \includegraphics[trim=0cm 0cm 0cm 0cm, clip=true, scale=0.16, angle = 0]{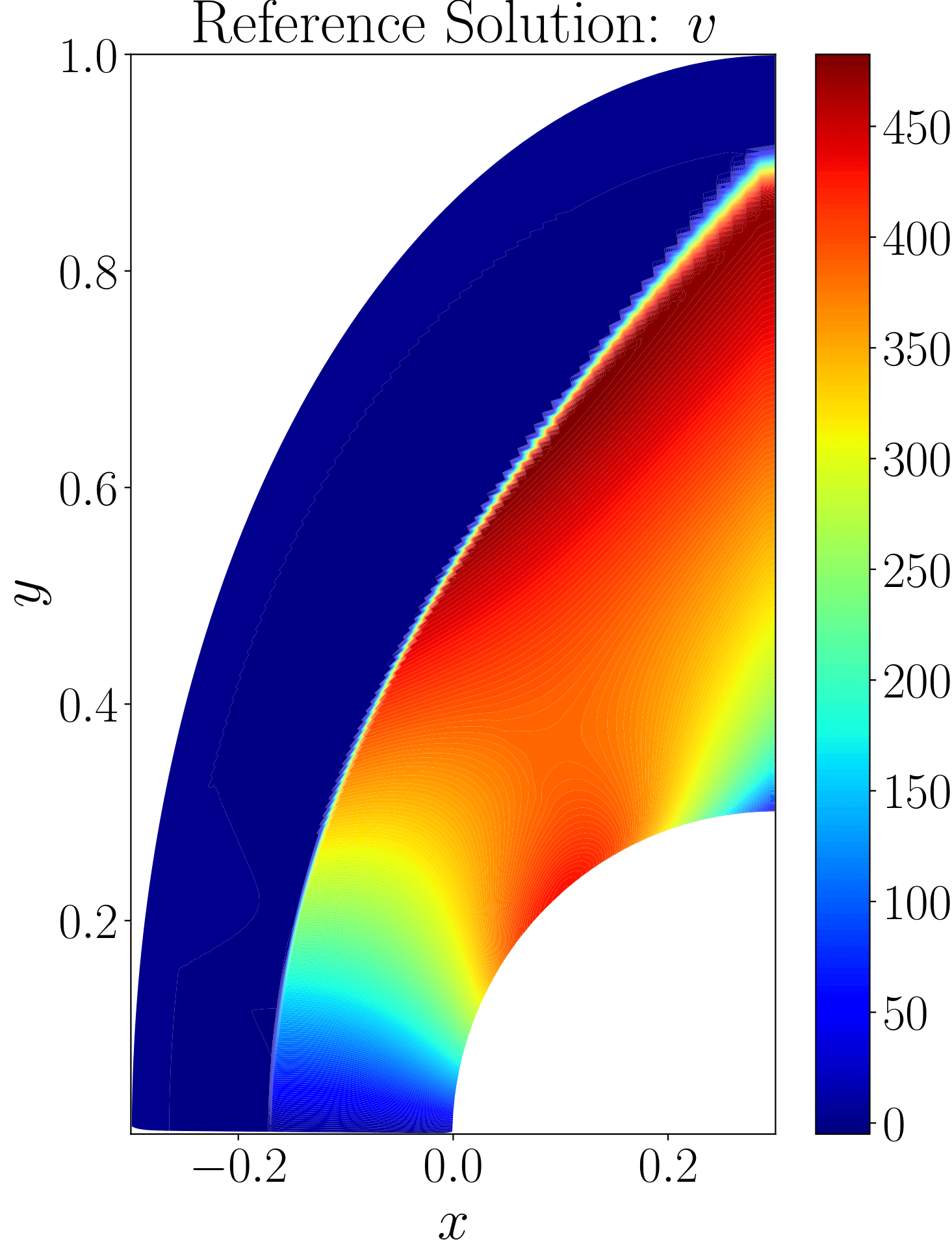}
    }
    \caption{Bow shock problem: Reference solutions for all primitive variables.}
    \label{fig:distribution:bow:exact}
\end{figure}

\begin{figure}[http]
    \centering
    \subfigure{
    \includegraphics[trim=0cm 0cm 0cm 0cm, clip=true, scale=0.16, angle = 0]{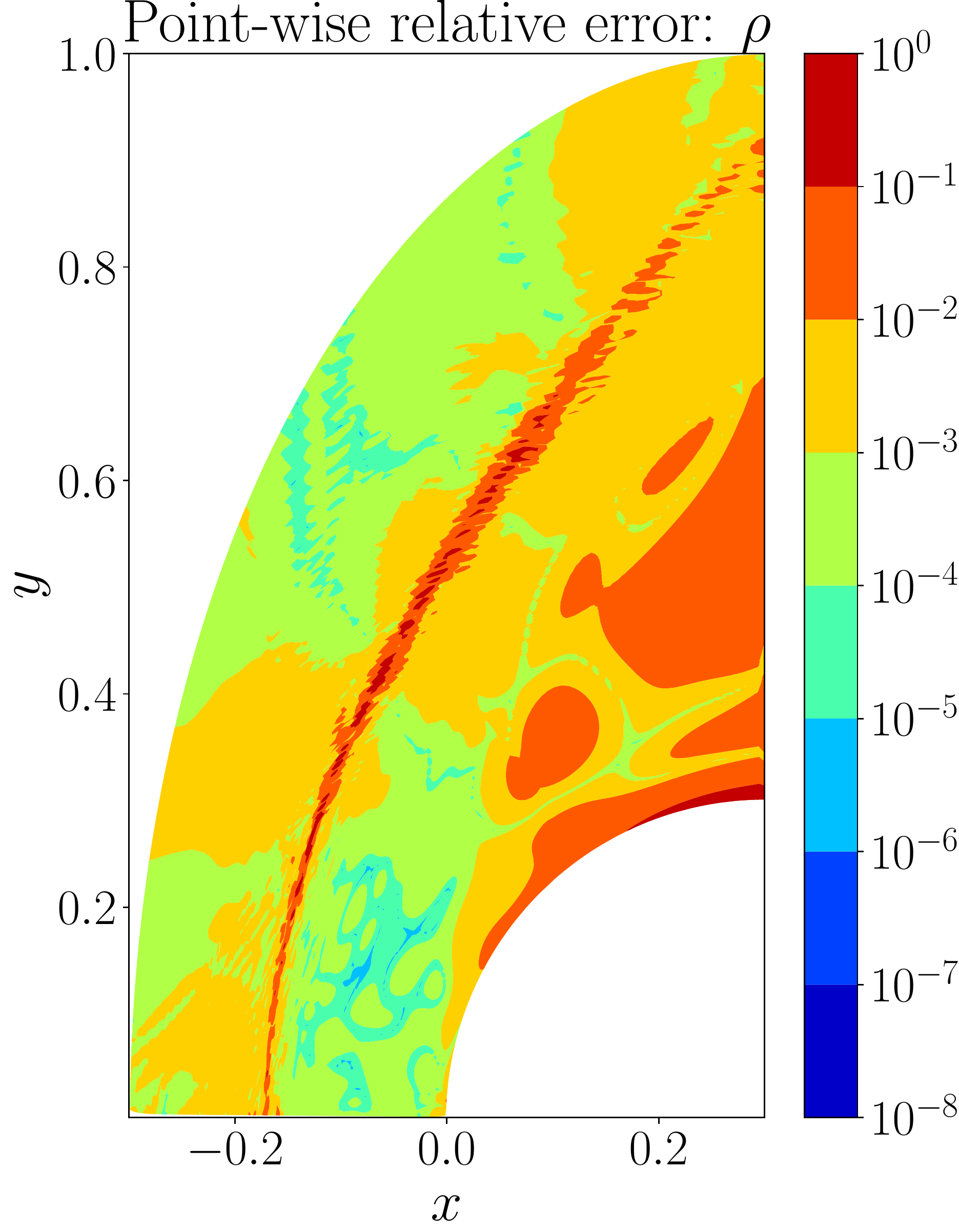}
    }
    \subfigure{
    \includegraphics[trim=0cm 0cm 0cm 0cm, clip=true, scale=0.16, angle = 0]{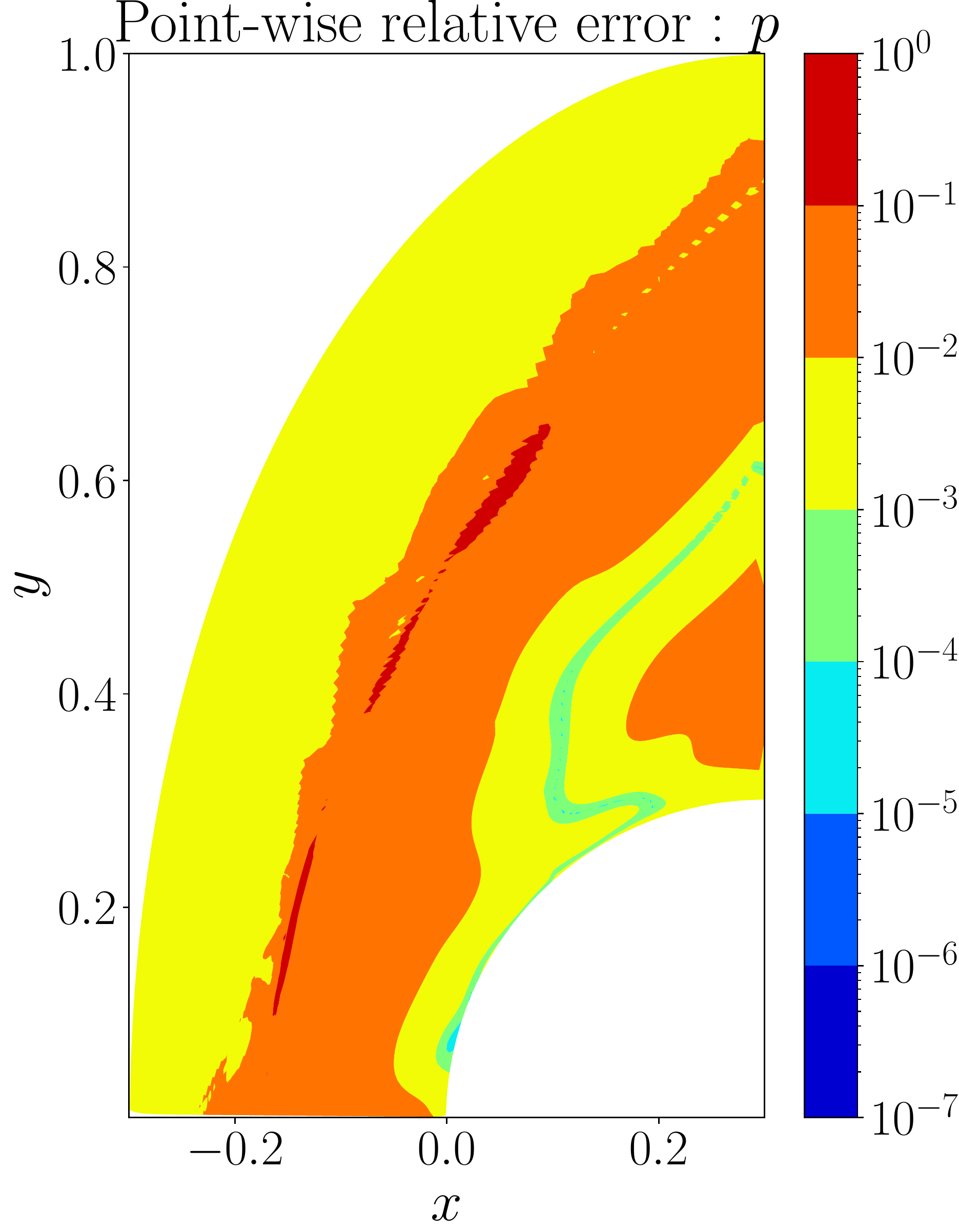}
    }
    \subfigure{
    \includegraphics[trim=0cm 0cm 0cm 0cm, clip=true, scale=0.16, angle = 0]{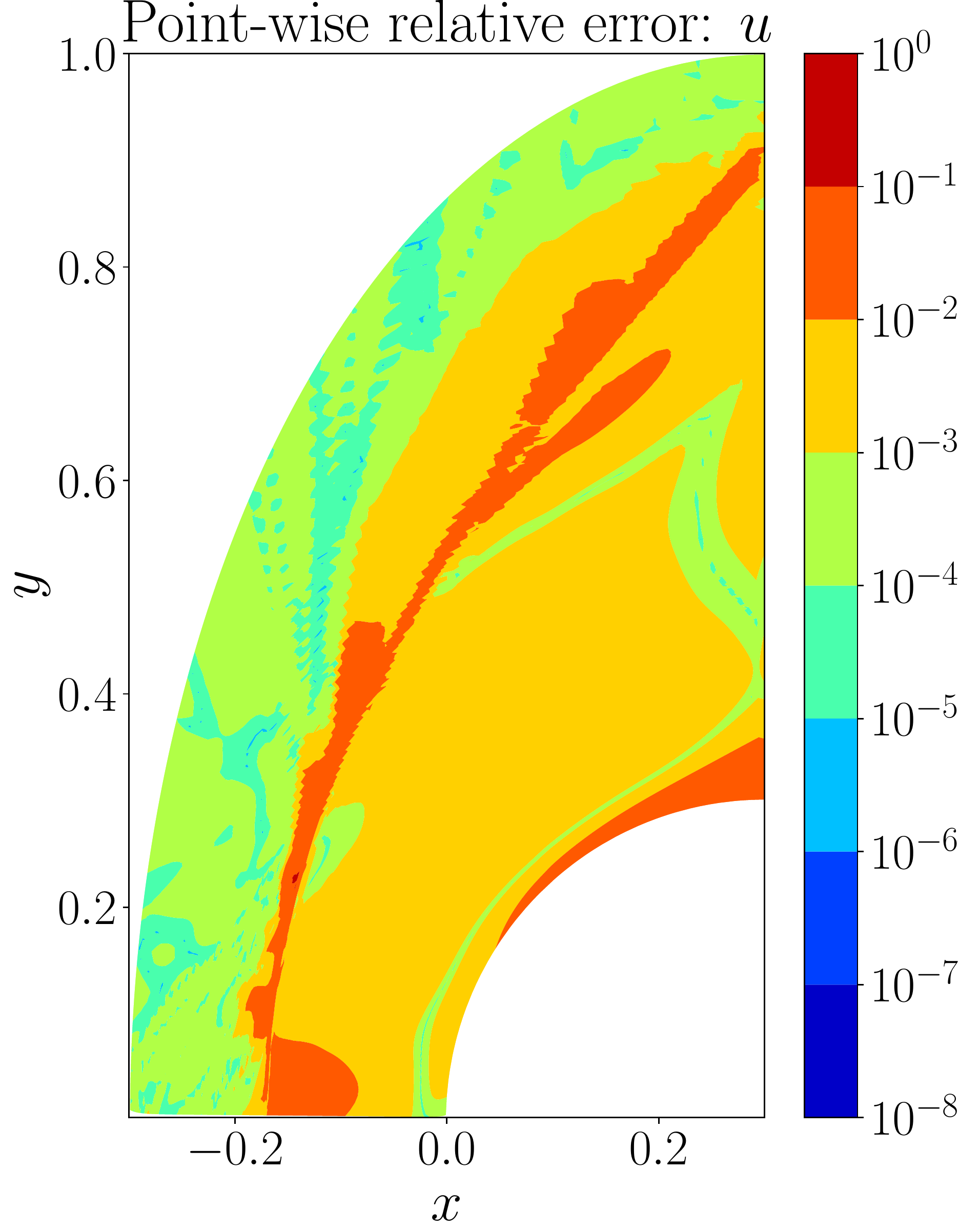}
    }
    \subfigure{
    \includegraphics[trim=0cm 0cm 0cm 0cm, clip=true, scale=0.16, angle = 0]{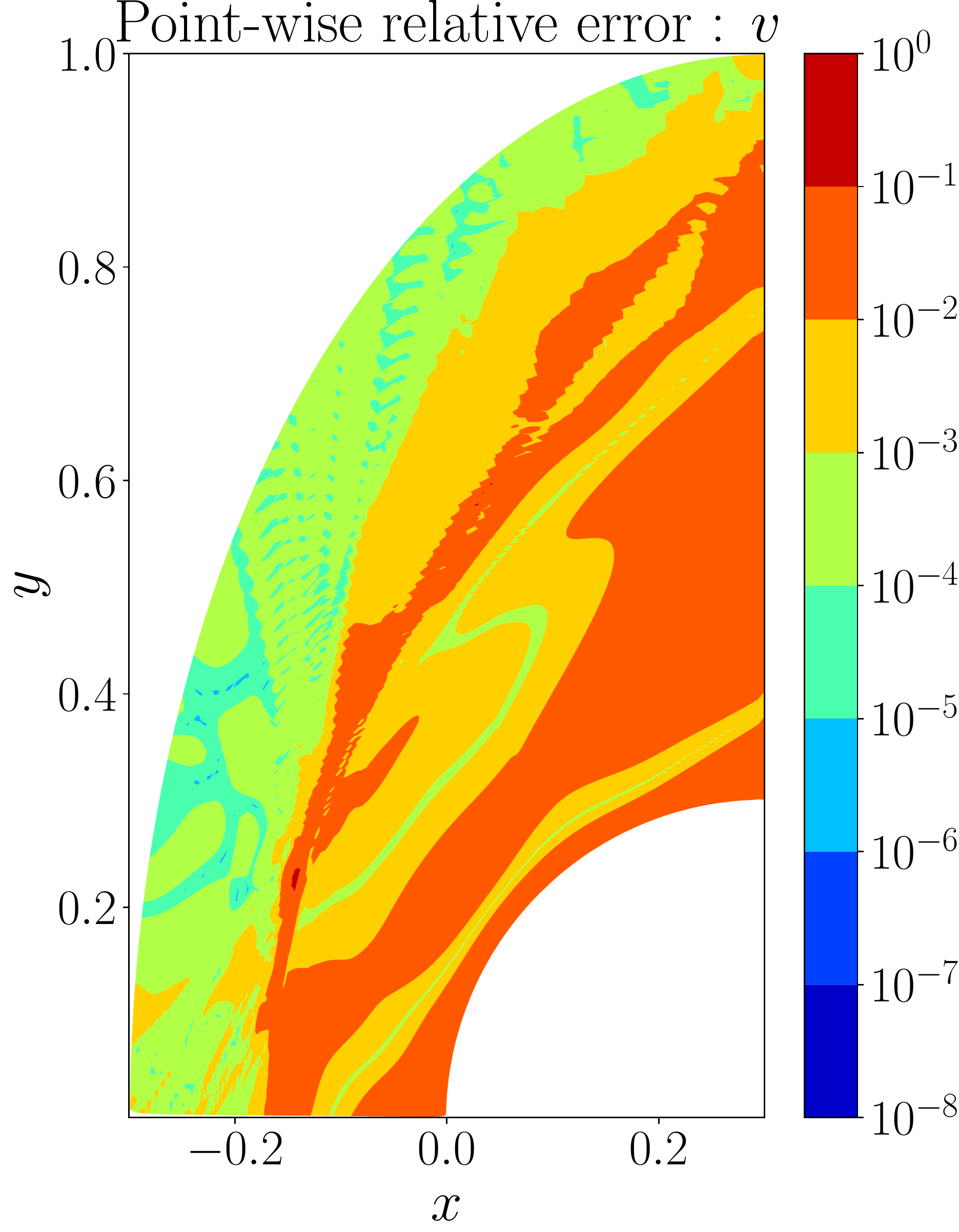}
    }\\
    \subfigure{
    \includegraphics[trim=0cm 0cm 0cm 0cm, clip=true, scale=0.16, angle = 0]{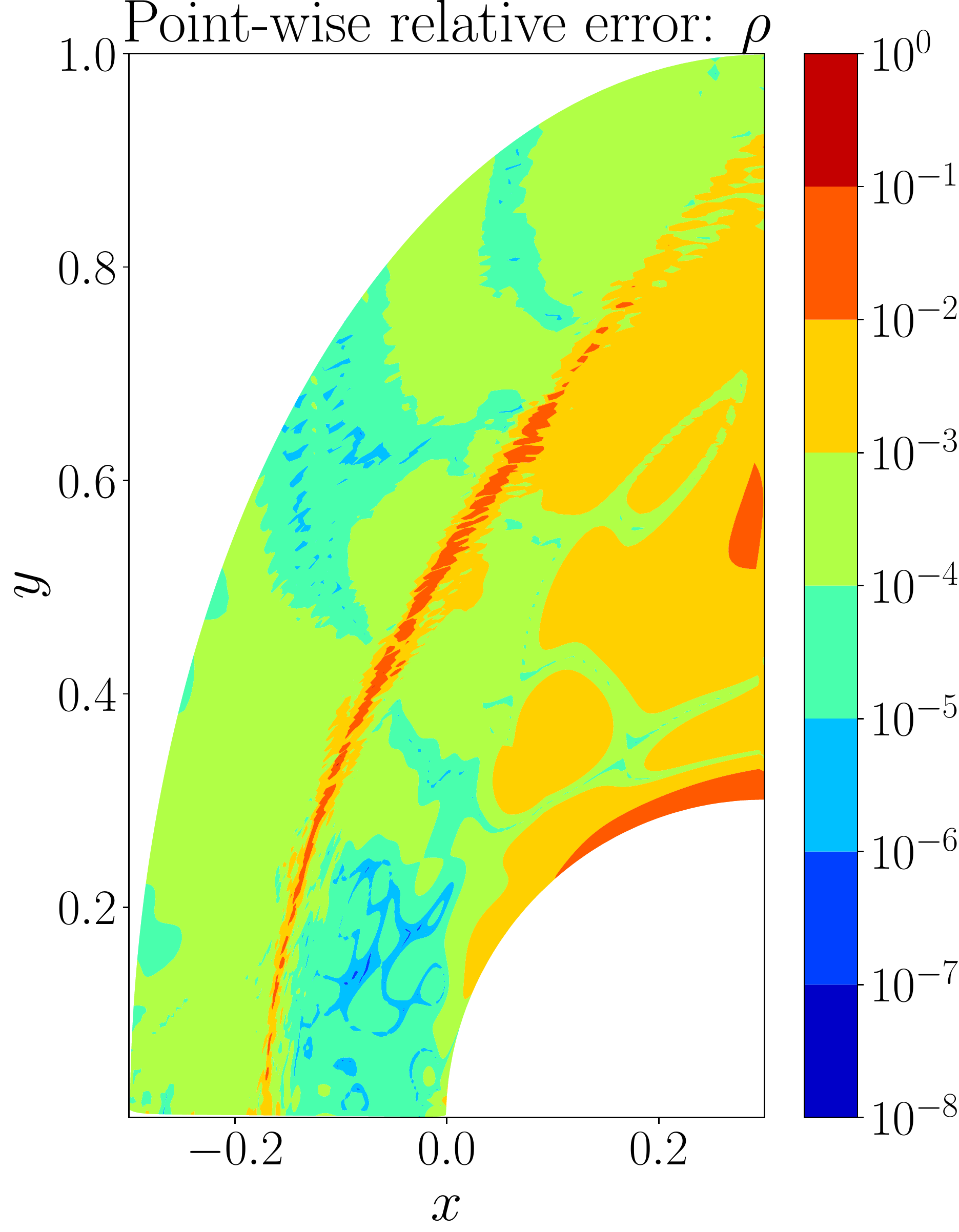}
    }
    \subfigure{
    \includegraphics[trim=0cm 0cm 0cm 0cm, clip=true, scale=0.16, angle = 0]{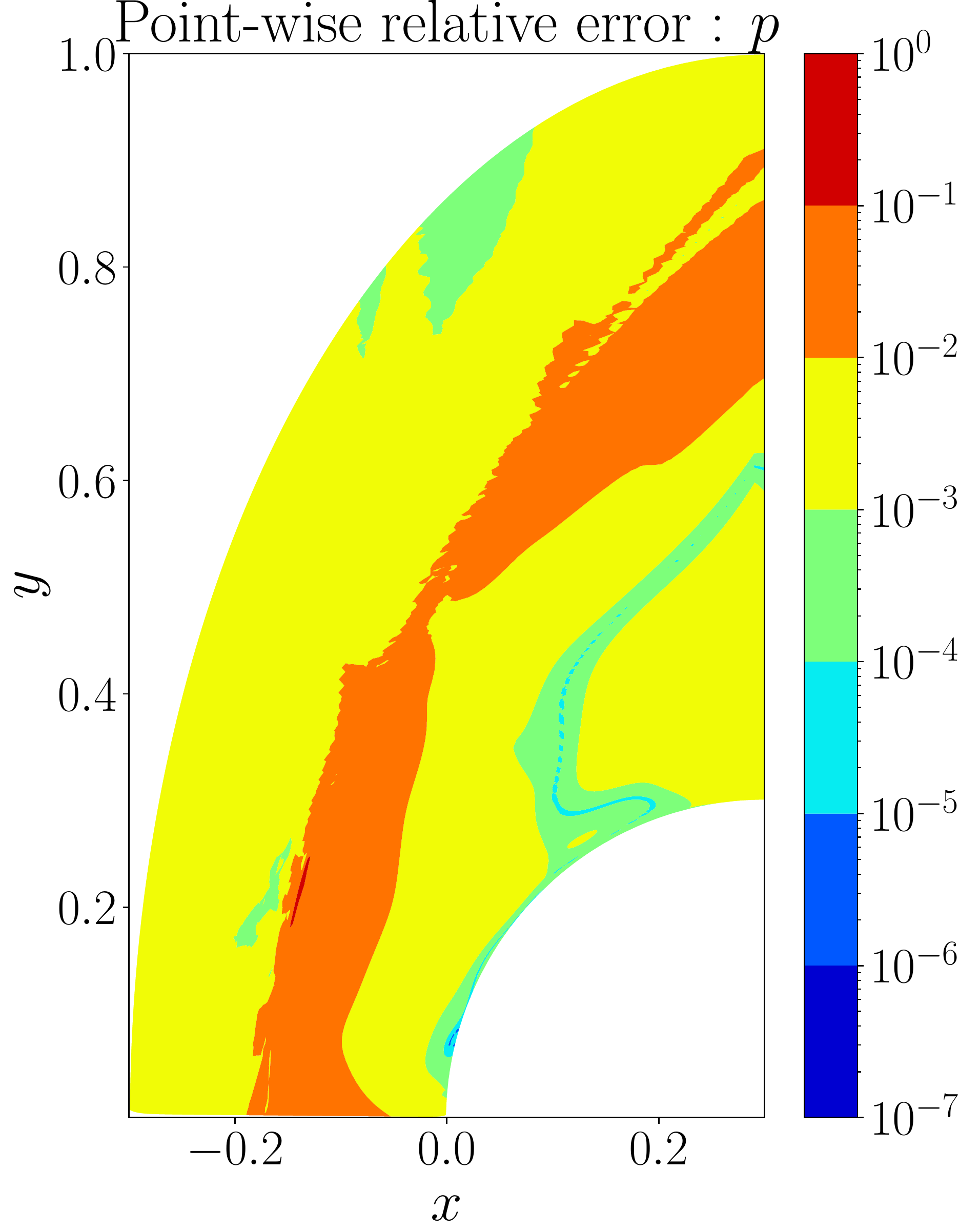}
    }
    \subfigure{
    \includegraphics[trim=0cm 0cm 0cm 0cm, clip=true, scale=0.16, angle = 0]{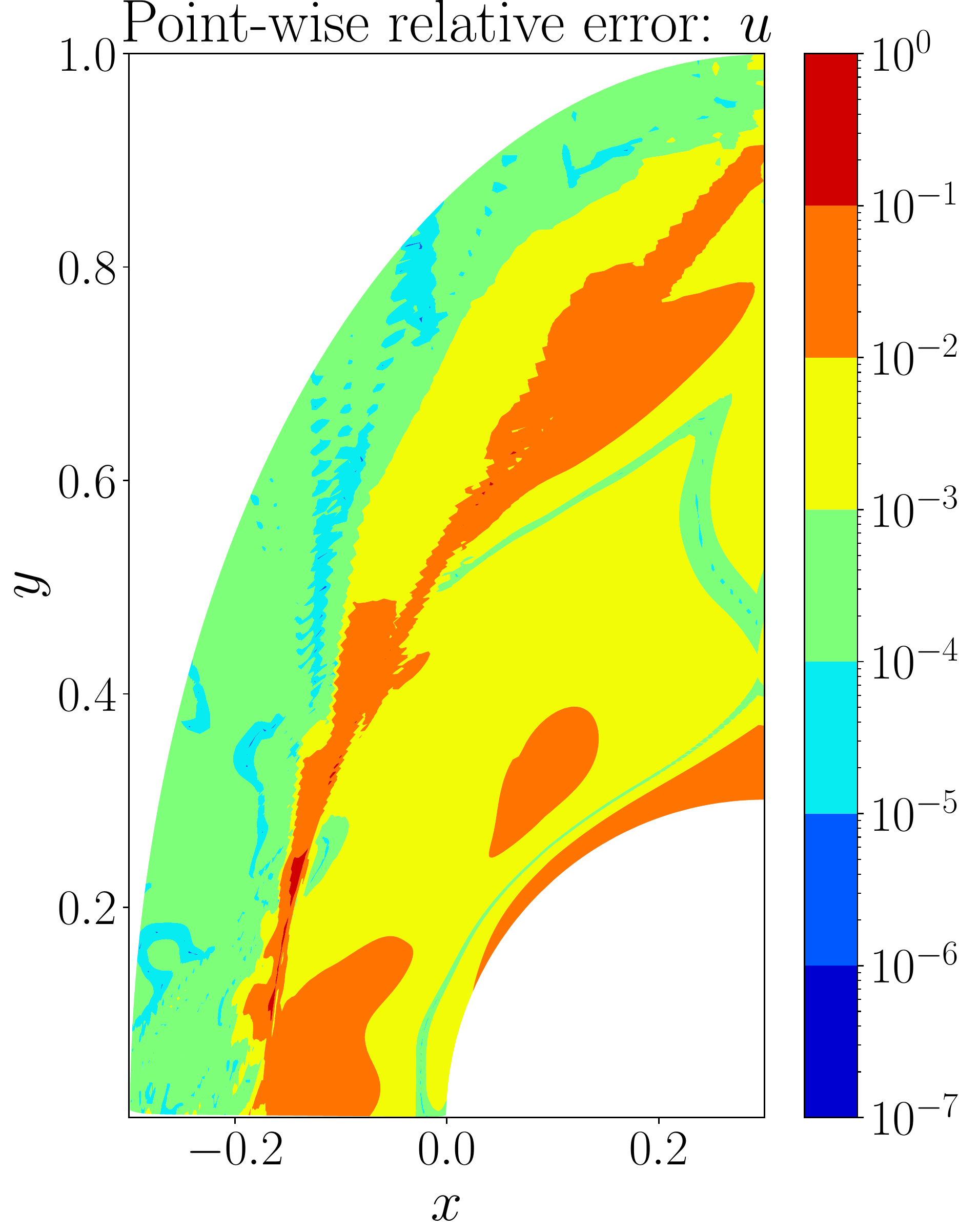}
    }
    \subfigure{
    \includegraphics[trim=0cm 0cm 0cm 0cm, clip=true, scale=0.16, angle = 0]{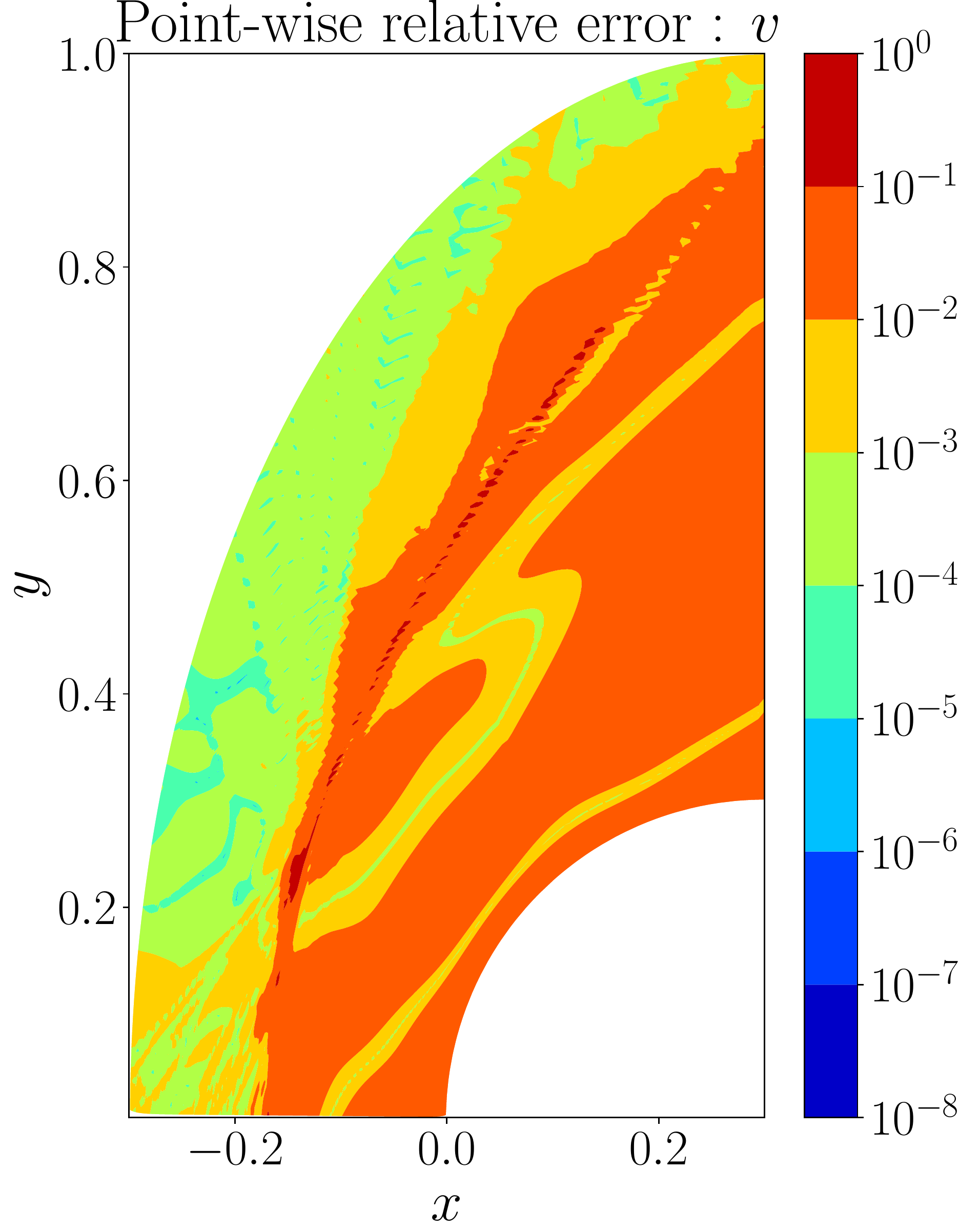}
    }
    \caption{Bow shock problem: Point-wise errors for PINNs (top row), and XPINNs (bottom row). The point-wise error in the density, pressure and velocity components are shown from left to right columns, respectively. }
    \label{fig:distribution:bow:pred}
\end{figure}
\begin{figure}[http]
    \centering
    \includegraphics[trim=0cm 0cm 0cm 0cm, clip=true, scale=0.4, angle = 0]{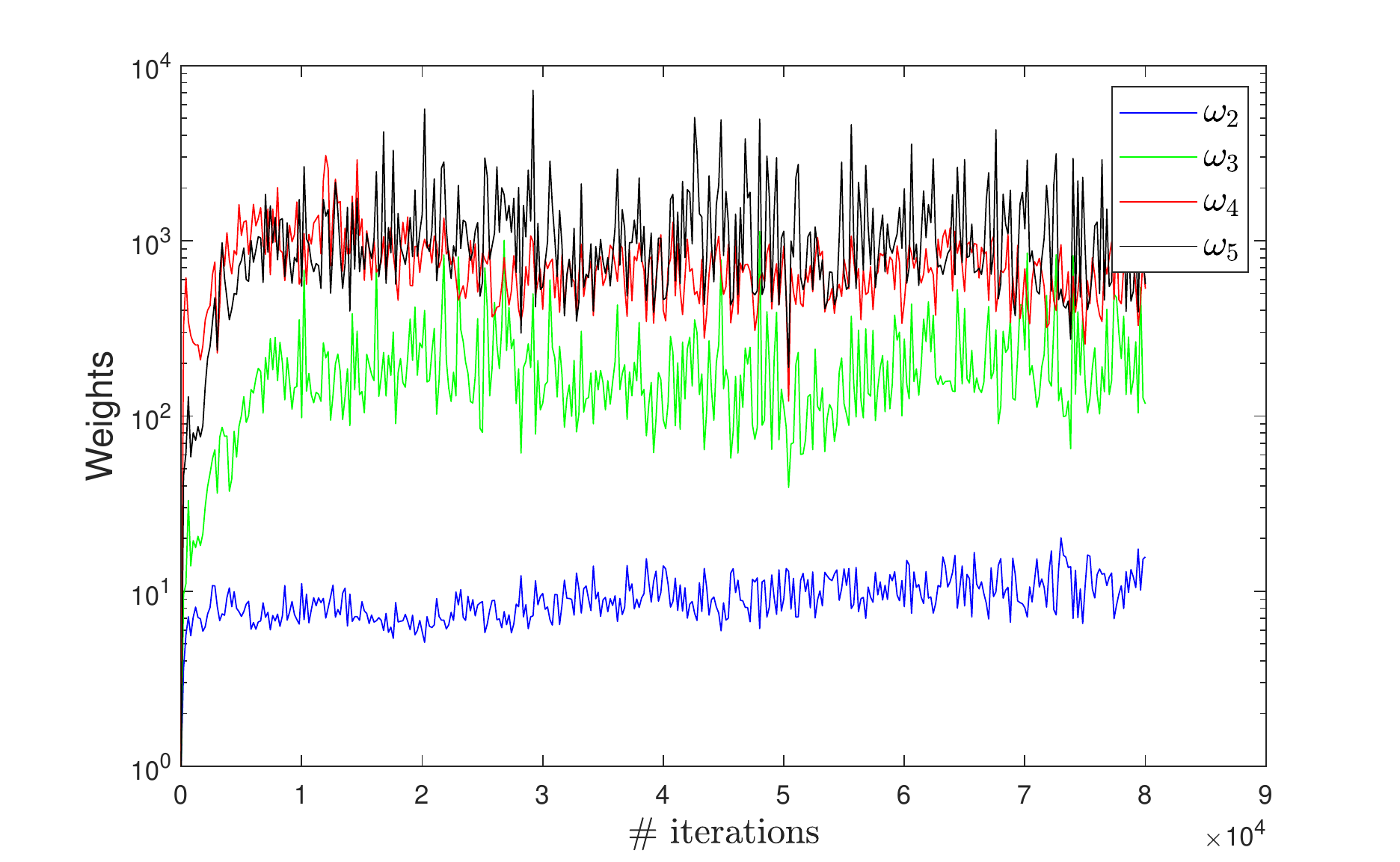}
    \caption{Bow shock problem: Variation of dynamic weights for PINNs.}
    \label{fig:BowWeights}
\end{figure}
\begin{figure}[http]
    \centering
        \subfigure{
    \includegraphics[trim=0cm 0cm 0cm 0cm, clip=true, scale=0.45, angle = 0]{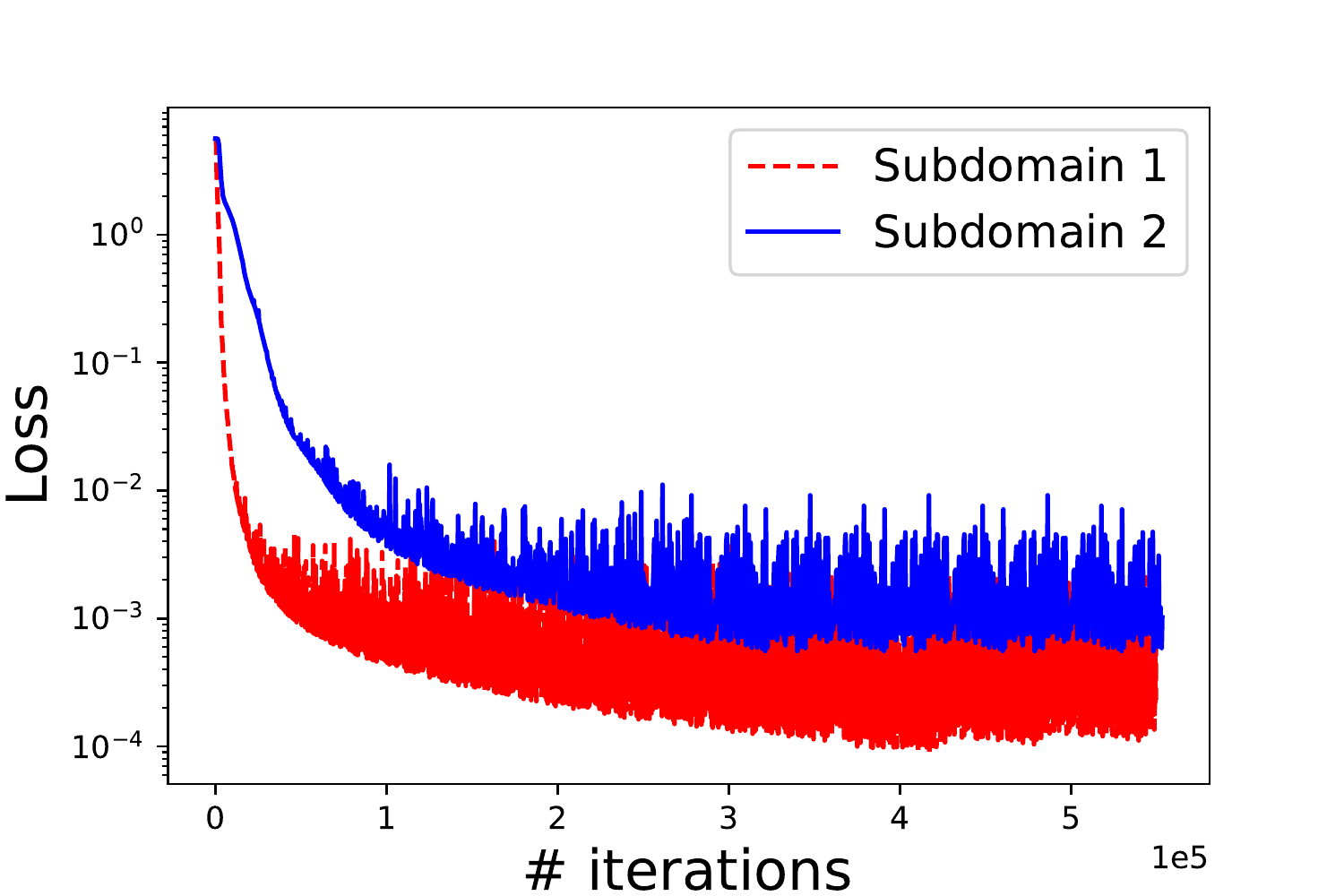}
    }
    \subfigure{
    \includegraphics[trim=0cm 0cm 0cm 0cm, clip=true, scale=0.45, angle = 0]{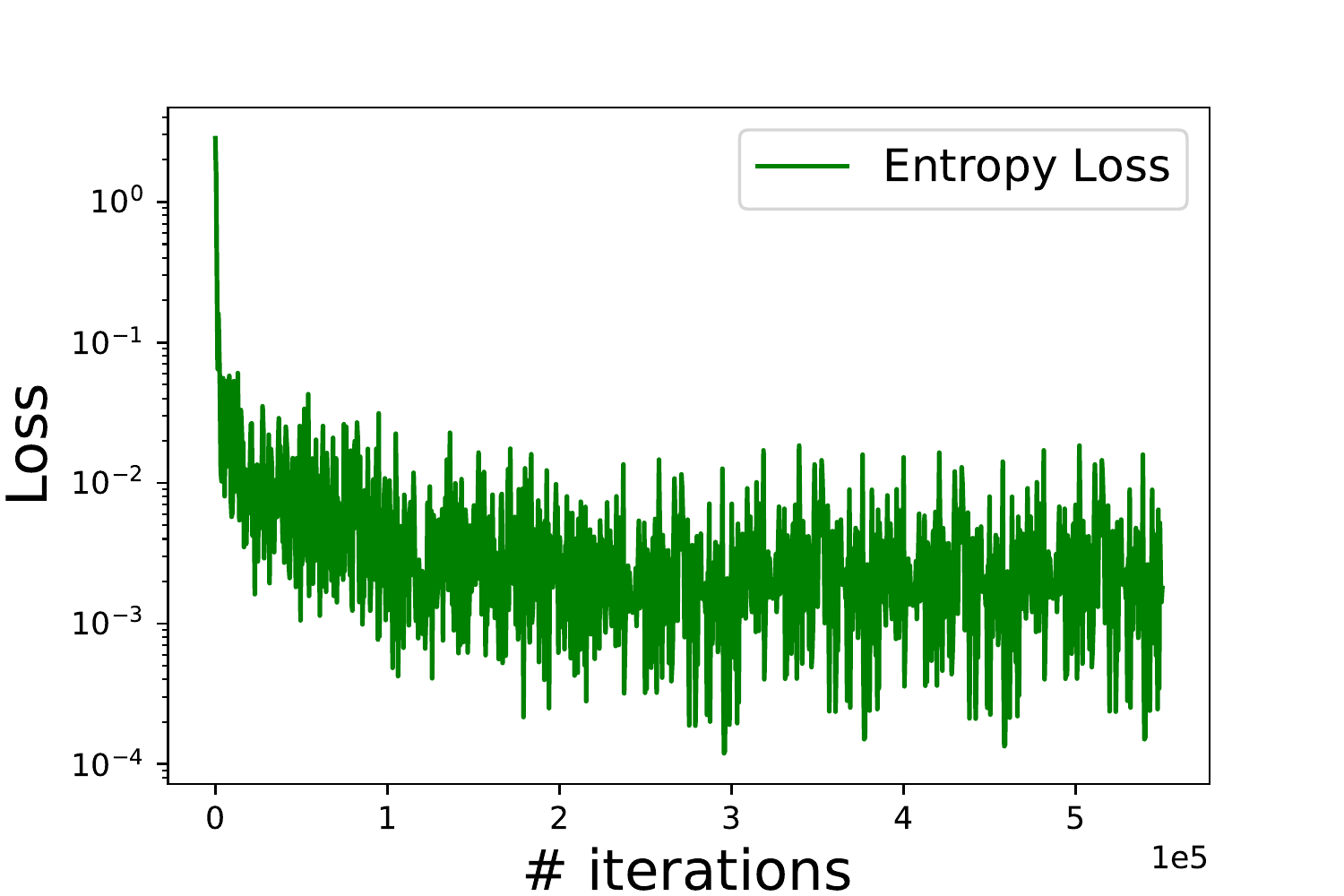}
    }
    \caption{Bow shock problem: Loss function variation for XPINNs in two subdomains (left) and the combined entropy loss in the entire computational domain (right).}
    \label{fig:distribution:bow:loss}
\end{figure}

We present in Fig. \ref{fig:distribution:bow:exact} the reference solutions for all primitive variables, and in Figure \ref{fig:distribution:bow:pred} the point-wise errors. The point-wise error in the density, pressure and velocity components are shown from left to right columns. It can be observed that the predictive accuracy of XPINN is better than that of PINNs. The dynamic weights for PINNs are shown in Fig. \ref{fig:BowWeights}. The dynamic weights $\omega_4$ and $\omega_5$ converge around $10^3$, whereas the weights $\omega_2$ and $\omega_3$ converge around $10^1$ and $10^2$. The loss functions versus number of iterations are plotted in Figure \ref{fig:distribution:bow:loss} (left), whereas the right side figure shows the combined entropy loss in the entire domain.

\section{Summary}
 In this work we solved the inverse supersonic flow problems involving expansion and compression (shock) waves. Accurate solutions to such inverse problems are often required for designing aerospace vehicles, which are otherwise difficult or even impossible to solve using the
traditional numerical solvers. For this purpose, we employ physics-informed neural networks (PINNs) as well as it's extended space-time domain decomposition-based version, named as extended PINNs (XPINNs). For the computation, we used the density gradient data obtained from the Schlieren technique, inflow boundary data and the pressure wall boundary data. Apart from the governing compressible Euler equations, the entropy conditions satisfied by the more relevant viscosity solutions are also enforced. Moreover, we also enforced the positivity
conditions on density and pressure.  We compared the results of PINNs and XPINNs in terms of predictive accuracy. We used the theoretical generalization bounds in XPINN and PINN models, which show the complexity of PINN and XPINNs methods as well as how well these methods generalize, particularly for the expansion wave and the oblique shock wave test cases.
We employed adaptive activation functions for both PINNs and XPINNs, which showed good predictive accuracy compared to fixed activation functions. Moreover, the use of dynamic weights shows accuracy improvements over fixed weights. For the problems with discontinuities, additional techniques, such as adaptive activate functions and dynamical weights, are important and can significantly improve the accuracy of the PINN solutions. XPINNs with locally tuned neural networks also gives more accurate results suggesting that XPINN is also a promising and scalable method for tackling problems with complex solutions involving local discontinuous features such as shocks.

\appendix

\section{Pedagogical example: Inverse problem for the two-dimensional Euler equations with smooth solutions}\label{sec:2dsmooth}
\begin{figure}[htpb]
    \centering
    \subfigure[Data for $\nabla \rho$]{\label{y=x}
    \includegraphics[trim=0cm 0cm 0cm 0cm, clip=true, scale=0.17, angle = 0]{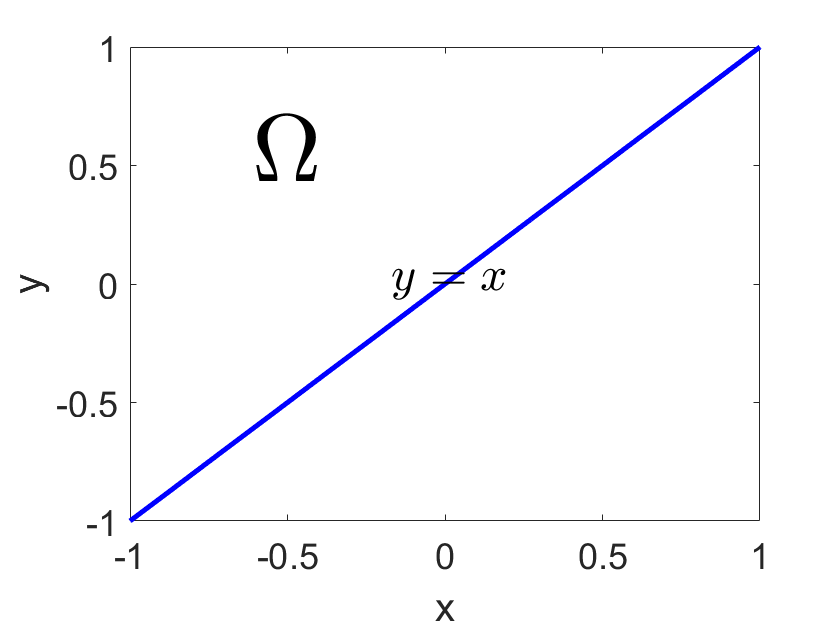}}
    \subfigure[Exact solution: density]{\label{smooth:ex:rho}
    \includegraphics[trim=0cm 0cm 0cm 0cm, clip=true, scale=0.17, angle = 0]{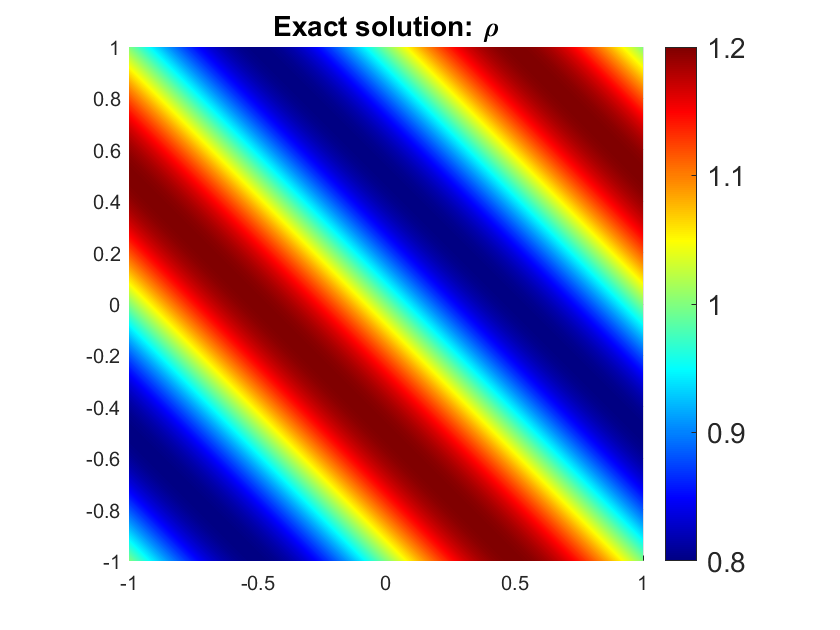}}
    \subfigure[PINN solution: density]{\label{smooth:pinn:rho}
    \includegraphics[trim=0cm 0cm 0cm 0cm, clip=true, scale=0.17, angle = 0]{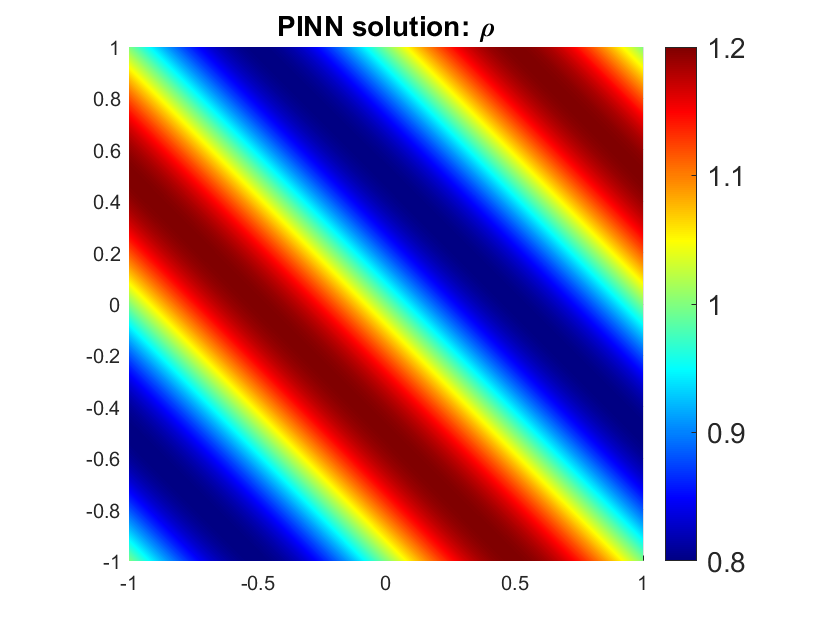}}\\
    \subfigure[PINN solution: pressure]{\label{smooth:pinn:p}
    \includegraphics[trim=0cm 0cm 0cm 0cm, clip=true, scale=0.17, angle = 0]{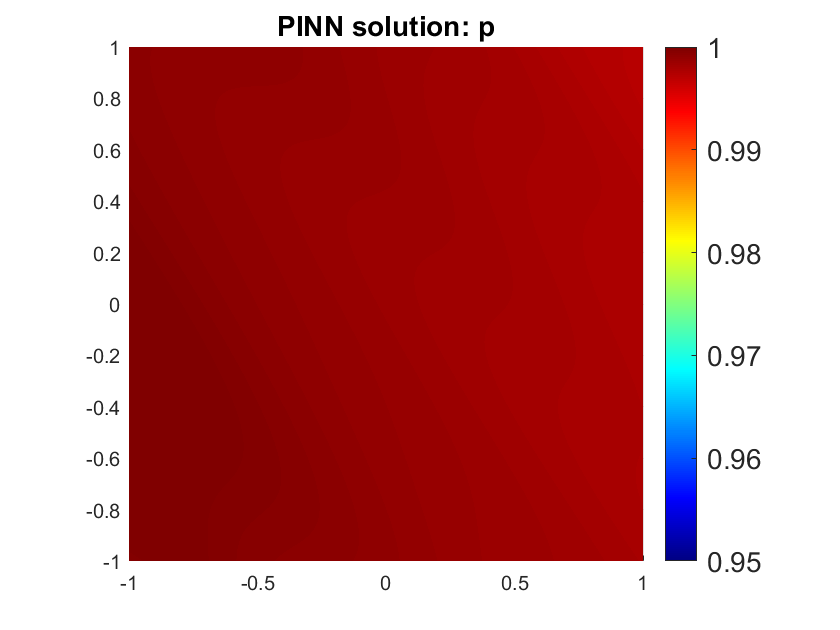}}
    \subfigure[PINN solution: velocity $u$]{\label{smooth:pinn:u}
    \includegraphics[trim=0cm 0cm 0cm 0cm, clip=true, scale=0.17, angle = 0]{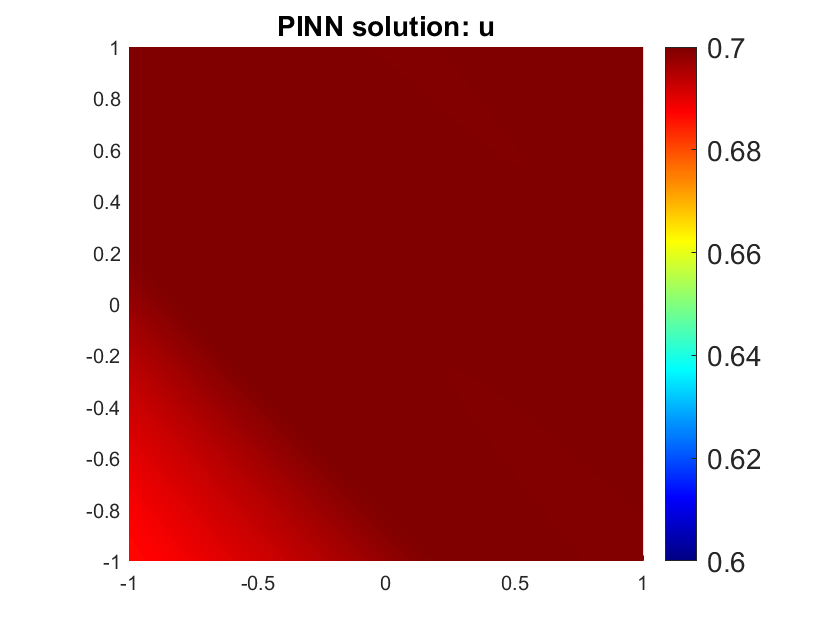}}
    \subfigure[PINN solution: velocity $v$]{\label{smooth:pinn:v}
    \includegraphics[trim=0cm 0cm 0cm 0cm, clip=true, scale=0.17, angle = 0]{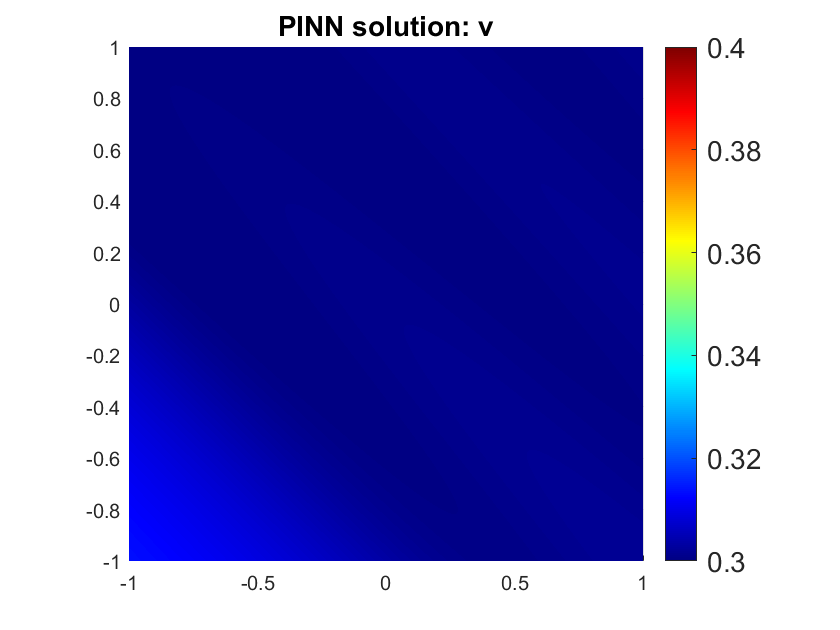}}
    \caption{(a) We only use the data of density gradient along the line $y=x$ shown in the plot. (b) Exact solution of the density. (c) PINN solution of the density. (d) PINN solution of the pressure. The exact solution of the pressure is $p(x,y) = 1$. (e) PINN solution of the velocity $u$. The exact solution of the velocity $u$ is $u(x,y) = 0.7$. (f) PINN solution of the velocity $v$. The exact solution of the velocity $v$ is $v(x,y) = 0.3$.}
    \label{fig:smooth}
\end{figure}
\emph{Consider the Euler equation \eqref{NCL} in the {two-dimensional} case with the periodic boundary conditions
and initial conditions
\begin{equation*}\label{IC:Smooth}
    U_0 = (\rho_0,~~ u_0,~~ v_0,~~ p_0) = \left(1.0 + 0.2\sin(\pi (x+y)),~~ 0.7, ~~0.3, ~~1.0 \right),
\end{equation*}
in which case we have the exact solutions
\begin{equation*}\label{Exact:Smooth}
    (\rho,~~ u,~~ p) = \left(1.0 + 0.2\sin(\pi (x +y -(u+v)t) ),~~0.7, ~~ 0.3, ~~1.0 \right),
\end{equation*}
where $(x,y)\in (-1,1)^2$.
}

In the one-dimensional case, we choose $D$ to be the whole computational domain. However, in the two-dimensional case, we use  density gradient only in a subdomain. In this example, we use the density gradient along the line $y = x$, see Figure \ref{y=x}. For the pressure, we only use one single point value of the pressure at different time, here we use $p(-1,-1, t)$. 
We also use the inflow conditions.
Therefore, the loss function is given by:
\begin{equation}
    \mathcal{J} = \text{MSE}_{\mathcal{F}} + \text{MSE}_{\nabla \rho|_{y=x}} + \text{MSE}_{p^*} + \text{MSE}_{Inflow}.
\end{equation}
The results are shown in Figs. \ref{fig:smooth} (c)-(f). Observe that we now have good predictions for all the states, i.e., the density, velocity and the pressure. 

\vspace{-0.2cm}
\section*{Acknowledgements}
This work was supported by the OSD/AFOSR MURI grant FA9550-
20-1-0358 and the Alexander von Humboldt fellowship to G.E. Karniadakis. The work of
Z. Mao was conducted while he was a postdoc at Brown University and also during his one month-visit to 
the Technical University of Munich working with Prof. N. Adams. 

\vspace{-0.2cm}

\bibliographystyle{unsrt}
\bibliography{ref}
\end{document}